\documentclass[a4paper,abstracton]{scrartcl}
\usepackage{amsfonts}
\usepackage{amssymb}
\usepackage{amsmath}
\usepackage{amsthm}

\usepackage{mathtools}
\usepackage{multirow}

\usepackage[ruled,vlined,linesnumbered]{algorithm2e}

\usepackage{graphicx}
\usepackage{color}

\usepackage{subfigure}

\usepackage{a4wide}

\usepackage{authblk}

\allowdisplaybreaks

\begin{document}

\newcommand{\X}{{\mathcal{X}}}
\newcommand{\cU}{{\mathcal{U}}}
\newcommand{\cI}{{\mathcal{I}}}
\newcommand{\cC}{{\mathcal{C}}}
\newcommand{\cS}{{\mathcal{S}}}

\newcommand{\MMD}[1]{\texttt{MM-D-{#1}}}
\newcommand{\MMB}[1]{\texttt{MM-B-{#1}}}
\newcommand{\MMRI}[1]{\texttt{MMR-I-{#1}}}
\newcommand{\MMRD}[1]{\texttt{MMR-D-{#1}}}
\newcommand{\TSTD}[1]{\texttt{2ST-D-{#1}}}
\newcommand{\TSTDB}[1]{\texttt{2ST-DB-{#1}}}
\newcommand{\TSTCB}[1]{\texttt{2ST-CB-{#1}}}
\newcommand{\RRD}[1]{\texttt{RR-D-{#1}}}
\newcommand{\RRDB}[1]{\texttt{RR-DB-{#1}}}
\newcommand{\RRCB}[1]{\texttt{RR-CB-{#1}}}
	
\title{Benchmarking Problems for Robust Discrete Optimization}
	
\author{Marc Goerigk}
\author{Mohammad Khosravi\thanks{Corresponding author. Email: mohammad.khosravi@uni-siegen.de}}
	
\affil{Network and Data Science Management, University of Siegen,\\Unteres Schlo{\ss}~3, 57072 Siegen, Germany}

\date{}
	
\maketitle

\begin{abstract}
Robust discrete optimization is a highly active field of research where a plenitude of combinations between decision criteria, uncertainty sets and underlying nominal problems are considered. Usually, a robust problem becomes harder to solve than its nominal counterpart, even if it remains in the same complexity class. For this reason, specialized solution algorithms have been developed. To further drive the development of stronger solution algorithms and to facilitate the comparison between methods, a set of benchmark instances is necessary but so far missing. In this paper we propose a further step towards this goal by proposing several instance generation procedures for combinations of min-max, min-max regret, two-stage and recoverable robustness with interval, discrete or budgeted uncertainty sets. Besides sampling methods that go beyond the simple uniform sampling method that is the de-facto standard to produce instances, also optimization models to construct hard instances are considered. Using a selection problem for the nominal ground problem, we are able to generate instances that are several orders of magnitudes harder to solve than uniformly sampled instances when solving them with a general mixed-integer programming solver. All instances and generator codes are made available online.
\end{abstract}
	

\noindent\textbf{Acknowledgements:} Supported by the Deutsche Forschungsgemeinschaft (DFG) through grant GO 2069/1-1.

\section{Motivation}
\label{sec:motivation}
	
In real-world decision making, problem data and the consequences of our actions can rarely be predicted with full accuracy. Instead, we must make use of imprecise measurements, forecasts and estimates to be capable of acting in an uncertain environment. Several methods have been developed to treat uncertainty in decision making, such as stochastic programming \cite{birge2011introduction} or fuzzy optimization \cite{klir1995fuzzy}. This paper focuses on the robust optimization paradigm for discrete optimization problems \cite{kouvelis2013robust,ben2009robust,kasperski2016robust}. While this paradigm can be realized in many different variants, the basic unifying idea is to define an uncertainty set of possible parameter outcomes, and to protect against the worst case over this set.
	
Robust optimization problems for discrete optimization problems have been studied for decades under various angles. While many publications focus on the complexity side of these problems (e.g., determining which variants can be solved in polynomial time, and which become hard), an active research community has emerged that is concerned with solving such problem in the most efficient manner \cite{goerigk2016algorithm}. Frequently, such algorithms are developed for a very specific type of problem. The diversity of robust optimization variants in combination with the trove of classic or real-world combinatorial optimization problems with which they can be combined have created a partitioned rather than unified landscape of methods. This has become visible by a lack of standardized benchmarks, problem generators, or algorithm libraries. As explained in the algorithm engineering paradigm  \cite{sanders2009algorithm,muller2001algorithm}, having a strong foundation of shared instances with easy access is instrumental in the development of better algorithms. Currently, many papers need to reimplement generator methods. At the same time, these generators usually show little sophistication, but instead sample parameter values uniformly from an interval.
	
The main goal of this paper is to provide a starting point in this direction, by establishing a first collection of instances that shall be further extended in the future. To this end, there are several obstacles that need to be overcome. Firstly, the large variety of types of robust optimization methods in connection with the possibilities of underlying nominal combinatorial problems means that we cannot consider all combinations simultaneously. As a response, one could fix on type of robust optimization approach (e.g., one-stage min-max with discrete uncertainty) and vary the underlying nominal problem. What we present in this paper instead is to focus on one fixed nominal problem and to vary possible ways to create a robust counterpart. This underlying problem is the so-called selection problem (see, e.g., \cite{deineko2013complexity}), where the only constraint states that a subset of exactly $p$ items have to be chosen from $n$ available items. The advantage of this problem is that it provides as little structure on its own as possible, thus disentangling the difficulty that may come from a complicated structure of feasible solutions (e.g., the structure of a graph for a shortest path problem) with the difficulty that arises from the robust objective function. Furthermore, it allows us to focus on uncertainty in the objective function.
	
Secondly, a benchmark set should contain challenging problems that drive the development of new and improved solution methods. But to determine the computational challenge of an instance, we need to know which algorithm is used to solve it. This difficulty is avoided by using what we believe to be the most common universal solution approach, that is, by formulating a compact mixed-integer program (MIP) and solving this through a general-purpose MIP solver.
	
Our contributions are as follows. We consider min-max, min-max regret, two-stage and recoverable robust selection problems in combination with discrete, interval and budgeted uncertainty sets. For each combination we consider different instance generators. We always include the simple approach of choosing parameter values uniformly from an interval as the baseline heuristic. Additionally, we consider alternative ways to sample problem instances, and compare the resulting computation times. We show that this way it is possible to create instances that are several orders of magnitude harder to solve. Furthermore, we build upon the HIRO (''hard instances for robust optimization'') approach from \cite{goerigk2020generating} to formulate the generation of hard instances as optimization problems. While such an approach is clearly more time-intensive to use than a generator based on random sampling, we show that for most problem variants, even harder instances can be found. All instances along with generator code are provided in an easy-to-use format online, in an attempt to stimulate further research for better solution methods.
	
The remainder of this paper is structured as follows. In Section~\ref{sec:problems}, we formally introduce all types of robust optimization problems and uncertainty sets that we use. Additionally, we recall the high-level idea of the HIRO approach from \cite{goerigk2020generating}. We then present methods for instance generation along with computational results for various problem variants: Min-max problems are presented in Section~\ref{sec:minmax}, min-max regret problems in Section~\ref{sec:minmaxregret}, two-stage problems in Section~\ref{sec:twostage}, and recoverable problems in Section~\ref{sec:recoverable}. We conclude this work by discussing the structure and format of our benchmark library in Section~\ref{sec:library} before summarizing our findings and pointing out further research questions in Section~\ref{sec:conclusions}. In Appendix~\ref{app:overview}, we give an overview to all sampling methods that are introduced in this paper, while Appendix~\ref{app:times} shows detailed generation times.

\section{Problem Definitions}
\label{sec:problems}
	
We write vectors in bold and use the notation $[n]$ to denote sets $\{1,\ldots,n\}$.
We focus on combinatorial problems with uncertain objective function. That is, given a set of feasible solutions $\X\subseteq \{0,1\}^n$, the so-called nominal problem is to solve
\[ \min_{\pmb{x}\in\X} \pmb{c}^t \pmb{x} \]
for a known cost vector $\pmb{c}\in\mathbb{R}^n_+$. In the selection problem that we use as a baseline nominal problem, we have $\X = \{\pmb{x}\in\{0,1\}^n : \sum_{i\in[n]} x_i = p\}$. It is easy to check that the optimal solution of the selection problem is a set of $p$ items with the smallest costs, which can be found in $O(n)$.  Due to their simplicity, selection problems have been playing a key role in the complexity analysis for robust combinatorial optimization, see, e.g., the results presented in \cite{busing2011phd,
dolgui2012min,
deineko2013complexity,
kasperski2015approximability,
kasperski2017robust,
chassein2018recoverable}.
	
If costs are uncertain, the robust optimization approach is to define an uncertainty set $\cU\subseteq \mathbb{R}^n_+$ containing all possible cost vectors that we wish to protect against. A range of such concepts have been defined; the most frequently used are listed here (for a more detailed overview, see, e.g., \cite{kasperski2016robust}).
	
In the min-max setting, we need to find a single solution $\pmb{x}\in\X$ that performs optimally under the worst-case decision criterion, that is, we need to solve
\begin{equation}
\min_{\pmb{x}\in\X} \max_{\pmb{c}\in\cU} \pmb{c}^t \pmb{x} \tag{Min-Max}
\end{equation}
In the min-max regret criterion, we compare against the best possible outcome under each scenario. Let $opt(\pmb{c}) = \min_{\pmb{x}\in\X} \pmb{c}^t \pmb{x}$ be the objective value if the cost scenario $\pmb{c}$ is known. The problem is then to solve
\begin{equation}
\min_{\pmb{x}\in\X} \max_{\pmb{c}\in\cU} \left( \pmb{c}^t \pmb{x} - opt(\pmb{c}) \right) \tag{Min-Max Regret}
\end{equation}
In case that the decision maker can respond to a scenario after it has been revealed, other approaches can be applied, which slightly differ in philosophy. For two-stage robust problems, the decision maker only fixes part of the solution ahead of the uncertainty. After the scenario has been revealed, he can complete the partial solution to a full solution. In recoverable robustness (which is essentially also a two-stage problem), the decision maker first fixes a complete solution, but is allowed to modify this solution after the scenario is known. More formally, a common definition is
\begin{equation}
\min_{\pmb{x}\in\X'} \max_{\pmb{c}\in\cU} \min_{\pmb{y}\in\X(\pmb{x})} \pmb{C}^t \pmb{x} + \pmb{c}^t \pmb{y} \tag{2St}
\end{equation}
for two-stage robustness with $\X(\pmb{x}) = \{\pmb{y}\in\{0,1\}^n : \pmb{x}+\pmb{y}\in\X\}$ and $\X' = \{\pmb{x}\in\{0,1\}^n : \X(\pmb{x})\neq \emptyset\}$, and
\begin{equation}
\min_{\pmb{x}\in\X} \max_{\pmb{c}\in\cU} \min_{\pmb{y}\in\X : d(\pmb{x},\pmb{y})\le\Delta} \pmb{C}^t \pmb{x} + \pmb{c}^t \pmb{y} \tag{Rec}
\end{equation}
for recoverable robustness, where $d: \X\times\X \to \mathbb{R}$ denotes some function measuring dissimilarity between solutions, and $\Delta$ is the maximum recovery distance.
	
Many more problem variants have been considered as well, which are not covered in this paper; these include for example the ordered weighted averaging criterion \cite{yager2012ordered,chassein2020approximating} or two-stage regret \cite{jiang2013two,goerigk2020combinatorial}.
	
In all decision problems introduced so far, we require a set $\cU$ of possible scenarios. The structure of this set is a deciding factor in the theoretical and practical hardness of instances. Amongst the most common sets for combinatorial problems are discrete uncertainty sets
\[ \cU_D = \{\pmb{c}^1, \ldots, \pmb{c}^N \} \]
where $N$ scenarios are simply given as an explicit list; interval sets
\[ \cU_I = \{\pmb{c}\in\mathbb{R}^n_+ : c_i \in [\underline{c}_i, \underline{c}_i + d_i] \} \]
where for each item $i\in [n]$ a lower bound $\underline{c}_i$ and a deviation $d_i$ are given; and budgeted uncertainty sets, which are essentially interval uncertainty sets where not all coefficients may deviate to the upper bound simultaneously, i.e.,
\[ \cU_{B}^\Gamma = \{ \pmb{c}\in\mathbb{R}^n_+ : c_i = \underline{c}_i + d_i\delta_i,\ \sum_{i\in[n]} \delta_i \le \Gamma,\ \pmb{\delta}\in[0,1]^n\} \]
for some parameter $\Gamma\in[0,n]$.
This case is also known as continuous budgeted uncertainty, as opposed to discrete budgeted uncertainty, where each $\delta_i$ is binary. 
Note that this does not impact min-max problems, but may make a difference in two-stage settings.
Sometimes, a variant of this set is considered, where the budget $\Gamma$ affects the sum of deviations, that is,
\[ \cU_{Bvar}^\Gamma = \{ \pmb{c}\in\mathbb{R}^n_+ : c_i = \underline{c}_i + \delta_i,\ \sum_{i\in[n]} \delta_i \le \Gamma,\ \delta_i\in[0,d_i]\ \forall i\in[n]\} \]
In the following, we may also write $\overline{c}_i =  \underline{c}_i + d_i$ to denote the upper bound on the costs of item $i$ in interval or budgeted uncertainty sets.
	
Many more uncertainty sets have been considered, including ellipsoidal uncertainty sets \cite{ben1998robust,chassein2017minmax}, data-driven polyhedral sets based on statistical testing \cite{bertsimas2018data} or machine learning \cite{shang2017data}, or even combinations of multiple sets simultaneously \cite{bertsimas2019two,dokka2020mixed}. 
	
In Table~\ref{tab:survey} we present a non-exhaustive list of papers that have performed experimental studies within this framework for different types of robust optimization problem and uncertainty set. We note that in particular the combinations of min-max robust optimization with discrete uncertainty and min-max regret robust optimization with interval uncertainty have been frequently considered.

\begin{table}
	\begin{center}
		\begin{tabular}{r|ccc}
			& Interval Uncertainty & Discrete Uncertainty & Budgeted Uncertainty \\
			\hline
			\multirow{2}{*}{Min-max} & \multirow{2}{*}{\textemdash} & \cite{goerigk2019representative},\cite{goerigk2020generating} & \cite{hansknecht2018fast},
			\cite{monaci2013exact} \\
			& & \cite{song2012incomplete},\cite{taniguchi2008heuristic} &  \cite{busing2019formulations}
			\\
			\hline
			\multirow{3}{*}{Min-max regret} & \cite{montemanni2007robust},\cite{montemanni2004exact} & \multirow{3}{*}{\cite{chassein2018scenario}} & \multirow{3}{*}{\textemdash} \\
			& \cite{kasperski2012tabu},\cite{furini2015heuristic} & & \\
			& \cite{wu2018exact} & & \\
			\hline
			Two-stage & \textemdash & \textemdash & \cite{chassein2018recoverable}, \cite{zeng2013solving} \\
			\hline
			\multirow{2}{*}{Recoverable} & \textemdash & \cite{busing2011recoverable} & \cite{chassein2018recoverable},\cite{chassein2016recoverable},\\
			& & &  \cite{busing2019formulations}
		\end{tabular}
	\end{center}
	\caption{Experimental studies from the literature, combining different uncertainty sets with types of robust optimization problems.}\label{tab:survey}
\end{table}

Note that min-max robust optimization problems with interval uncertainty are equivalent to nominal problems, as
\[ \min_{\pmb{x}\in\X} \max_{\pmb{c}\in\cU_I} \pmb{c}^t\pmb{x} = \min_{\pmb{x}\in\X} \overline{\pmb{c}}^t\pmb{x}. \]
For this reason, min-max problems with interval uncertainty are excluded from the remainder of this study. Similarly, two-stage and recoverable robust problems with interval uncertainty are equivalent to problems with only a single scenario and are therefore not considered here, though they find some interest in the recent literature (see, e.g., \cite{kasperski2017robust,fischer2020investigation,bold2021recoverable}). 
Furthermore, while some papers consider two-stage problems with discrete scenarios, they are most commonly not given by an explicit list, but described implicitly, see, for example, \cite{hashemi2021exploiting,subramanyam2021lagrangian}.
	
We conclude this section by briefly recalling the high-level ideas of the HIRO approach from \cite{goerigk2020generating} to generate uncertainty sets. The aim is to avoid optimal robust solutions that clearly outperform alternative solutions. To this end, a given uncertainty set is modified, such that the value of the optimal robust solution is increased. This creates a shorter left-hand tail in the distribution of robust objective values. Computational experience shows that such approaches can increase the computational difficulty of instances by orders of magnitude for min-max problems with discrete uncertainty sets. In the following sections, we extend this methodology to other types of robust optimization problem and other types of uncertainty sets.
	
More formally, given an uncertainty set $\tilde{\cU}$ (e.g., produced by a random sampling method), the HIRO approach is to solve
\[ \max_{\cU \in N(\tilde{\cU})} \min_{\pmb{x}\in\X} \max_{\pmb{c}\in\cU} \pmb{c}^t \pmb{x} \]\label{HIRO}
where $N(\tilde{\cU})$ denotes a neighborhood of set $\tilde{\cU}$. In the original paper, the authors considered $\tilde{\cU} = \{\tilde{\pmb{c}}^1,\ldots,\tilde{\pmb{c}}^N\}$ and
$N(\tilde{\cU}) = \cU(\tilde{\pmb{c}}^1) \times \cU(\tilde{\pmb{c}}^2) \times \ldots \times \cU(\tilde{\pmb{c}}^N)$ with
\[ \cU(\tilde{\pmb{c}}^j) = \left\{ \pmb{c}\in\mathbb{R}^n_+ : c_i \in [ \max\{\tilde{c}^j_i-b,0\}, \min\{\tilde{c}^j_i + b, C\}]\ \forall i\in[n],\ \sum_{i\in[n]} c_i \le \sum_{i\in[n]} \tilde{c}^j_i \right\} \]
where $b\ge 0$ is a constant (the budget for change) and $C$ is a desired global upper bound on costs, i.e., each scenario can be modified by changing costs by up to $b$ units and not increasing the total sum of costs in this scenario. In this paper, we always use $C=100$. Due to the added maximization layer, HIRO is at least as hard as the robust optimization problem that is being changed. In \cite{goerigk2020generating} the authors show that the HIRO problem is $\Sigma^p_2$-complete. To solve it, an iterative solution procedure is possible, where we choose $\cU\in N(\tilde{\cU})$ in a master problem, and determine a robust solution $\pmb{x}$ in a subproblem.

\section{Min-Max Problems}
\label{sec:minmax}
	
We now consider robust optimization problems with the min-max criterion. We discuss how to generate instances with discrete or with budgeted uncertainty and compare solution times when solving these problems with CPLEX. For all experiments we use an Intel(R) Xeon(R) Gold 6154 CPU @ 3.00GHz computer with 754 GB RAM. In addition, only one core has been used for each instance generation or evaluation. Recall that all sampling methods are summarized in Appendix~\ref{app:overview} and generation times are presented in Appendix~\ref{app:times}.

\subsection{Discrete Uncertainty}
\label{subsec:minmax-discrete}
	
\subsubsection{Problem Statement}
	
We assume that the uncertainty set is modeled as a list of $N$ scenarios with $\cU_D = \{\pmb{c}^1,\ldots,\pmb{c}^N\}$. Using the epigraph reformulation, the corresponding mixed-integer program in case of the selection problem is as follows.
\begin{align*}
\min \; & \; t \\
\mathrm{s.t.} \; & \; t \geq \sum_{i \in [n]} c_i^j x_i & \forall j \in [N] \\
& \sum_{i \in [n]} x_i = p \\
& x_i \in \{0,1\} & \forall i\in[n]
\end{align*}
This problem is known to be weakly NP-hard even for $N=2$, and strongly NP-hard when $N$ is unbounded (see \cite{averbakh2001complexity,kasperski2009randomized}).

\subsubsection{Sampling}
\label{sec:mmsampling}
	
To generate uncertainty sets $\cU_D$, a basic (and often-used) approach is to choose values randomly and uniformly iid over some interval. We refer to this method as \MMD{U}. In this case, we choose each value $c^j_i$ from the set $\{1,\ldots,100\}$. We consider two additional sampling methods. In the first method \MMD{1}, we choose $c^j_i$ from the set $\{1,\ldots,10\}\cup\{91,\ldots,100\}$, i.e., values are either small or large. In the second method \MMD{2}, we create a symmetry in the items. The first $\lfloor n/2 \rfloor$ items have costs generated as with \MMD{U}, while the last $\lceil n/2\rceil$ items have costs $c^j_i = 100-c^j_{i-\lfloor n/2\rfloor}$.
	
\subsubsection{HIRO}
\label{sec:mm-d-hiro}
	
We use the HIRO approach to increase the hardness of instances for robust discrete optimization problems first introduced in \cite{goerigk2020generating}, where details of this method can be found. To solve the resulting max-min-max problem, an iterative solution method is used. For a current set of $K$ candidate solutions $\pmb{x}^k$, a scenario set is constructed by solving the following optimization problem:
\begin{align*}
\max \; & t \\
\mathrm{s.t.} \; &  t \leq \sum_{j\in[N]} \sum_{i\in[n]} d_{ijk} x^k_i & \forall k \in [K] \\
& \sum_{j\in[N]} \lambda_j^k = 1 & \forall k \in [K] \\
& d_{ijk} \leq c^j_i & \forall i\in[n], j\in[N], k\in[K] \\
& d_{ijk} \leq \overline{c}^j_i \lambda_j^k & \forall i\in[n], j\in[N], k\in[K] \\
& \pmb{c}^j \in \cU(\tilde{\pmb{c}}^j) & \forall j\in[N] \\
& \lambda_j^k \in \{0,1\} & \forall j\in[N], k\in[K]
\end{align*}
where sets $\cU(\tilde{\pmb{c}}^j)$ allow a deviation from the input scenario through a budget $b$ as defined in Section~\ref{sec:problems}. Here, $\pmb{\lambda}$ variables represent the assignment of worst-case scenario to solutions. Having constructed an uncertainty set, we solve the corresponding robust problem to find the next candidate solution $\pmb{x}$. This process is repeated until convergence or a time limit is reached. We refer to instances of this type as \MMD{XH-B}, where \texttt{X} refers to the sampling type of the input uncertainty set, and \texttt{B} to the budget.
	
\subsubsection{Experimental Setup}
	
We perform four experiments to test different parameter combinations. In Exp1, we consider the pairs $(20,11)$, $(25,13)$, $(30,15)$, $(35,17)$ and $(40,21)$ for $(n,p)$ and set $N=n$ in all cases. We chose $p$ to be odd as preliminary experiments showed higher computation times for \MMD{2} in such cases. In Exp2, we fix $n=N=30$ but change $p\in\{5,11,15,21,25\}$. In Exp3, we fix $n=30$ and $p=15$ but change $N\in\{5,10,15,20,25,30,35,40\}$. Finally, Exp4 considers large-scale problems with $N\in\{100,500,1000,5000,10000\}$ also using $n=30$ and $p=15$.
	
In each experiment, we generate 50 instances using \MMD{U}, \MMD{1} and \MMD{2}. Using each of these sets as input, we also apply \MMD{XH-B} with budgets $b\in\{1,2,5\}$ (recall that parameter $b$ controls by how much coefficients of an input set may be changed), except for Exp4. For instance generation, we apply a 600-second time limit. Thus, the number of instances in the four experiments are $5\times50\times3\times4=3000$ , $5\times50\times3\times4=3000$, $8\times50\times3\times4=4800$, and $5\times50\times3=750$, respectively.
	
All instances are solved using CPLEX-20-10 with a 600-second time limit.
	
\subsubsection{Experimental Results}
	
\begin{figure}[htbp]
	\includegraphics[width=\textwidth]{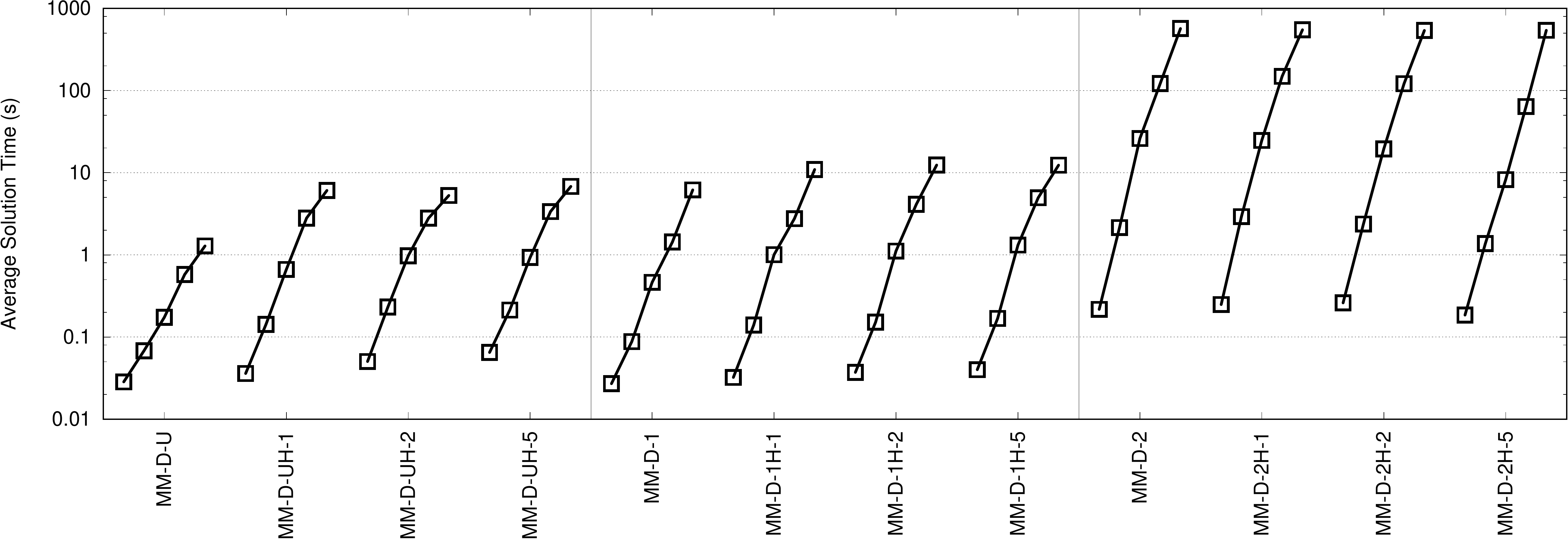}
	\caption{Min-max problems with discrete uncertainty for varying values of $(n,p)$, Exp1.}
	\label{fig:mm-dis-exp2-dataset1}
\end{figure}

\begin{figure}[htbp]
	\includegraphics[width=\textwidth]{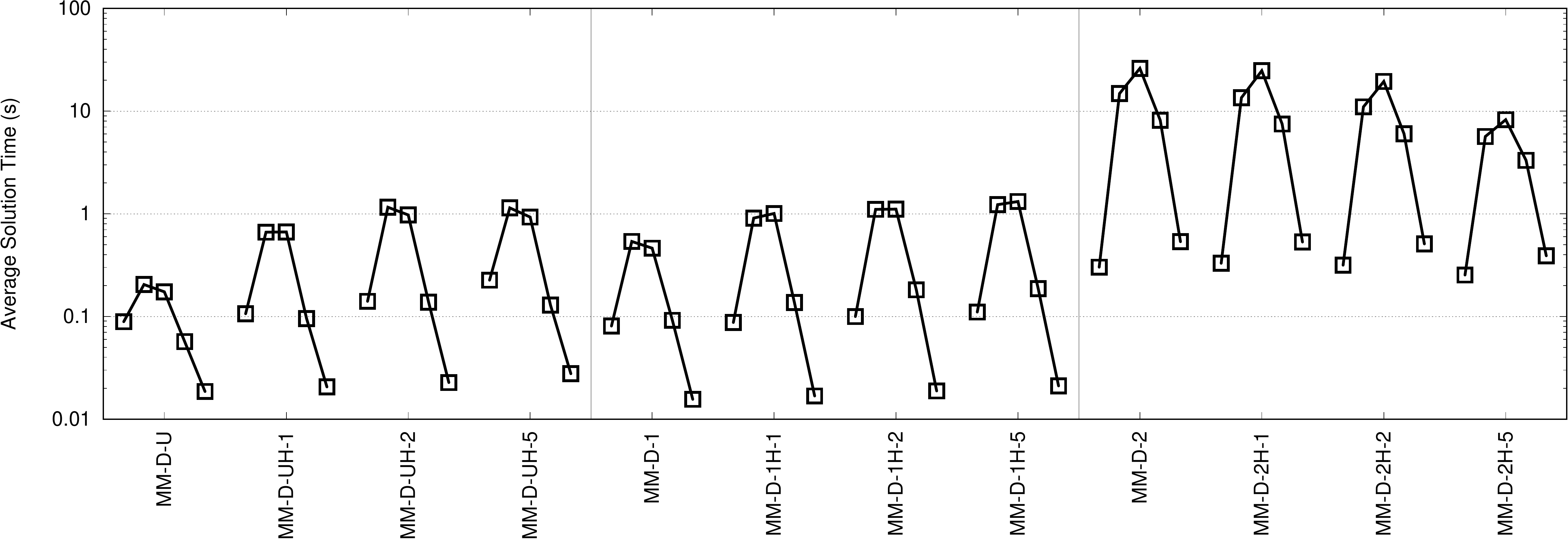}
	\caption{Min-max problems with discrete uncertainty for varying values of $p$, Exp2.}
	\label{fig:mm-dis-exp2-dataset2}
\end{figure}

\begin{figure}[htbp]
	\includegraphics[width=\textwidth]{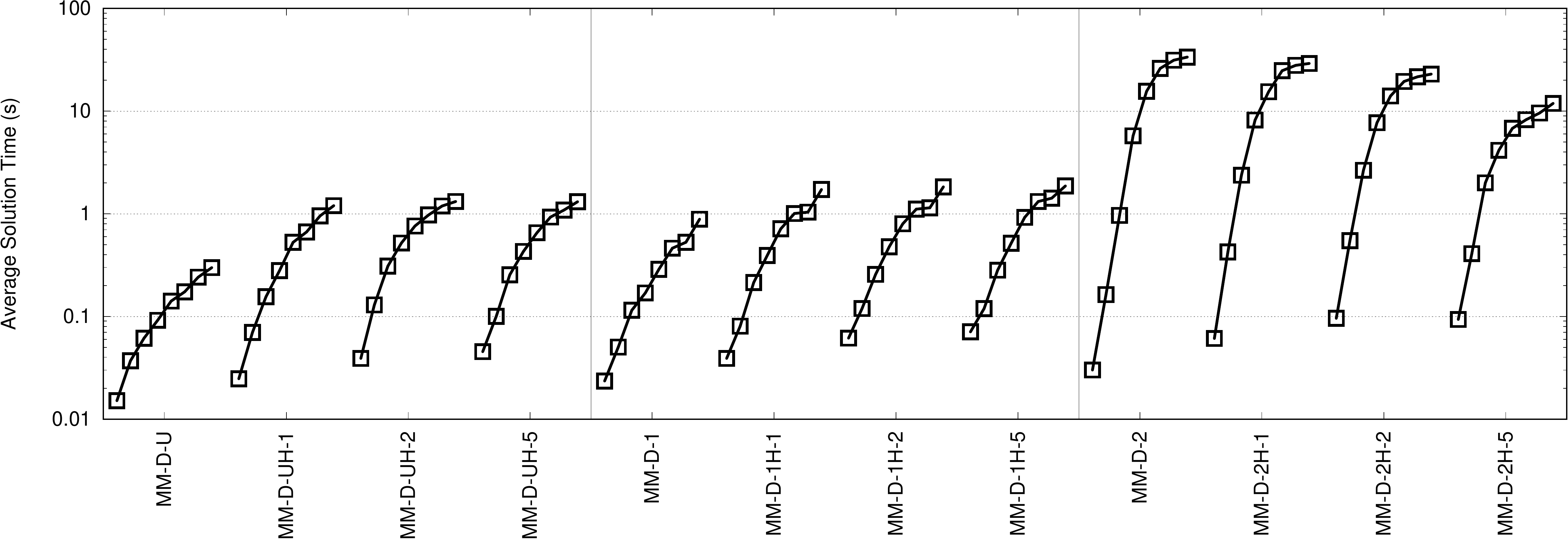}
	\caption{Min-max problems with discrete uncertainty for varying values of $N$, Exp3.}
	\label{fig:mm-dis-exp2-dataset3}
\end{figure}

\begin{figure}[htbp]
	\includegraphics[width=\textwidth]{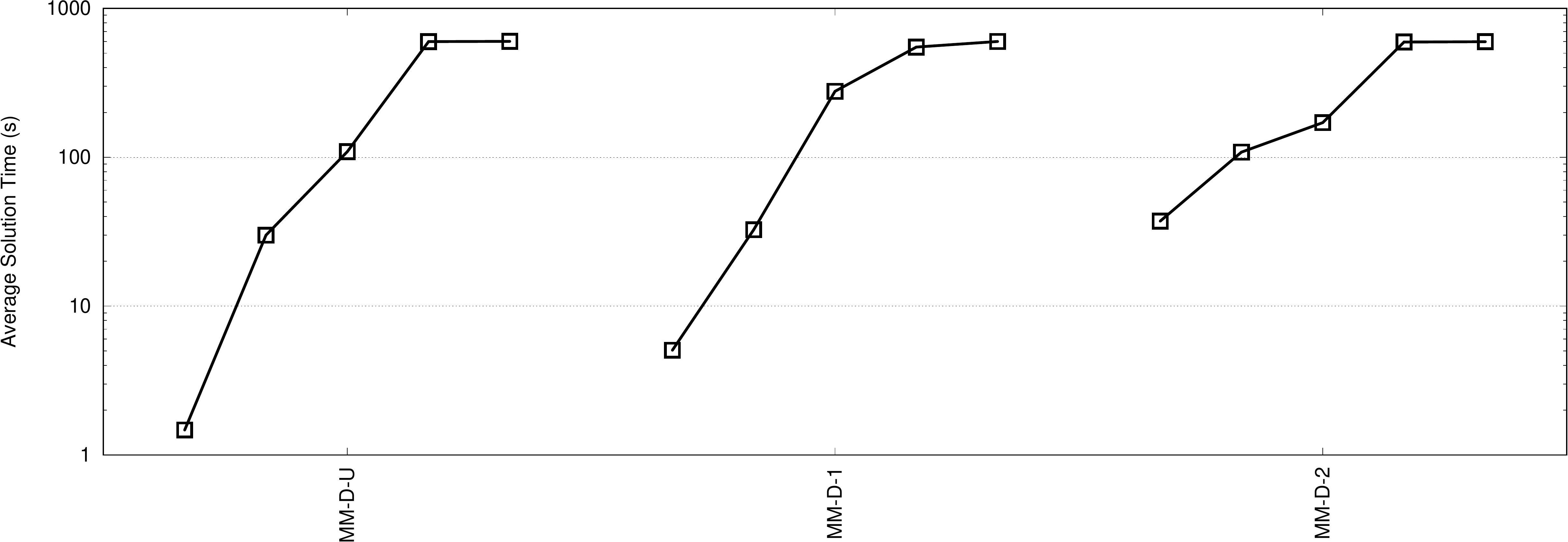}
	\caption{Min-max problems with discrete uncertainty for varying values of $N$, Exp4.}
	\label{fig:mm-dis-exp2-dataset4}
\end{figure}

Figures~\ref{fig:mm-dis-exp2-dataset1}--\ref{fig:mm-dis-exp2-dataset4} summarize the results of the four experiments. Each line connects average solution times for the parameter that is changed in the respective experiment.
		
Figure~\ref{fig:mm-dis-exp2-dataset1} demonstrates that while \MMD{U} instances can be solved in less than one second, \MMD{2} instances can hit the time limit of 600 seconds when $n=N=40 , p=21$. Also, it is observed that HIRO is effective for both \MMD{U} and \MMD{1} with all budgets $b$. As for \MMD{2}, the HIRO approach can slightly increase the hardness of instances with $b=1$, while it decreases the evaluation time of the instances with $b=2$ and $b=5$.  
	
Figure~\ref{fig:mm-dis-exp2-dataset2} shows that instances generated using sampling methods for the cases $p=11$ and $p=15$ are harder to solve in comparison to the other values of $p$. HIRO shows a similar performance as in Exp1.
	
Figure~\ref{fig:mm-dis-exp2-dataset3} implies that for the instances with lower number of scenarios HIRO also works for \MMD{2}. However, when the number of scenarios grows, HIRO decreases the difficulty of the instances, as before. Problem hardness increases with $N$, as can be expected.
	
Finally, Figure~\ref{fig:mm-dis-exp2-dataset4} illustrates that for large-scale problems (in terms of the number of scenarios), the gap between hardness of instances generated using \MMD{U}, \MMD{1} and \MMD{2} decreases. In some cases (e.g., $N=1000$), performance of \MMD{1} even exceeds performance of \MMD{2}. In addition, all instances generated with the three sampling methods hit the time limit for $N=5000$ and $N=10000$.
	
\subsection{Budgeted Uncertainty}
	
\subsubsection{Problem Statement}
	
By dualizing the inner worst-case problem, a compact MIP of the robust min-max selection problem under $\cU_B^\Gamma$ can  be formulated as follows:
\begin{align*}
\min_{\pmb{x}\in\X} \max_{\pmb{c}\in\cU_B^\Gamma} \pmb{c}^t \pmb{x} = \min\ & \sum_{i \in [n]} {\underline c}_{i} x_{i} \; + \; \Gamma \pi \; + \; \sum_{i \in [n]} \delta_{i}\\
\text{s.t. }  & \sum_{i\in[n]} x_i = p\\
& \pi \; + \; \delta_{i} \; \geq \; d_{i} x_{i} & \forall{i} \in [n]\\
& \pi \geq 0\\
& \delta_i \geq 0 & \forall{i} \in [n]\\
& x_i \in \{0,1\} & \forall{i} \in [n]
\end{align*}
In \cite{bertsimas2003robust} it is noted that candidate values of variable $\pi$ can be enumerated. As the resulting problems with fixed value of $\pi$ are equivalent to nominal problems, it is possible to solve this problem in polynomial time. More precisely, enumerating $\pi$ from the set $P = \{d_i : i\in[n]\} \cup \{0\} =: \{\pi^1, \dots, \pi^K\}$ with $K=|P|$ and substituting $\delta_i$, we can rewrite the problem as follows.
\[\min_{\pmb{x}\in\X} \max_{\pmb{c}\in\cU_B^\Gamma} \pmb{c}^t \pmb{x} = \min_{k\in[K]} \min_{\pmb{x}\in\X} (\sum_{i\in[n]} \underline{c}_{i} x_{i} + \Gamma \pi^{k} + \sum_{i\in[n]} [d_i - \pi^{k}]_{+} x_i)\]
	
\subsubsection{Sampling}
	
As a baseline sampling method, we use \MMB{U}, where values $\underline{c}_i$ and $d_i$ are chosen uniformly iid from $\{1,\ldots,100\}$. Additionally, we consider an approach where each $\underline{c}_i$ is chosen uniformly iid from $\{1,\ldots,100\}$, but we set $d_i = 100 - \underline{c}_i$, i.e., all items have the same upper bound $\underline{c}_i + d_i = 100$. We refer to this approach as \MMB{1}. Finally, we also choose $\underline{c}_i$ uniformly iid from $\{1,\ldots,10\}$ and $d_i$ uniformly from $\{100-\underline{c}_i-1,\ldots,100\}$. This means that items with smaller nominal costs $\underline{c}_i$ are more likely to have higher deviations $d_i$. This approach is denoted as \MMB{2}.

\subsubsection{HIRO}
	
We follow the HIRO idea to construct an optimization model to modify a given uncertainty set. For a fixed value $\Gamma$, note that the uncertainty is defined through the lower bounds $\underline{\pmb{c}}$ and deviations $\pmb{d}$. Let $\tilde{\underline{\pmb{c}}}$ and $\tilde{\pmb{d}}$ be given lower and upper bounds, respectively. A modification might affect just one of these vectors, or both. We consider all three possibilities. Whenever $\tilde{\underline{\pmb{c}}}$ or $\pmb{d}$ are changed, we require that the sum of lower bounds (or deviations) must not increase, and that the lower bounds (or deviations) are changed by up to $b$ units in each item.
	
\begin{enumerate}
		
	\item Modification of lower bounds.
		
	As the set $P=\{d_i : i\in[n]\} \cup \{0\}$ is only affected by the choice of $\pmb{d}$ we can dualize the inner nominal problems to obtain a compact linear programming formulation of HIRO:
	\begin{align*}
	\max_{\underline{\pmb{c}}} \min_{\pmb{x}\in\X} \max_{\pmb{c}\in\cU^\Gamma_B} \pmb{c}^t \pmb{x} = \max\ & t\\
	\text{s.t. } & t \leq \Gamma \pi^k + p \cdot \alpha^{k} - \sum_{i\in[n]} \beta_{i}^{k} & \forall k\in[K]\\
	& \alpha^{k} - \beta_{i}^{k} \leq {\underline{c}}_{i} + [d_i - \pi^k]_{+} & \forall i\in[n],k\in[K]\\
	& \underline{\pmb{c}} \in \cU(\underline{\tilde{\pmb{c}}}) \\
	& \alpha^k \geq 0 & \forall k\in[K]\\
	& \beta_{i}^{k} \geq 0 & \forall i\in[n], k\in[K]
	\end{align*}
		
	\item Modification of deviations.
		
	As $P$ depends on $\pmb{d}$, we enumerate $\pi$ through the index $k$ of the item where $\pi=d_k$. The resulting HIRO problem is then given as:
	\begin{align*}
	\max_{\pmb{d}} \min_{\pmb{x}\in\X} \max_{\pmb{c}\in\cU^\Gamma_B} \pmb{c}^t \pmb{x} = \max\ & t\\
	\text{s.t. } & t \leq \Gamma d_k + p \cdot \alpha^{k} - \sum_{i\in[n]} \beta_{i}^{k} & \forall k\in[n]\\
	& \alpha^{k} - \beta_{i}^{k} \leq {\underline{c}}_{i} + [d_{i} - d_{k}]_{+} & \forall i,k\in[n]\\
	& t \leq p \cdot \alpha^{0} - \sum_{i\in[n]} \beta_{i}^{0}\\
	& \alpha^{0} - \beta_{i}^{0} \leq {\underline{c}}_{i} + d_{i} & \forall i\in[n]\\
	& \pmb{d} \in\cU(\tilde{\pmb{d}}) \\
	& \alpha^0, \alpha^k \geq 0 & \forall k\in[n]\\
	& \beta^0_i, \beta_{i}^{k} \geq 0 & \forall i\in[n], k\in[n]
	\end{align*}
	To avoid the nonlinearity in $[d_i - d_k]_+$, we assume without loss of generality that all $\tilde{d}_i$ values are sorted in non-decreasing order. In the HIRO model, we require that this order must be preserved. That is, we can simplify 
	\begin{align*}
	[d_i - d_k]_+ = 
	\begin{cases}
	d_i - d_k & : \; k < i \\
	0 & : \; i \leq k \\
	\end{cases}
	\end{align*}
	and add the constraint
	\[d_1 \leq \cdots \leq d_n. \]
\end{enumerate}
	
The modification of both lower bounds and deviations is possible through a model that combines both models above.

\subsubsection{Experimental Setup}
We use sampling methods \MMB{U}, \MMB{1} and \MMB{2} to generate instances with $n=40$, $p=20$ and $\Gamma \in \{5,10,15,20\}$. For all combinations of method and parameter we generate 50 instances. Then, we use generated instances as input to each of the HIRO formulations (solving with respect to $\underline{\pmb{c}}$, $\pmb{d}$, or both) and solve them using CPLEX-20-10 with budgets $b\in\{1,2,5,10,20\}$ and a 600-second time limit. All problem instances are then solved using a 600-second time limit as well.

\subsubsection{Experimental Results}
	
\begin{figure}[htbp]
	\includegraphics[width=\textwidth]{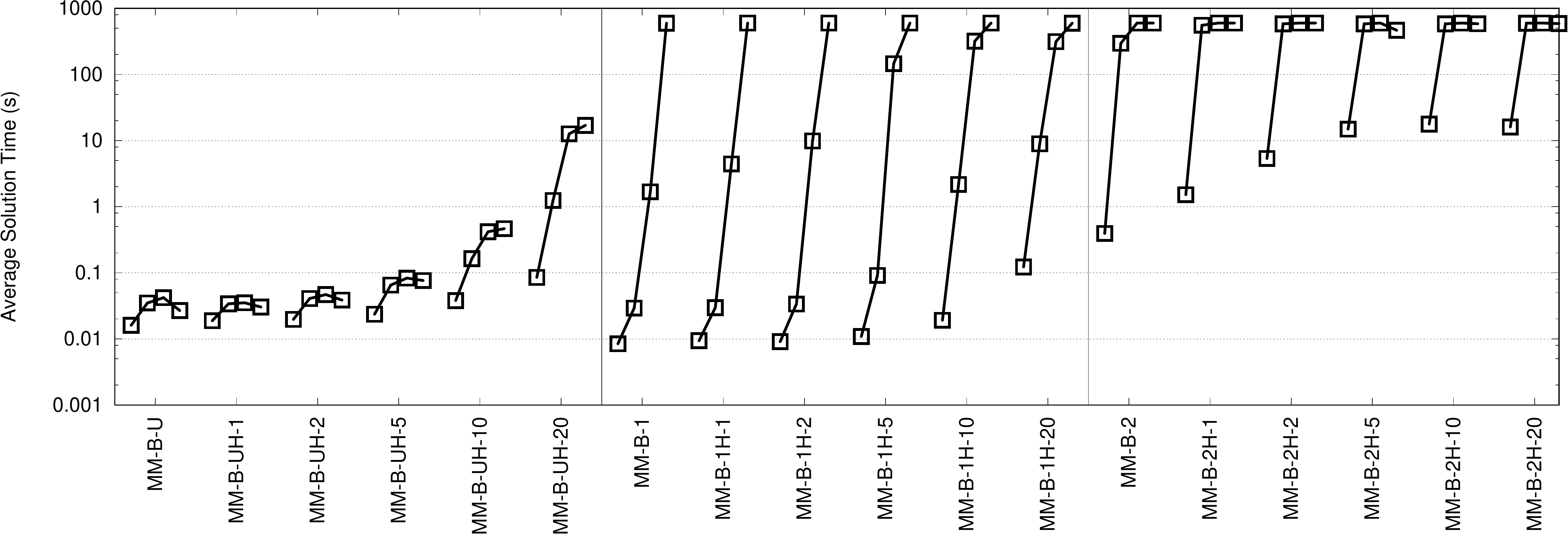}
	\caption{Min-max problems with budgeted uncertainty for varying $\Gamma$. HIRO modifies $\underline{\pmb{c}}$.}
	\label{fig:mm-bud-c}
\end{figure}

\begin{figure}[htbp]
	\includegraphics[width=\textwidth]{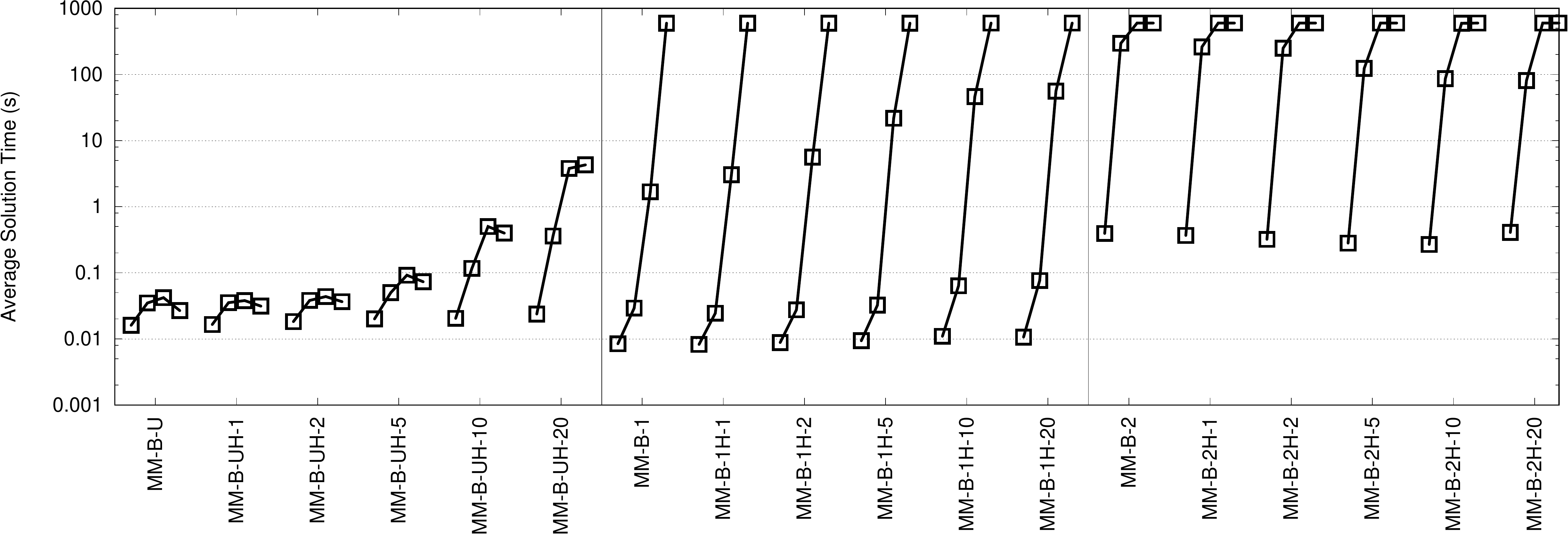}
	\caption{Min-max problems with budgeted uncertainty for varying $\Gamma$. HIRO modifies $\pmb{d}$.}
	\label{fig:mm-bud-d}
\end{figure}

\begin{figure}[htbp]
	\includegraphics[width=\textwidth]{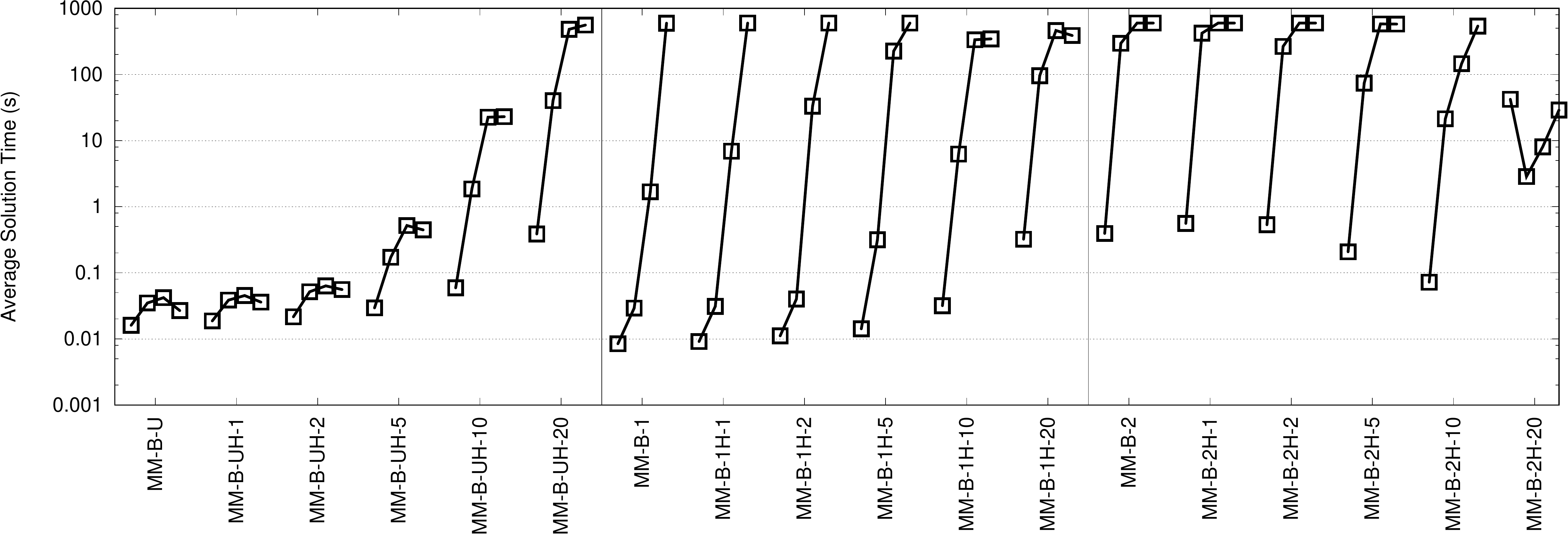}
	\caption{Min-max problems with budgeted uncertainty for varying $\Gamma$. HIRO modifies $\underline{\pmb{c}}$ and $\pmb{d}$.}
	\label{fig:mm-bud-cd}
\end{figure}

In Figures~\ref{fig:mm-bud-c}-\ref{fig:mm-bud-cd}, each line connects the average solution times for increasing values of $\Gamma$.
It can be seen that while \MMB{U} instances are relatively easy to solve ($<0.1$ seconds on average), finding the solution of instances generated by \MMB{1} and \MMB{2} are considerably more time consuming ($>1000$ seconds on average). In particular, when $\Gamma$ increases the instances hit the given time limit.
	
All three HIRO methods are effective over both \MMB{U} and \MMB{1} instances. The figures represent better results when $b$ increases. The best case for \MMB{U} happens when modification is done for both $\underline{\pmb{c}}$ and $\pmb{d}$ with $b=\Gamma=20$. However, HIRO gives different results when applied to \MMB{2} instances. Although modifying $\underline{\pmb{c}}$ always leads to the better results by increasing $b$, modification over $\pmb{d}$ or both $\underline{\pmb{c}}$ and $\pmb{d}$ for these instances makes them generally easier to solve, with few exceptions.

\section{Min-Max Regret Problems}
\label{sec:minmaxregret}

\subsection{Interval Uncertainty}
	
\subsubsection{Problem Statement}
	
To calculate the regret $\max_{\pmb{c}\in\cU} (\pmb{c}^t\pmb{x} - opt(\pmb{c}))$, we can find a worst-case scenario (i.e., a maximizer of this problem) by setting
\[\pmb{c}^{rwc}_i(\pmb{x}) = 
\begin{cases}
\underline{c}_i + d_i & \text{ if } x_i=1\\
\underline{c}_i & \text{ if } x_i=0
\end{cases}
\]
see, e.g. \cite{aissi2009min}. Accordingly, the min-max regret problem can be reformulated:
\[
\min_{\pmb{x}\in\X} \max_{\pmb{c}\in\cU} \left( \pmb{c}^t \pmb{x} - opt(\pmb{c}) \right) = \min_{\pmb{x}\in\X} \max_{\pmb{y}\in\X} \left( (\underline{c}_i + d_i) x_i - \sum_{i \in [n]} (\underline{c}_i + d_i x_i) y_i  \right)
\]
If $\X$ can be represented by a polyhedron with compact description, we can apply duality in the inner maximization problem to find a compact MIP formulation. In case of the selection problem, this leads to the following regret problem.
\begin{align*}
\min\ & \sum_{i\in[n]} (\underline{c}_i+d_i) x_i - p\pi + \sum_{i\in[n]} \rho_i \\
\text{s.t. } & \pi - \rho_i \le d_i x_i + \underline{c}_i & \forall i\in[n] \\
& \sum_{i\in[n]} x_i = p \\
& x_i \in \{0,1\} & \forall i\in[n] \\
& \pi \ge 0 \\
& \rho_i \ge 0 & \forall i\in[n] 
\end{align*}
It is possible to solve this problem in polynomial time, see \cite{averbakh2001complexity}.
	
\subsubsection{Sampling}
	
We use a uniform sampling method (\MMRI{U}) as a baseline, where $\underline{c}_i$ and $d_i$ are chosen iid from $\{1,\ldots,100\}$. Two additional sampling methods are considered. In \MMRI{1}, we choose $\underline{c}_i$ uniformly from $\{1,\ldots,10\}$ and $d_i$ from $\{91,\ldots,100\}$ with probability $0.5$. Otherwise, we choose $\underline{c}_i \in\{91,\ldots,100\}$ and $d_i \in \{1,\ldots,10\}$, i.e., items either have low costs $\underline{c}_i$ but high deviations $d_i$ or vice versa. In \MMRI{2}, we follow a similar approach, but sample items that either have both small $\underline{c}_i$ and $d_i$, or both large $\underline{c}_i$ and $d_i$.

\subsubsection{HIRO}
	
In the MIP formulation of the selection problem, note that there always exists an optimal solution where $\pi \in P = \{0\} \cup \{\underline{c}_i : i\in[n]\} \cup \{d_i : i\in[n]\}$. Let us write $P=\{\pi^1,\ldots,\pi^K\}$ with $K=|P|$. For fixed choice $\pi^k\in P$, the problem becomes
\begin{align*}
\min & \sum_{i\in[n]} (\underline{c}_i + d_i) x_i - p\pi^k + \sum_{i\in[n]}([\pi^k - \underline{c}_i - d_i]_+ - [\pi^k - \underline{c}_i]_+) x_i + \sum_{i\in[n]} [\pi^k - \underline{c}_i]_+\\
\mathrm{s.t.} & \sum_{i\in[n]} x_i \geq p\\
& x_i \in \{0,1\} &\forall i\in[n]
\end{align*}
The dual of this problem can be written as
\begin{align*}
\max \quad & p \alpha^k - \sum_{i \in [n]} \beta^k_{i} - p \pi^k + \sum_{i \in [n]} [\pi^k - \underline{c}_i]_+\\
\mathrm{s.t.} \quad &\alpha^k - \beta^k_i \leq \underline{c}_i + d_i + [\pi^k - \underline{c}_i - d_i]_+ - [\pi^k - \underline{c}_i]_+ & \forall i \in [n]\\ 
& \alpha^k \geq 0 \\
& \beta_{i}^{k} \geq 0 & \forall i\in[n]
\end{align*}
To formulate the HIRO approach, we introduce binary variables to linearize the nonlinearities of this model. Let $z^k_i=1$ if $\pi^k-\underline{c}_i \ge 0$ and $q^k_i=1$ if $\pi^k-\underline{c}_i-d_i \ge 0$. The linearized model is then as follows.
\begin{align*}
\max\ & t \\
\text{s.t. } &t \leq  \sum_{i\in[n]} (z_i^k \pi^k - \hat{z}_i^k) - p\pi^k + p\alpha^k - \sum_{i\in[n]} \beta_i^k & \forall k\in[K]\\
& \alpha^k - \beta_i^k \leq \underline{c}_i + d_i + \pi^k q^k_i - \hat{q}^k_i - \tilde{q}^k_i - s^k_i & \forall i\in[n], k\in[K] \\
& s_i^k \geq \pi^k - \underline{c}_i & \forall i\in[n], k\in[K]\\
& \hat{z}_i^k \geq \underline{c}_i - ({\underline{\tilde{c}}}_{i} + b) (1-z_i^k) & \forall i\in[n], k\in[K]\\
& \hat{q}_i^k \geq \underline{c}_i - ({\underline{\tilde{c}}}_{i} + b) (1-q_i^k) & \forall i\in[n], k\in[K]\\
& \tilde{q}_i^k \geq \underline{c}_i - ({\tilde{d}}_{i} + b) (1-q_i^k) & \forall i\in[n], k\in[K]\\
& \underline{\pmb{c}} \in \cU(\tilde{\underline{\pmb{c}}}) \\
& \pmb{d} \in \cU(\tilde{\pmb{d}}) \\
& z_i^k , q_i^k \in \{0,1\} & \forall i\in[n], k\in[K]\\
& \alpha^k \geq 0 & \forall k\in[K]\\
& s_i^k , \beta_{i}^{k} , \hat{z}_i^k , \hat{q}_i^k , \tilde{q}_i^k \geq 0 & \forall i\in[n], k\in[K]\\
\end{align*}

\subsubsection{Experimental Setup}\label{Exp-Reg-Int}
	
We generate instances with $n=100$ and vary $p\in\{10,20,\ldots,90\}$. For each sampling method, 50 instances are generated. These instances are also used as an input to the HIRO method with budgets $b\in\{1,2,5\}$ and 600-second generation time limit. The total number of instances generated is therefore equal to $3 \times 9 \times 50 \times 4 = 5400$. All instances are solved using CPLEX-20-10 with a time limit of 600 seconds.

\subsubsection{Experimental Results}
	
\begin{figure}[htbp]
	\includegraphics[width=\textwidth]{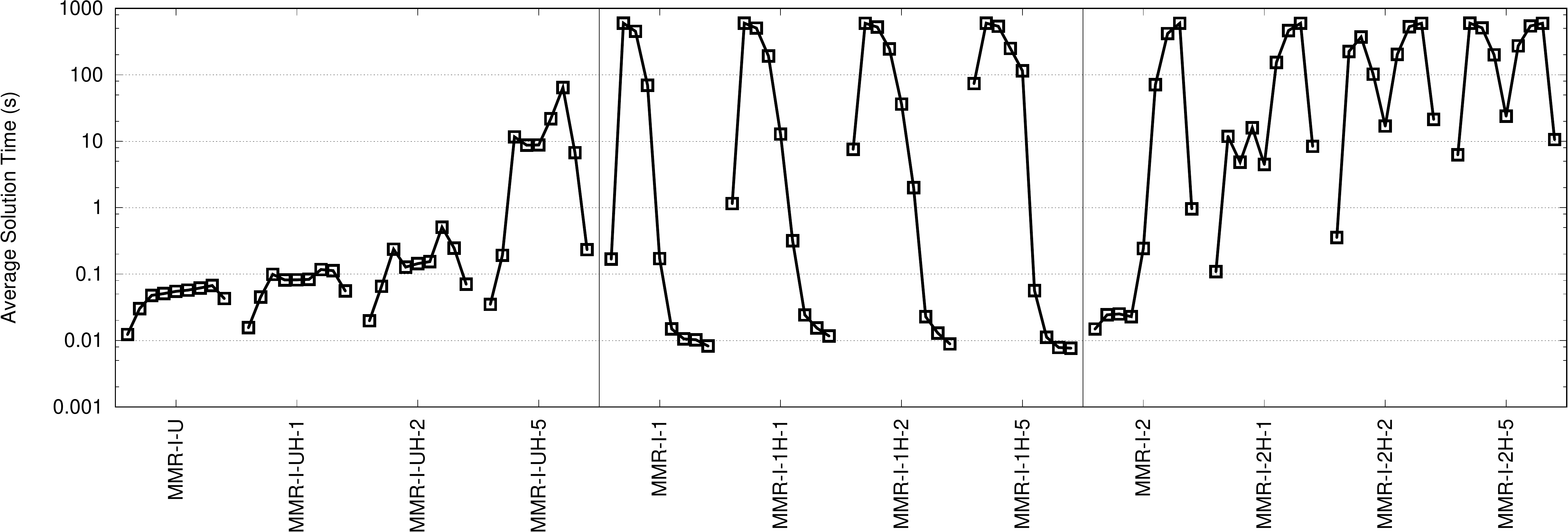}
	\caption{Min-max regret problems with interval uncertainty for varying values of $p$.}
	\label{fig:reg-int}
\end{figure}

Figure~\ref{fig:reg-int} summarizes the results of this experiment. Each line connects the average solution times for increasing values of $p$. For small values of $p$, \MMRI{1} generates instances that are considerably harder to solve than those generated by \MMRI{U} (time limit versus $<0.1$ seconds), while for larger values of $p$, method \MMRI{2} produces harder instances. The HIRO approach can successfully increase the hardness of instances. In particular, HIRO increases the solution time of \MMRI{U} more than 100 times for $p \in \{30,40,\ldots,80\}$ when $b=5$. In addition, \MMRI{2H-5} generates instances that are hard for both small and large values of $p$, thus resulting in instances that combine the best features of the other methods.
	
\subsection{Discrete Uncertainty}
	
\subsubsection{Problem Statement}
	
We consider using the min-max regret criterion for the selection problem with discrete uncertainty sets $\cU_D = \{\pmb{c}^1,\ldots,\pmb{c}^N\}$. By calculating the optimal objective values $opt(\pmb{c}^j)$ for each scenario $j\in[N]$, the corresponding MIP formulation is as follows.
\begin{align*}
\min\ & t \\
\text{s.t. } & t \geq \sum_{i\in[n]} c_i^j x_i - opt(\pmb{c}^j)  & \forall{j}\in[N]\\
& \sum_{i\in[n]}x_i = p \\
&x_i \in \{0,1\} & \forall{i}\in[n]
\end{align*}
Similar to the min-max case, this problem is weakly NP-hard even for $N=2$, and strongly NP-hard when $N$ is unbounded (see \cite{kasperski2016robust}).
	
\subsubsection{Sampling}
		
We consider the same sampling approaches as in Section~\ref{sec:mmsampling}. These are: uniform sampling in $\{1,\ldots,100\}$ (\MMRD{U}), sampling in $\{1,\ldots,10\}\cup\{91,\ldots,100\}$ (\MMRD{1}) and sampling with item symmetry (\MMRD{2}).
	
\subsubsection{HIRO}
	
We modify the model introduced in \cite{goerigk2020generating} for min-max problems with discrete uncertainty (see Section~\ref{sec:mm-d-hiro}). To this end, we introduce a new variable $\pmb{y}^j$ which represents an optimal solution under scenario $j\in[N]$. For a set $\{\pmb{x}^1,\ldots,\pmb{x}^K\}$ of candidate solutions, constructing an uncertainty set that maximizes the best robust objective value among $\pmb{x}^k$, $k\in[K]$, is formulated as follows:
\begin{align*}
\max \; & t \\
\mathrm{s.t.} \; & t \leq \sum_{j\in[N]} \lambda_j^k (\pmb{c}^j \pmb{x}^k - \pmb{c}^j \pmb{y}^j) & \forall k \in [K]\\
& \sum_{j\in[N]} \lambda_j^k = 1 & \forall k \in [K]\\
& \pmb{c}^j \in \cU(\tilde{\pmb{c}}^j) & \forall j \in [N]\\
& \pmb{y}^j \in \X & \forall j \in [N]\\
& \lambda_j^k \in \{0,1\} & \forall j \in [N], k \in [K] 
\end{align*}
Note that the formulation is nonlinear due to the products between variables $\pmb{\lambda}$, $\pmb{c}$, and $\pmb{y}$. We introduce new variables $d_{ijk} = \lambda_j^k c^j_i$ and 
$\alpha_{ijk} = \lambda_j^k c_i^j y_i^j$. The linearized MIP formulation we use is then as follows.
\begin{align*}
\max \; & t \\
\mathrm{s.t.} \; & t \leq \sum_{j\in[N]} \sum_{i\in[n]} (x_i^k d_{ijk} - \alpha_{ijk}) & \forall k \in [K]\\
& \sum_{j\in[N]} \lambda_j^k = 1 & \forall k \in [K]\\
& d_{ijk} \leq c_i^j & \forall i \in [n], j \in [N], k \in [K] \\
& d_{ijk} \leq \overline{c}_i^j \lambda_j^k & \forall i \in [n], j \in [N], k \in [K] \\
& \alpha_{ijk} \geq c_i^j - \overline{c}_i^j (2 - \lambda_j^k - y_i^j) & \forall i \in [n], j \in [N], k \in [K] \\
& \pmb{c}^j \in \cU(\tilde{\pmb{c}}^j) & \forall j \in [N]\\
& \pmb{y}^j \in \X & \forall j \in [N]\\
& \lambda_j^k \in \{0,1\} & \forall j \in [N], k \in [K] \\
& d_{ijk},\alpha_{ijk} \geq 0 & \forall i \in [n], j \in [N], k \in [K]
\end{align*}
	
\subsubsection{Experimental Setup} 
	
We follow a similar setup as in Section~\ref{subsec:minmax-discrete} and consider four experiments. In Exp1 all inputs change, which means that for parameters $(N,n,p)$, we consider combinations $(30,30,15)$, $(40,40,20)$ and $(40,40,21)$. In Exp2 we fix $n=N=30$ but change $p\in\{10,11,15,20,21\}$. In Exp3 we fix $n=30$ and $p=15$ but change $N\in\{20,30,40\}$. Finally, Exp4 has a similar setup as Exp3, but we consider large-scale problems with $N\in\{100,200,500,1000,2000,5000\}$.
	
For each parameter choice, we generate 50 instances using \MMRD{U}, \MMRD{1} and \MMRD{2}. All instances of Experiments 1,2 and 3 are used as an input to the HIRO formulation with a budget $b\in\{1,2,5\}$ and a 600-second generation time limit. Thus, Exp1 considers $3\times50\times3\times4=1800$ instances, Exp2 considers $5\times50\times3\times4=3000$ instances, and Exp3 considers $3\times50\times3\times4=1800$ instances. For Exp4, only the sampling methods are used, which means that we consider $6\times50\times3=900$ instances. All instances are solved using CPLEX-20-10 with a 600-second time limit.
	
\subsubsection{Experimental Results}
	
\begin{figure}[htbp]
	\includegraphics[width=\textwidth]{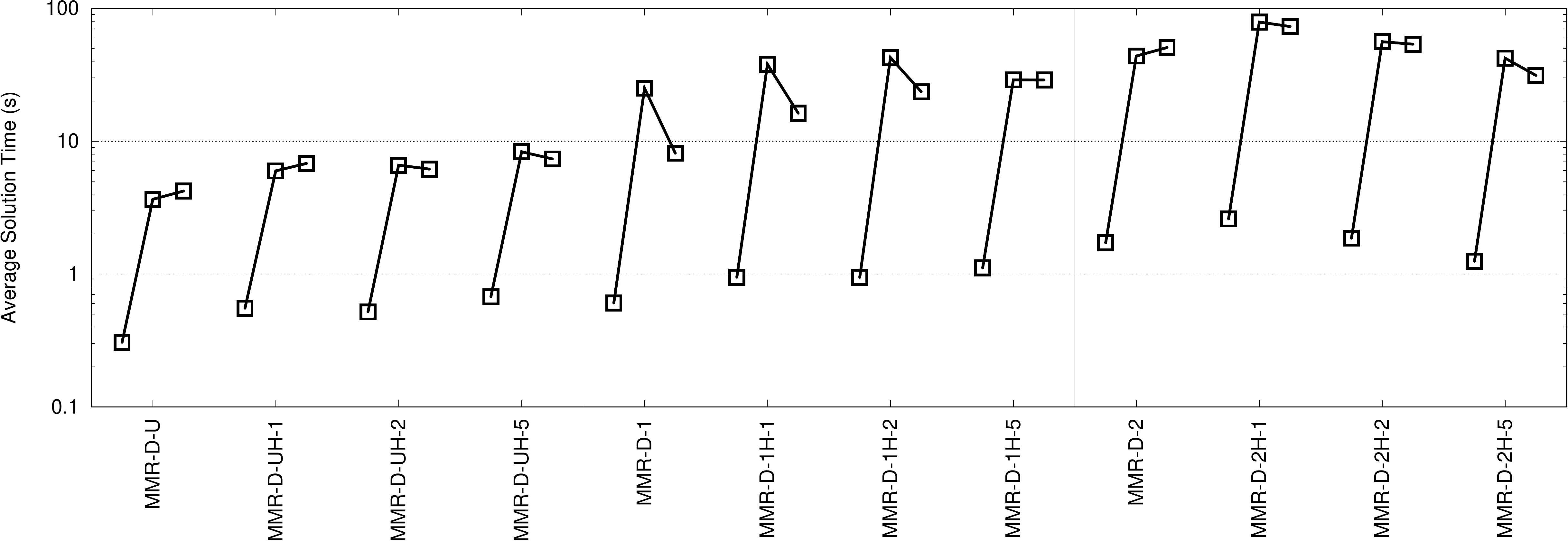}
	\caption{Min-max regret problems with discrete uncertainty for varying values of $(n,p,N)$, Exp1.}
	\label{fig:reg-dis-exp1}
\end{figure}

\begin{figure}[htbp]
	\includegraphics[width=\textwidth]{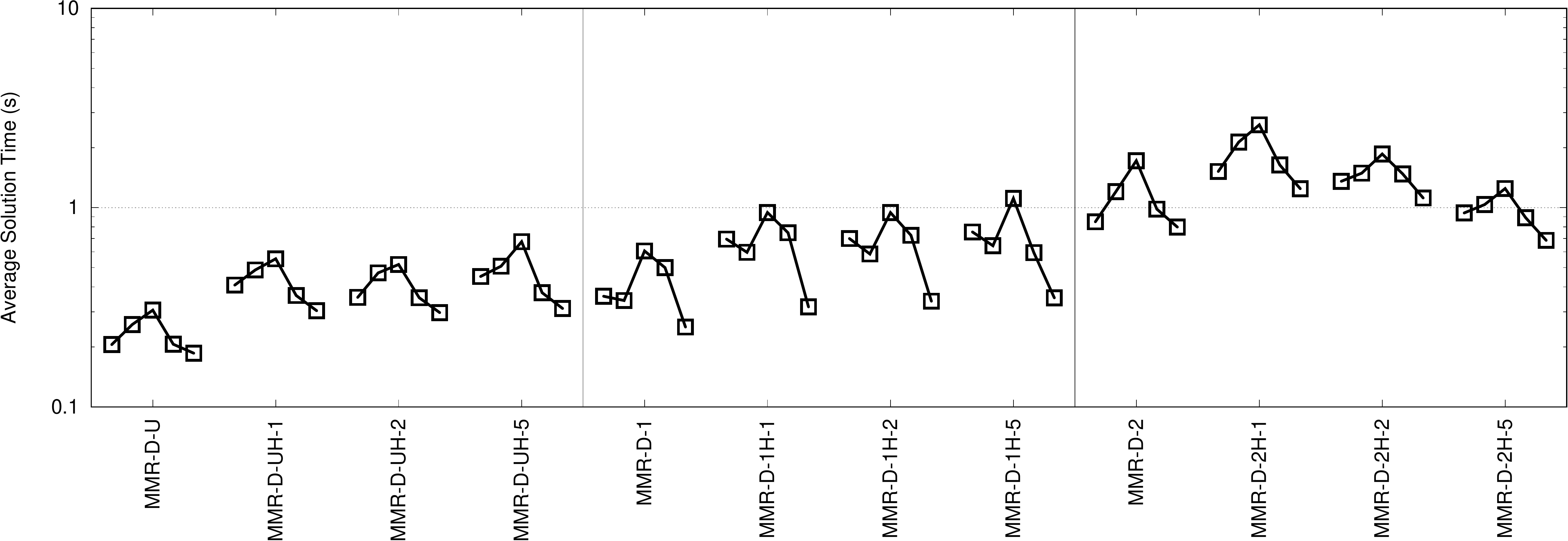}
	\caption{Min-max regret problems with discrete uncertainty for varying values of $p$, Exp2.}
	\label{fig:reg-dis-exp2}
\end{figure}

\begin{figure}[htbp]
	\includegraphics[width=\textwidth]{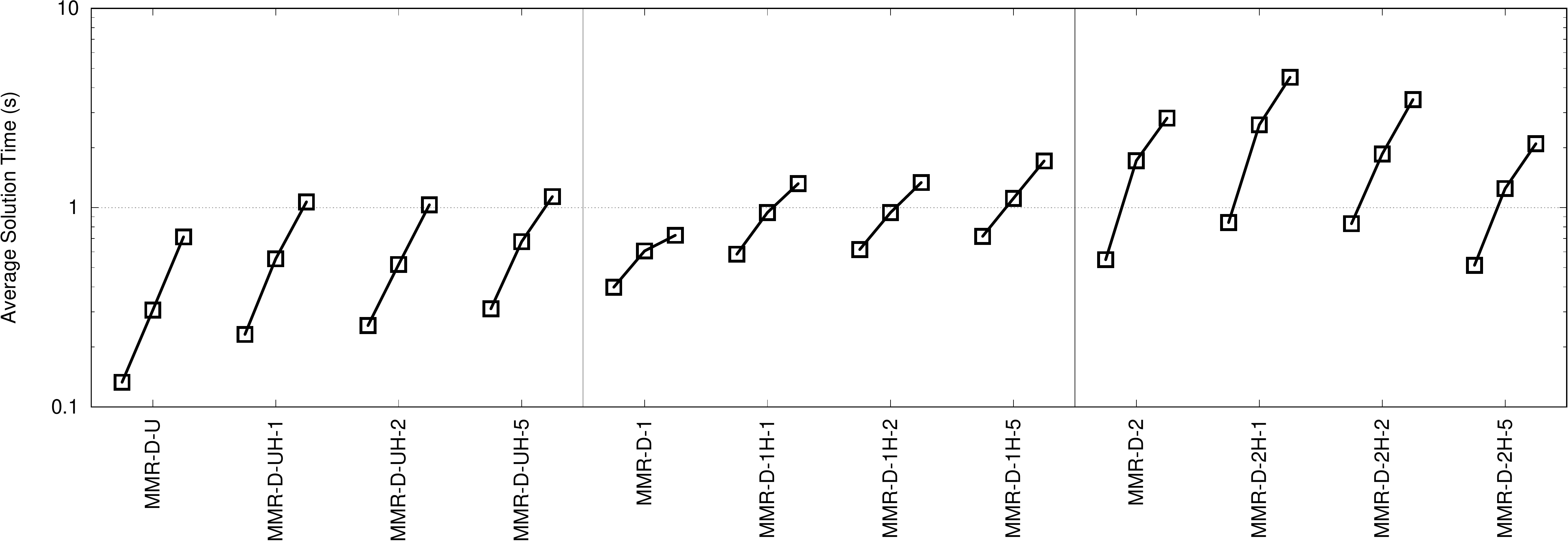}
	\caption{Min-max regret problems with discrete uncertainty for varying values of $N$, Exp3.}
	\label{fig:reg-dis-exp3}
\end{figure}

\begin{figure}[htbp]
	\includegraphics[width=\textwidth]{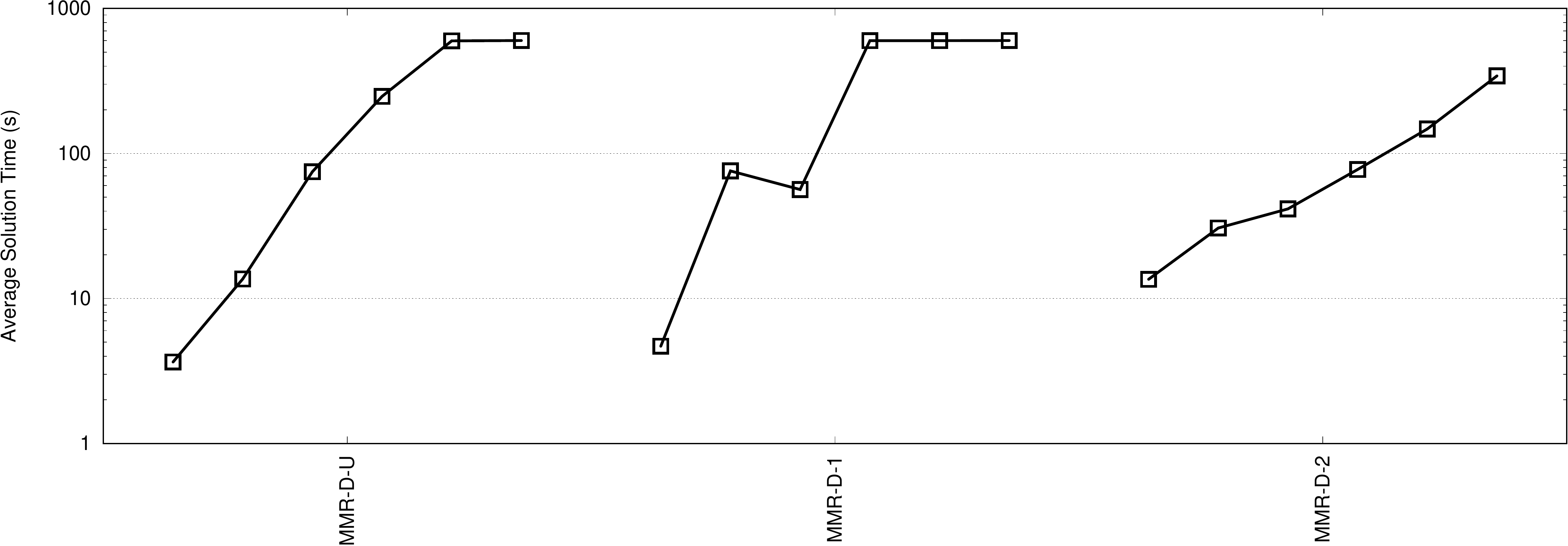}
	\caption{Min-max regret problems with discrete uncertainty for varying values of $N$, Exp4.}
	\label{fig:reg-dis-exp4}
\end{figure}

Figures~\ref{fig:reg-dis-exp1}, \ref{fig:reg-dis-exp2}, \ref{fig:reg-dis-exp3} and \ref{fig:reg-dis-exp4} summarize the results of Exp1, Exp2, Exp3 and Exp4, respectively. 
It can be seen that HIRO can slightly increase solution times of some instances
Similarly, HIRO shows mixed results for the other experiments as well.
In particular instances generated by \MMRD{2} can increase solution times significantly over \MMRD{U}. Interestingly, this does not apply for large-scale instances from Exp4, where \MMRD{1} tends to create instances that are harder to solve.

\section{Two-Stage Problems}\label{sec:twostage}
	
\subsection{Discrete Uncertainty}\label{subsec:twostage-discrete}
	
\subsubsection{Problem Statement}
	
By introducing a set of variables $\pmb{y}^j$ for each scenario $j\in[N]$, the two-stage problem with discrete uncertainty is formulated as follows.
\begin{align*}
\min\ & \sum_{i\in[n]} C_i x_i + t \\
\text{s.t. } & t \ge \sum_{i\in[n]} c^j_i y^j_i & \forall j\in[N] \\
& \sum_{i\in[n]} x_i + y^j_i = p & \forall j\in[N] \\
& x_i + y^j_i \leq 1 & \forall i\in[n], j\in[N] \\
& x_i \in\{0,1\} & \forall i\in[n] \\
& y^j_i \in\{0,1\} & \forall i\in[n], j\in[N]
\end{align*}
	
\subsubsection{Sampling}\label{subsubsec:twostage-discrete-sampling}

We consider the following sampling methods to determine first-stage costs $\pmb{C}$ and second-stage costs $\pmb{c}^j$ for all $j\in[N]$. As a baseline (\TSTD{U}), all values are chosen uniformly from the set $\{1,\ldots,100\}$. In \TSTD{1}, we choose $C_i$ uniformly either from from $\{45,\ldots,55\}$ with probability $0.5$, or from set $\{25,\ldots,75\}$ with probability $0.5$. For each second-stage cost $c^j_i$, we choose uniformly from one of the sets $\{C_i-5,\ldots,C_i+5\}$, $\{1,\ldots,10\}$ or $\{91,\ldots,100\}$ with probabilities $0.5$, $0.25$, and $0.25$, respectively. Finally, in \TSTD{2}, $C_i$ is chosen uniformly from $\{1,\ldots,100\}$ with probability $0.5$, or is set to 50 with probability $0.5$. If $C_i = 50$, then with probability equal to 0.5 we choose $c_i^j$ from $\{1,\ldots,10\}$ and with probability equal to 0.5 we choose $c_i^j$ from $\{91,\ldots,100\}$; otherwise, if $C_i\neq 50$, we choose $c_i^j$ from $\{C_i - 5,\ldots,C_i + 5\}$ and set negative values to be equal to zero.

\subsubsection{HIRO}\label{subsubsec:twostage-discrete-HIRO}
	
Given a set of candidate first-stage solutions $\pmb{x}^k$, $k\in[K]$, we would like to find first-stage costs $\pmb{C}$ and second-stage costs $\pmb{c}^j$, $j\in[N]$, such that the worst-case costs of the best candidate solution is maximized, i.e., to solve
\begin{align*}
\max\ & t \\
\text{s.t. } & t \le \sum_{i\in[n]} C_i x^k_i + \max_{j\in[N]} Q(\pmb{x}^k,\pmb{c}^j) & \forall k\in[K] \\
&\pmb{c}^j \in \cU(\tilde{\pmb{c}}^j) & \forall j\in[N] \\
& \pmb{C} \in \cU(\tilde{\pmb{C}})
\end{align*}
where $Q(\pmb{x}^k,\pmb{c}^j)$ denotes the second-stage costs of $\pmb{x}^k$ under scenario $\pmb{c}^j$. To calculate these for some scenario $\pmb{c}^j$, we solve
\begin{align*}
\min\ &\sum_{i\in[n]} c^j_i y^k_i \\
\text{s.t. } & \sum_{i\in[n]} x^k_i + y^k_i  = p\\
& x^k_i + y^k_i \le 1 & \forall  i\in[n] \\
& y^k_i \in\{0,1\}& \forall i\in[n]
\end{align*}
By relaxing this problem, dualizing it and then using variables $\pmb{\lambda}$ to decide which scenario is assigned to which solution to maximize the objective value, we find the following formulation for HIRO.
\begin{align*}
\max\ & t \\
\text{s.t. } & t \le \sum_{i\in[n]} C_i x^k_i + (p-\sum_{i\in[n]} x^k_i) \beta^k + \sum_{i\in[n]} (x^k_i - 1)\gamma^k_i &\forall k\in[K] \\
& \sum_{j\in[N]} \lambda^k_j = 1 & \forall k\in[K]\\
& \beta^k \le \gamma^k_i + \sum_{j\in[N]} c^j_i \lambda^k_j & \forall i\in[n],k\in[K] \\
& \lambda^k_j \in\{0,1\} & \forall j\in[N], k\in[K] \\
& \beta^k \ge 0 & \forall k\in[K] \\
& \gamma^k_i \ge 0  & \forall i\in[n],k\in[K] \\
&\pmb{C} \in \cU(\tilde{\pmb{C}})  \\
&\pmb{c}^j \in \cU(\tilde{\pmb{c}}^j) & \forall j\in[N]
\end{align*}
If we only optimize over first-stage costs $\pmb{C}$, we use nominal values $\tilde{\pmb{c}}^j$ instead of treating them as variables. However, if we optimize over both $\pmb{C}$ and $\pmb{c}^j$, the nonlinearity in $c^j_i \lambda^k_j$ needs to be replaced using variables $d_{ijk}$ (see Section~\ref{sec:mm-d-hiro}).

\subsubsection{Experimental Setup}\label{subsubsec:twostage-discrete-experiments}
	
We use a similar setup as for the single-stage min-max problem from Section~\ref{subsec:minmax-discrete}. We create four datasets to assess the effect of different parameters. In Exp1, we consider the two cases when $N=n=50,p=25$ and $N=n=100,p=50$. In Exp2, we fix $n=N=50$ but change $p\in\{10,20,25,30,40\}$. In Exp3, we fix $n=50$ and $p=25$ but change $N\in\{10,20,\ldots,60\}$. Finally, in Exp4, we consider problems with a large-scale number of scenarios $N\in\{100,200,500,1000,2000,5000\}$ with $n=50$ and $p=25$.
	
For each parameter configuration, we generate 50 instances using \TSTD{U}, \TSTD{1} and \TSTD{2}. For Experiments~1-3, we also use HIRO with budget $b\in\{1,2,5\}$ and a 600-second generation time limit. We use HIRO for modifying only first-stage costs $\pmb{C}$, or both first- and second-stage costs $\pmb{C}$ and $\pmb{c}^j$. Thus, we consider $2\times50\times3\times7=2100$ in Exp1, $5\times50\times3\times7=5250$ instances in Exp2, $6\times50\times3\times7=6300$ instances in Exp3, and $6\times50\times3=900$ instances in Exp4, respectively. All instances are solved using CPLEX-20-10 with a 600-second time limit.
	
\subsubsection{Experimental Results}\label{subsubsec:twostage-discrete-results}

\begin{figure}[htbp]
	\includegraphics[width=\textwidth]{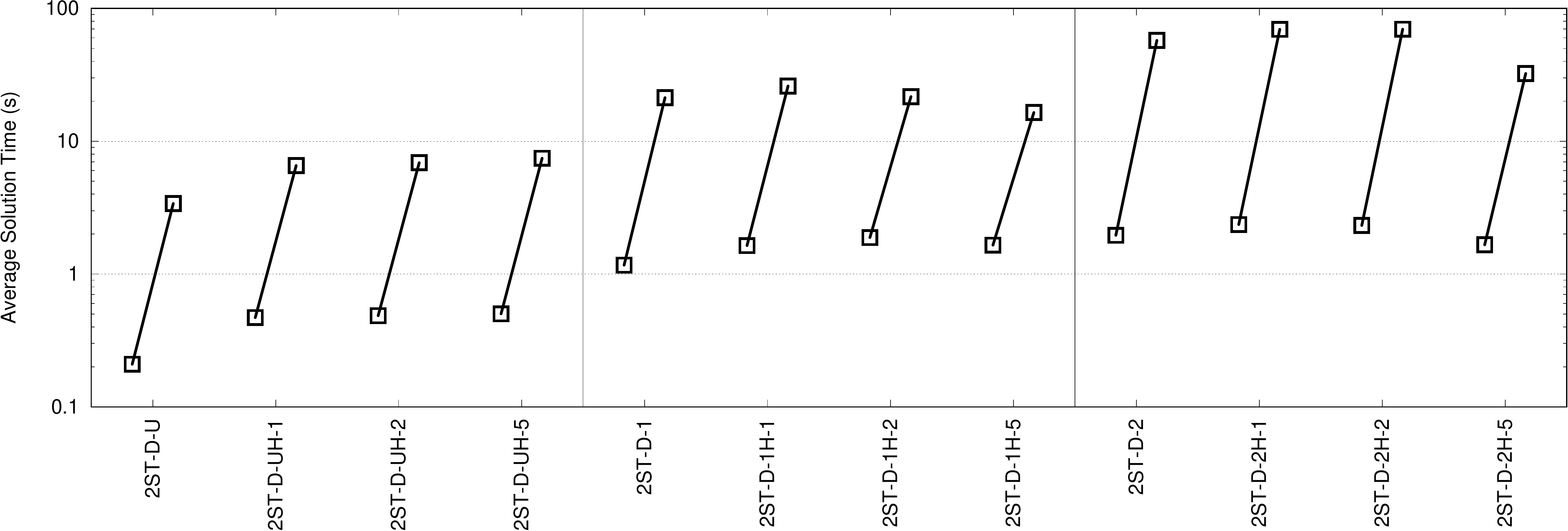}
	\caption{Two-stage problems with discrete uncertainty for varying values of $(n,p,N)$, Exp1. HIRO modifies $\pmb{C}$.}
	\label{fig:ts-dis-exp1-C}
\end{figure}

\begin{figure}[htbp]
	\includegraphics[width=\textwidth]{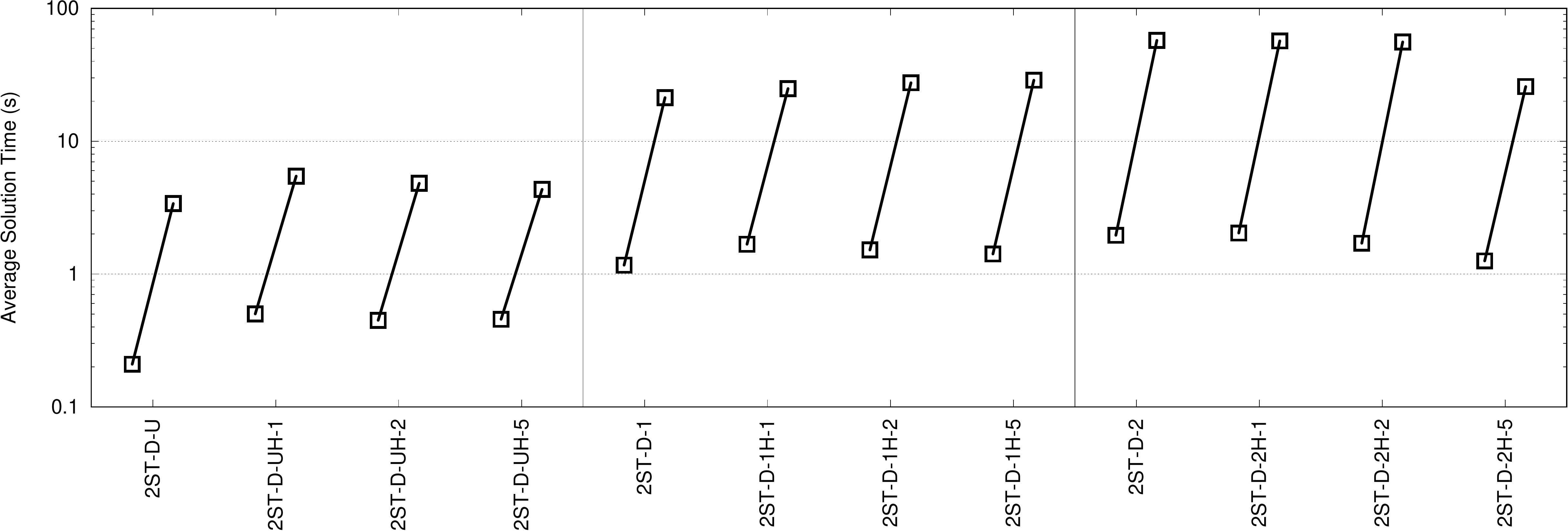}
	\caption{Two-stage problems with discrete uncertainty for varying values of $(n,p,N)$, Exp1. HIRO modifies $\pmb{C}$ and $\pmb{c}^j$.}
	\label{fig:ts-dis-exp1-Cc}
\end{figure}

In Figures~\ref{fig:ts-dis-exp1-C} and \ref{fig:ts-dis-exp1-Cc}, we show results on Exp1.
HIRO modifying $\pmb{C}$ or both $\pmb{C}$ and $\pmb{c}^j$ for all $j\in[N]$ can increase the solution time of \TSTD{U} instances. In addition, in most cases it is possible to slightly increase the solution time of \TSTD{1} instances, while this approach is less effective for most \TSTD{2} instances, which result in highest solution times overall.
	
\begin{figure}[htbp]
	\includegraphics[width=\textwidth]{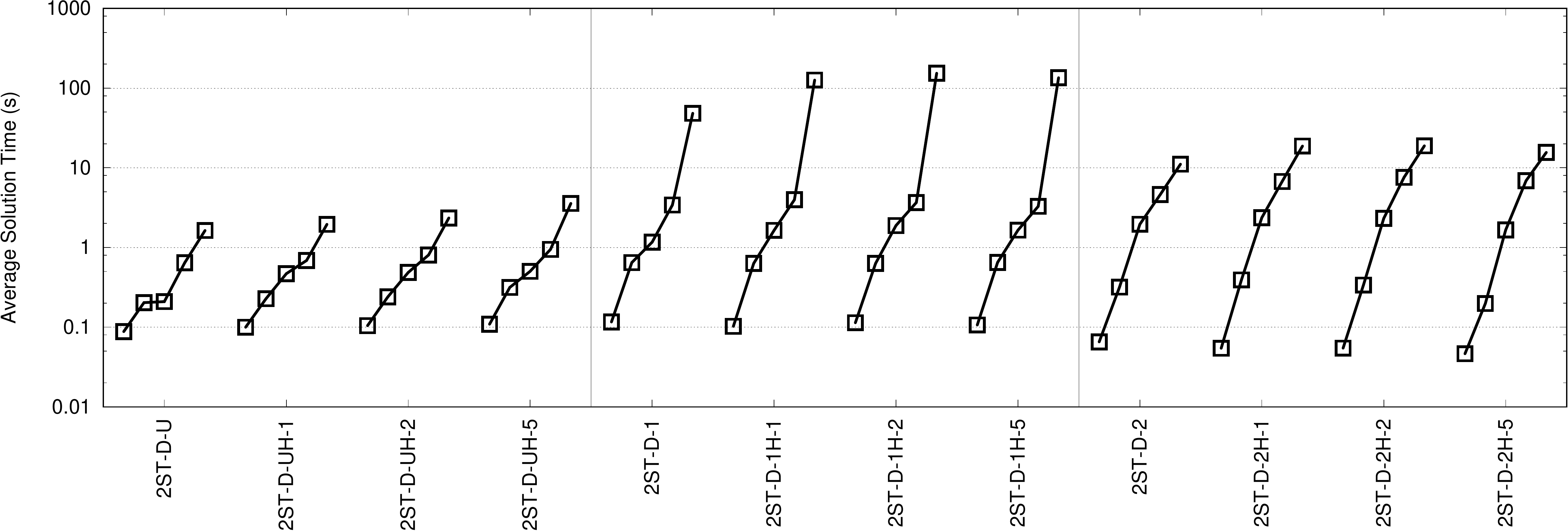}
	\caption{Two-stage problems with discrete uncertainty for varying values of $p$, Exp2. HIRO modifies $\pmb{C}$.}
	\label{fig:ts-dis-exp2-C}
\end{figure}

\begin{figure}[htbp]
	\includegraphics[width=\textwidth]{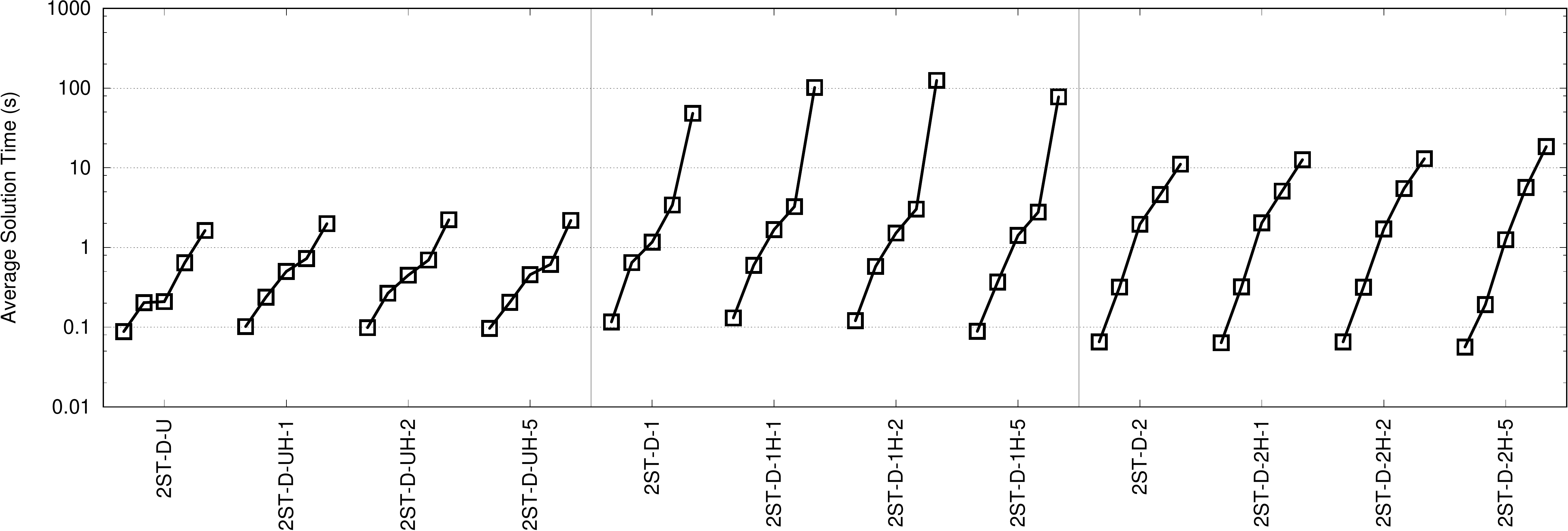}
	\caption{Two-stage problems with discrete uncertainty for varying values of $p$, Exp2. HIRO modifies $\pmb{C}$ and $\pmb{c}^j$.}
	\label{fig:ts-dis-exp2-Cc}
\end{figure}

The results shown in Figures~\ref{fig:ts-dis-exp2-C} and \ref{fig:ts-dis-exp2-Cc} on Exp2 show that larger values of $p$ result in higher solution times. HIRO can slightly increase solution times on most instances, while particularly \TSTD{1} gives good results for large values of $p$.
	
\begin{figure}
	\includegraphics[width=\textwidth]{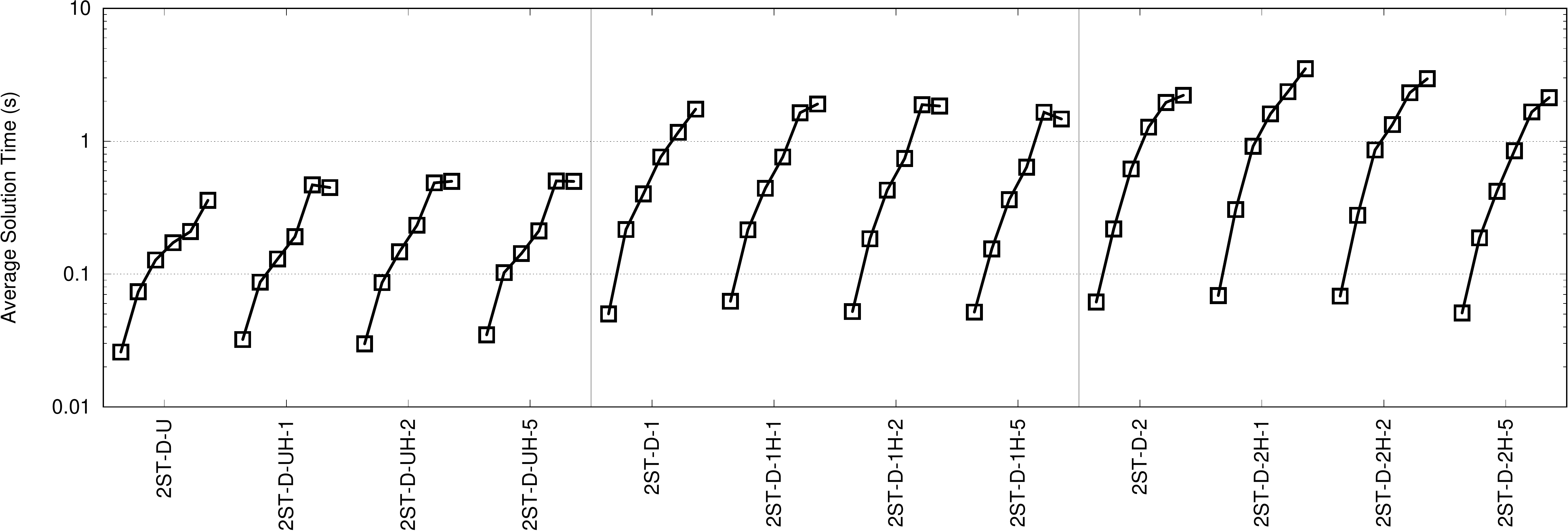}
	\caption{Two-stage problems with discrete uncertainty for varying values of $N$, Exp3. HIRO modifies $\pmb{C}$.}
	\label{fig:ts-dis-exp3-C}
\end{figure}

\begin{figure}[htbp]
	\includegraphics[width=\textwidth]{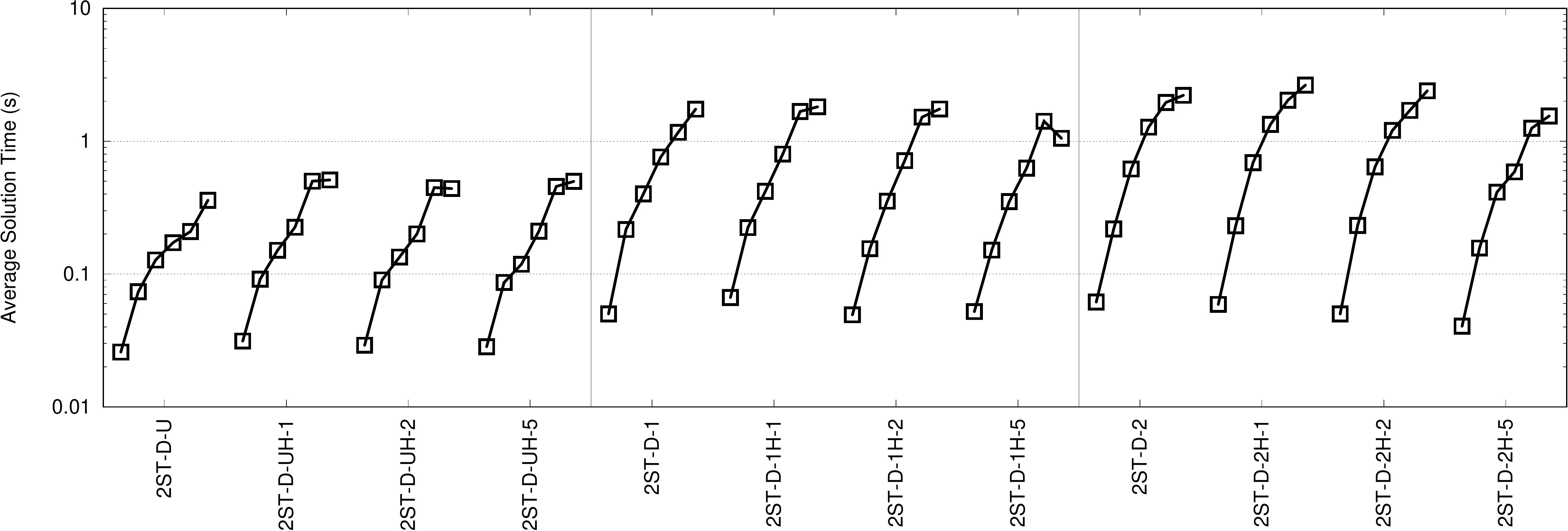}
	\caption{Two-stage problems with discrete uncertainty for varying values of $N$, Exp3. HIRO modifies $\pmb{C}$ and $\pmb{c}^j$.}
	\label{fig:ts-dis-exp3-Cc}
\end{figure}

\begin{figure}[htbp]
	\includegraphics[width=\textwidth]{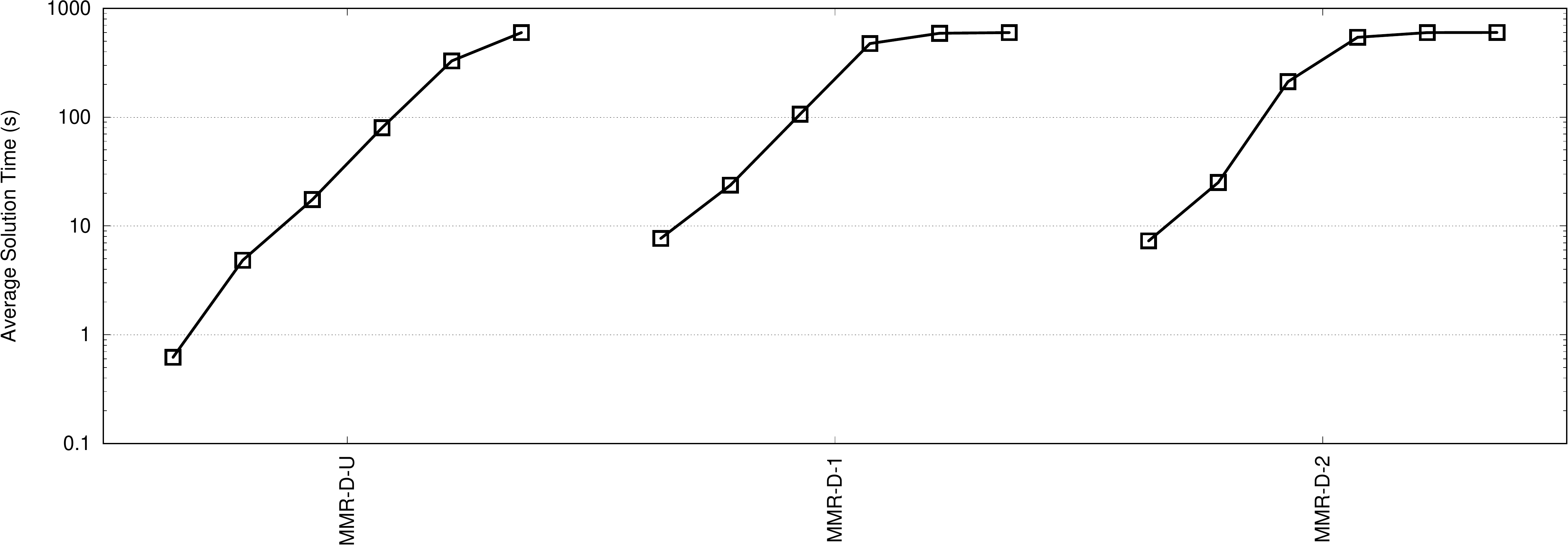}
	\caption{Two-stage problems with discrete uncertainty for varying values of $N$, Exp4.}
	\label{fig:ts-dis-exp4}
\end{figure}

In Figures~\ref{fig:ts-dis-exp3-C} and \ref{fig:ts-dis-exp3-Cc}, we show results on Exp3 and finally, we show the results for large values of $N$ (Exp4) in Figure~\ref{fig:ts-dis-exp4}. We note that in this setting, \TSTD{2} gives instances with the highest solution times, followed by \TSTD{1}.

\subsection{Discrete Budgeted Uncertainty}
\label{Two Stage-Discrete Budgeted Uncertainty}
	
\subsubsection{Problem Statement}
	
In this section we consider the two-stage selection problem under discrete budgeted uncertainty set. We use the linear formulation introduced in \cite{chassein2018recoverable} to assess the performance of different instance generators. Let 
\[ \cS = \{0\} \cup \{\underline{c}_i : \; i\in[n]\} \cup \{\overline{c}_i : \; i\in[n]\} =: \{\alpha^1 , \dots , \alpha^K\} \]
with $K = |\cS|$, the model is as follows.
\begin{align*}
\min\ & t\\
\text{s.t. } & t \geq \sum_{i\in[n]} C_i x_i + (p-\sum_{i\in[n]}x_i)\alpha^k - \sum_{i\in[n]}(1-x_i)[\alpha^k-\underline{c}_i]^+ + \Gamma \pi^k + \sum_{i\in[n]} \rho_i^k & k\in[K]\\
& \sum_{i\in[n]} x_i \leq p\\
& \pi^k + \rho^k_i \geq (1-x_i) ([\alpha^k-\underline{c}_i]^+ - [\alpha^k-\underline{c}_i-d_i]^+) & i\in[n] , k\in[K]\\
& x_i\in\{0,1\} & i\in[n]\\
& \pi^k \geq 0 & k\in[K]\\
& \rho^k_i \geq 0 & i\in[n] , k\in[K]
\end{align*}
	
\subsubsection{Sampling}\label{Two Stage-Discrete Budgeted Uncertainty-Sampling}
	
As a baseline sampling method, we use \TSTDB{U}, where all values $C_i$, $\underline{c}_i$ and $d_i$ are chosen uniformly from $\{1,\ldots,100\}$. In \TSTDB{1}, we differ from this setting by sampling small lower costs $\underline{c}_i$ from $\{1,\ldots,10\}$ and high deviations $d_i$ from $\{100-c_i,\ldots,100\}$. Finally, in \TSTDB{2}, we differ from \TSTDB{U} by setting $\underline{c}_i = 100 - C_i$ and sampling $d_i$ from $\{c_i,\ldots,100\}$.
	
\subsubsection{Experimental Setup}\label{Two Stage-Discrete Budgeted Uncertainty-Experiment}
	
In this experiment we consider the cases when $n=100$ and $p \in \{25,50,75\}$. We also use different $\Gamma$ for each value of $p$. If $p=25$, we choose $\Gamma\in\{5,10,15,20\}$;  if $p=50$, we choose $\Gamma\in\{10,20,30,40\}$; finally, we choose $\Gamma\in\{15,30,45,60\}$ for $p=75$. For all inputs we generate 50 instances using the sampling methods \TSTDB{U}, \TSTDB{1} and \TSTDB{2}. As a consequence, $3\times4\times50\times3=1800$ instances are generated. All instances are solved using CPLEX-20-10 with a 600-second time limit. Problem sizes were chosen so that results can be compared with those on recoverable robust optimization provided in Section~\ref{subsec:recoverable-discrete-budgeted}.
	
\subsubsection{Experimental Results}\label{Two Stage-Discrete Budgeted Uncertainty-Results}
	
\begin{figure}[htbp]
	\includegraphics[width=\textwidth]{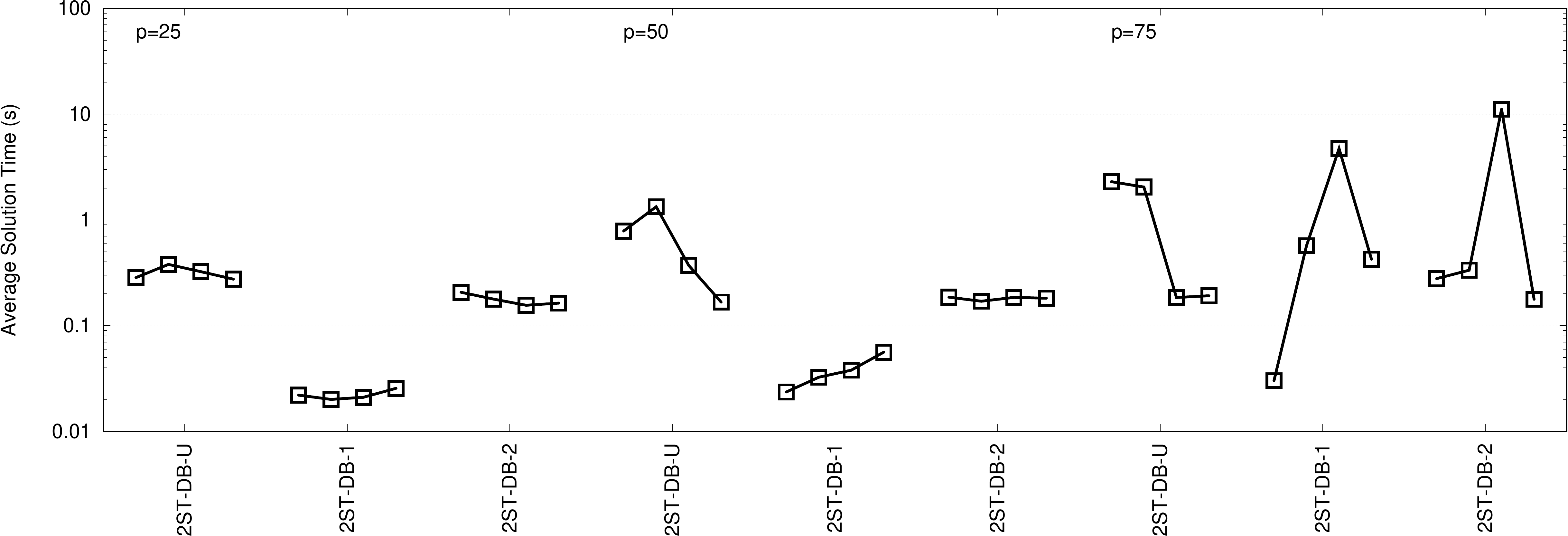}
	\caption{Two-stage problems with discrete budgeted uncertainty for varying $p$ (noted in top left corners) and $\Gamma$.}
	\label{fig:ts-dis-bud}
\end{figure}

The results shown in Figure~\ref{fig:ts-dis-bud} show that computation times are very small and even \TSTDB{1} and \TSTDB{2} fail to generate harder instances compared to \TSTDB{U} in most cases. Although the parameter $\Gamma$ have no considerable impact on the solution time when $p=25$, it can be effective for \TSTDB{U} when $p=50$; in particular, $\Gamma=20$ increases the solution time to over one second. The best results of \TSTDB{1} and \TSTDB{2} could be observed for $p=75$ and $\Gamma=30$ when the solution time passes 5 and 10 seconds, respectively.

\subsection{Continuous Budgeted Uncertainty}\label{Two Stage-Continuous Budgeted Uncertainty}
	
\subsubsection{Problem Statement}
	
In this section we consider the two-stage selection problem under continuous budgeted uncertainty $\cU_{Bvar}^\Gamma$. Using the formulation provided in \cite{chassein2018recoverable}, the problem can be formulated as follows.
\begin{align*}
\min\ & \sum_{i\in[n]} C_i x_i + \sum_{i\in[n]} \underline{c}_i y_i + \Gamma \pi + \sum_{i\in[n]} d_i \rho_i\\
\text{s.t. } & \sum_{i\in[n]} (y_i + x_i) = p\\
& x_i + y_i \leq 1 & \forall i\in[n]\\
& \pi + \rho_i \geq y_i & \forall i\in[n]\\
& \pi \geq 0 \\
& \rho_i \geq 0 & \forall i\in[n]\\
& x_i\in\{0,1\} & \forall i\in[n]\\
& y_i\in [0,1] & \forall i\in[n]
\end{align*}

\subsubsection{Sampling}\label{Two Stage-Continuous Budgeted Uncertainty-Sampling}
	
We use the same instance generators as in the case of discrete budgeted uncertainty, described in Section~\ref{Two Stage-Discrete Budgeted Uncertainty-Sampling}. These are \TSTCB{U} (uniform sampling), \TSTCB{1} (small lower bounds $\underline{\pmb{c}}$ and large deviations $\pmb{d}$), and \TSTCB{2} (lower bounds plus first-stage costs are constant).

\subsubsection{Experimental Setup}\label{Two Stage-Continuous Budgeted Uncertainty-Experiment}
	
We use a similar experiment as for discrete budgeted uncertainty, see Section~\ref{Two Stage-Discrete Budgeted Uncertainty-Experiment}. Thus, we consider selection problems with $n=100$ and $p\in\{25,50,75\}$. Note that $\Gamma$ needs to be chosen differently, as it now denotes the budget of deviation instead of the number of items that can deviate. Hence, we choose $\Gamma \in \{400,800,1000,1200\}$. For each parameter choice and sampling method, 50 instances are generated. As a result, $3\times4\times50\times3=1800$ instances are generated in total. We use CPLEX-20-10 with a 600-second time limit to solve each instance.
	
\subsubsection{Experimental Results}\label{Two Stage-Continuous Budgeted Uncertainty-Results}
	
\begin{figure}[htbp]
	\includegraphics[width=\textwidth]{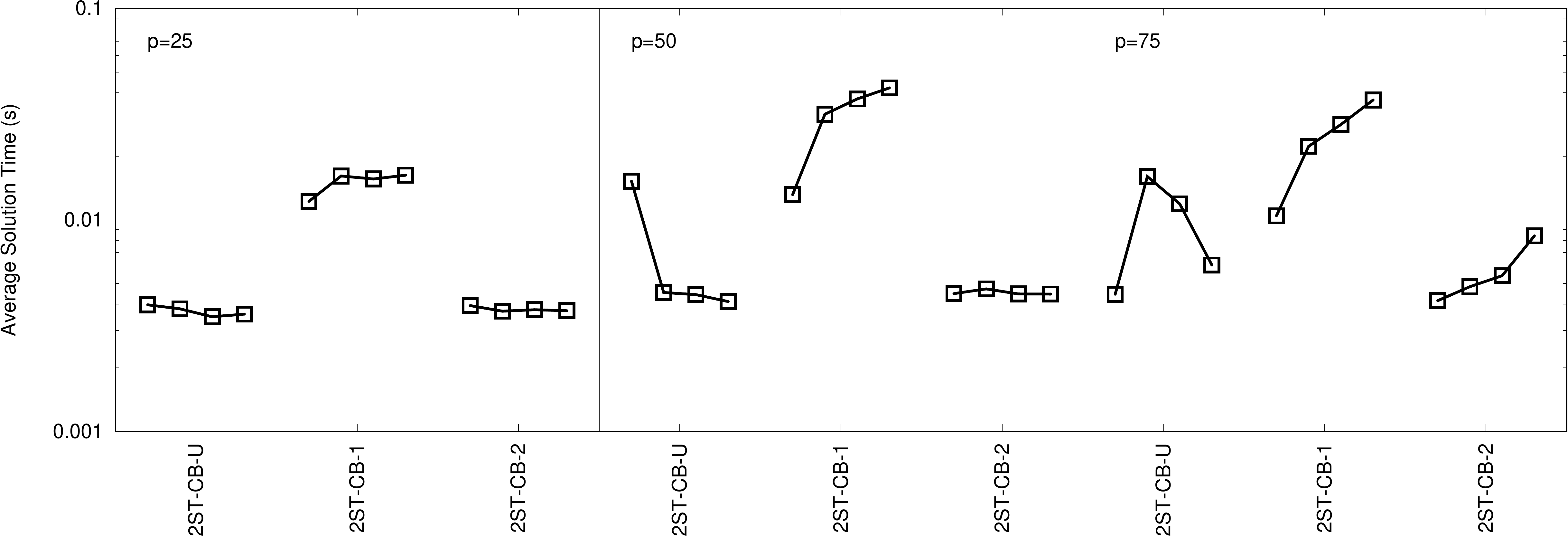}
	\caption{Two-stage problems with continuous budgeted uncertainty for varying $p$ (denoted in the top left corners) and $\Gamma$.}
	\label{fig:ts-con-bud}
\end{figure}

The results presented in Figure~\ref{fig:ts-con-bud} show that all three methods generate instances which can be solved fast (less that 0.1 second), but \TSTCB{1} stands out, as it can increase the solution time by almost a factor 10, compared to \TSTCB{U} when $p=50$ and $\Gamma=800,1000,1200$. 
	
\section{Recoverable Problems}\label{sec:recoverable}
	
\subsection{Discrete Uncertainty}\label{subsec:recoverable-discrete}
	
\subsubsection{Problem Statement}
	
We consider recoverable robust selection problems under discrete uncertainty sets, which we formulate using the model presented in \cite{chassein2018recoverable}:
\begin{align*}
\min\ & \sum_{i\in[n]} C_i x_i + t \\
\text{s.t. } & t \ge \sum_{i\in[n]} c^j_i y^j_i & \forall j\in[N] \\
& \sum_{i\in[n]} x_i = p \\
& \sum_{i\in[n]} y^j_i = p & \forall j\in[N] \\
& z^j_i \le x_i & \forall i\in[n], j\in[N] \\
& z^j_i \le y^j_i & \forall i\in[n], j\in[N] \\
& \sum_{i\in[n]} z^j_i \ge \Delta & \forall j\in[N] \\
& x_i \in\{0,1\} & \forall i\in[n] \\
& y^j_i \in\{0,1\} & \forall i\in[n], j\in[N] \\
& z^j_i \in\{0,1\} & \forall i\in[n], j\in[N] 
\end{align*}
Here, $\Delta$ is called the recovery factor and means that at least $\Delta$ items of the first-stage solution must remain in the second-stage solution. Variables $z^j_i$ measure if first-stage solution $\pmb{x}$ and second-stage solution $\pmb{y}^j$ both pack the same item $i\in[n]$.
	
\subsubsection{Sampling}\label{subsubsec:recoverable-discrete-sampling}

We use the same sampling methods as provided for two-stage problems in Section~\ref{subsubsec:twostage-discrete-sampling}, which we denote here as \RRD{U}, \RRD{1}, and \RRD{2}, respectively.

\subsubsection{HIRO}\label{subsubsec:recoverable-hiro}
	
To formulate the HIRO approach, we would like to optimize $\pmb{C}$ and $\pmb{c}^j$, $j\in[N]$, with respect to a set of candidate solution $\pmb{x}^k$, $k\in[K]$, as follows:
\begin{align*}
\max\ & t \\
\text{s.t. } & t \le \sum_{i\in[n]} C_i x^k_i + \max_{j\in[N]} Q(\pmb{x}^k,\pmb{c}^j) & \forall k\in[K] \\
&\pmb{c}^j \in \cU(\tilde{\pmb{c}}^j) & \forall j\in[N] \\
&\pmb{C} \in \cU(\tilde{\pmb{C}}) 
\end{align*}
where, as before, $Q(\pmb{x}^k,\pmb{c}^j)$ denotes the second-stage costs of $\pmb{x}^k$ under scenario $\pmb{c}^j$. To determine these costs, we need to solve the second-stage recovery problem, which is given as
\begin{align*}
\min\ &\sum_{i\in[n]} c^j_i y^k_i \\
\text{s.t. } & \sum_{i\in[n]} y^k_i  = p\\
& \sum_{i\in[n]} x^k_i y^k_i  \geq p - \Delta \\
& x^k_i + y^k_i \le 1 & \forall  i\in[n] \\
& y^k_i \in\{0,1\}& \forall i\in[n]
\end{align*}
Note that there always exists an integer optimal solution for the LP-relaxation of this problem. Hence, using duality and variables $\pmb{\lambda}$ that decide which scenario is assigned to each candidate solution, the HIRO problem is formulated as follows.
\begin{align*}
\max\ & t \\
\text{s.t. } & t \le \sum_{i\in[n]} C_i x^k_i + p \beta^k + (p - \Delta) \eta^k - \sum_{i\in[n]} \gamma^k_i & \forall k\in[K] \\
& \sum_{j\in[N]} \lambda^k_j = 1 & \forall k\in[K]\\
& \beta^k + \eta^k x^k_i \leq \gamma^k_i + \sum_{j\in[N]} c^j_i \lambda^k_j & \forall i\in[n],k\in[K] \\
& \lambda^k_j \in\{0,1\}  & \forall j\in[N],k\in[K]\\
& \beta^k, \eta^k \geq 0 & \forall k\in[K]\\
& \gamma^k_i \geq 0 & \forall i\in[n],k\in[K]\\
&\pmb{c}^j \in \cU(\tilde{\pmb{c}}^j) & \forall j\in[N] \\
&\pmb{C} \in \cU(\tilde{\pmb{C}}) 
\end{align*}
If we only optimize over $\pmb{C}$, we can use nominal values $\tilde{\pmb{c}}^j$ for second-stage costs instead of treating them as variables. However, if we optimize over both $\pmb{C}$ and $\pmb{c}^j$, the nonlinearity in $c^j_i \lambda^k_i$ needs to be replaced with additional variables.

\subsubsection{Experimental Setup}\label{subsubsec:recoverable-experiments}
	
This setup is similar to the one introduced in for two-stage problems in Section~\ref{subsubsec:twostage-discrete-experiments}. We consider four experiments to assess the effect of different parameters. 
In Exp1, we consider the cases when $N=n=50$, $p=25$, $\Delta\in\{13,20\}$ and $N=n=100$, $p=50$, $\Delta\in\{25,40\}$. 
In Exp2, we fix $n=N=50$, but change $p$ and $\Delta$, using $p=25$, $\Delta\in\{13,20\}$, $p=30$, $\Delta\in\{15,25\}$, and $p=40$, $\Delta\in\{20,30\}$. 
In Exp3, we fix $n=50$ and $p=25$, but change $N\in\{40,50,60\}$ and $\Delta\in\{13,20\}$. 
Finally, in Exp4, we use the parameters of Exp3, but with $N\in\{100,200,500,1000,2000\}$.
	
For each parameter choice and sampling method, we generate 50 instances. For Exp1, Exp2 and Exp3, we also apply the HIRO method to each instance, where we vary only the first-stage costs, or both first-stage and second-stage costs, using varying  budgets $b\in\{1,2,5\}$ and a 600-second generation time limit. Thus, Exp1 contains $2\times2\times50\times3\times7=4200$ instances, Exp2 contains $5\times2\times50\times3\times7=10500$ instances, Exp3 contains $6\times2\times50\times3\times7=12600$ instances, and there are$6\times2\times50\times3=1800$ instances in Exp4. We solve all generated instances using CPLEX-20-10 with a 600-second time limit.
	
\subsubsection{Experimental Results}\label{subsubsec:recoverable-discrete-results}
	
\begin{figure}[htbp]
	\includegraphics[width=\textwidth]{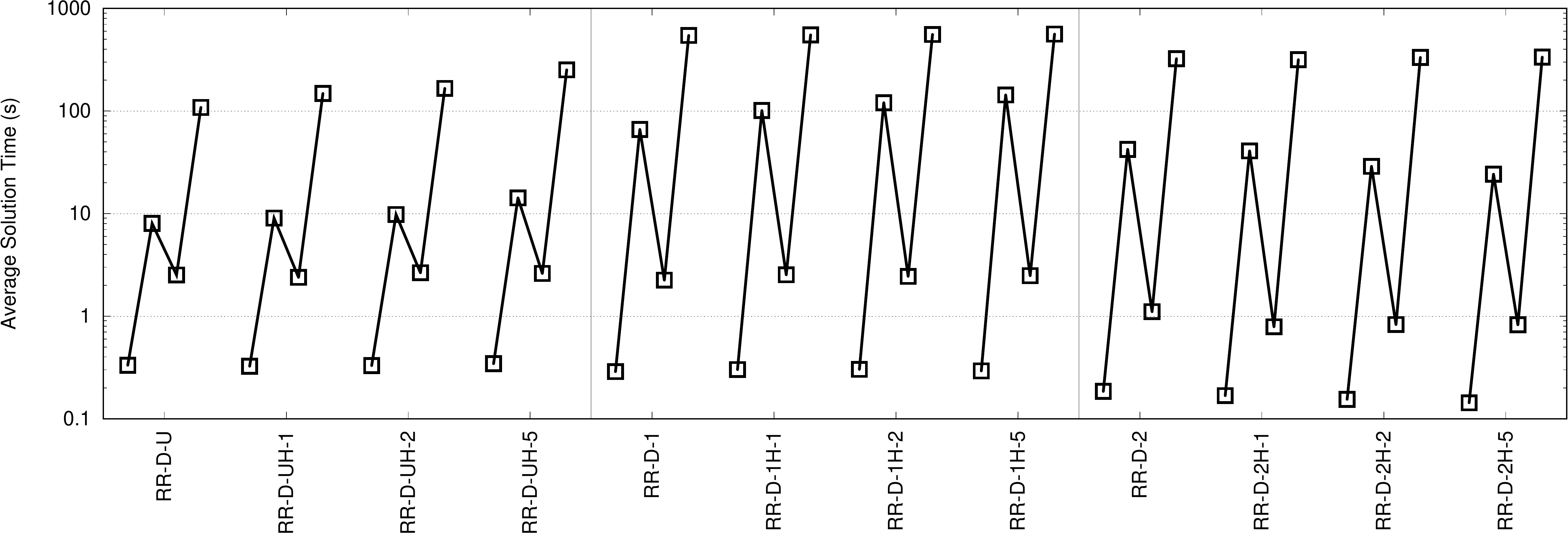}
	\caption{Recoverable robust problems with discrete uncertainty for varying values of $(N,n,p,\Delta)$, Exp1. HIRO modifies $\pmb{C}$.}
	\label{fig:rec-dis-exp1-C}
\end{figure}

\begin{figure}[htbp]
	\includegraphics[width=\textwidth]{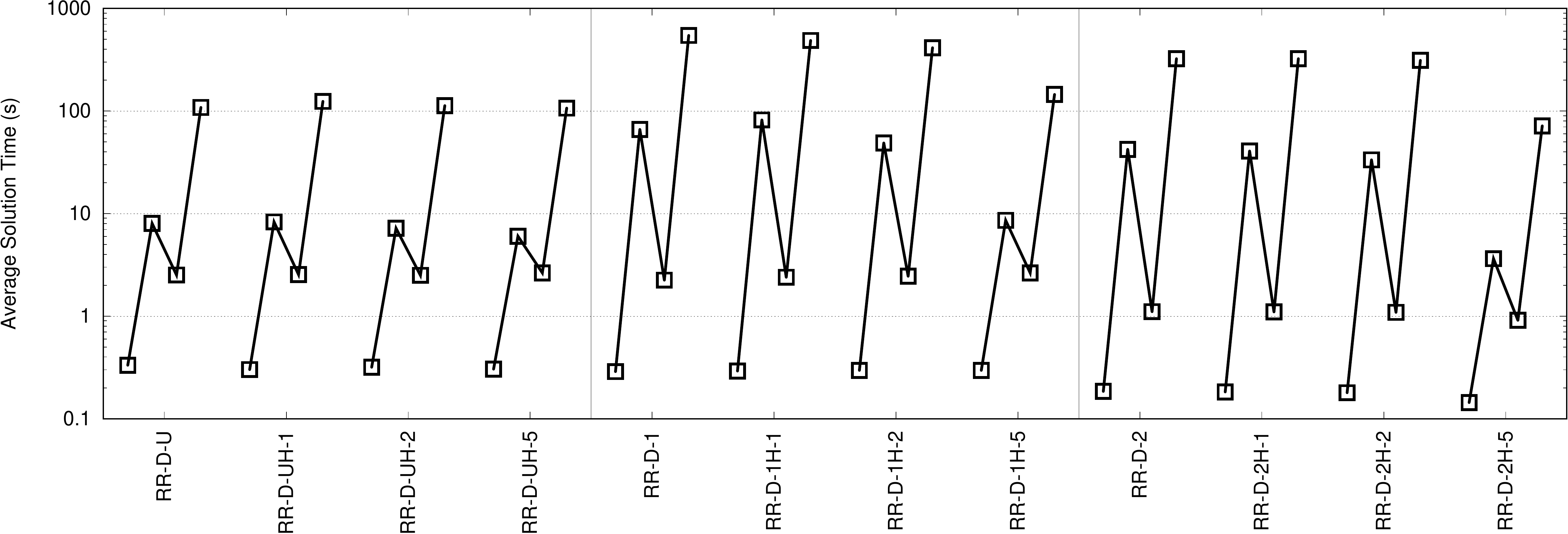}
	\caption{Recoverable robust problems with discrete uncertainty for varying values of $(N,n,p,\Delta)$, Exp1. HIRO modifies $\pmb{C}$ and $\pmb{c}^j$.}
	\label{fig:rec-dis-exp1-Cc}
\end{figure}

Figures~\ref{fig:rec-dis-exp1-C} and \ref{fig:rec-dis-exp1-Cc} show results for Exp1. While HIRO methods can increase solution times for most instances generated by \RRD{U} or \RRD{1} when we only modify first-stage costs $\pmb{C}$, results are mixed in the other cases. Both \RRD{1} and \RRD{2} generate instances that are considerably harder to solve than those generated by \RRD{U}, in particular for higher values of $\Delta$.
	
\begin{figure}[htbp]
	\includegraphics[width=\textwidth]{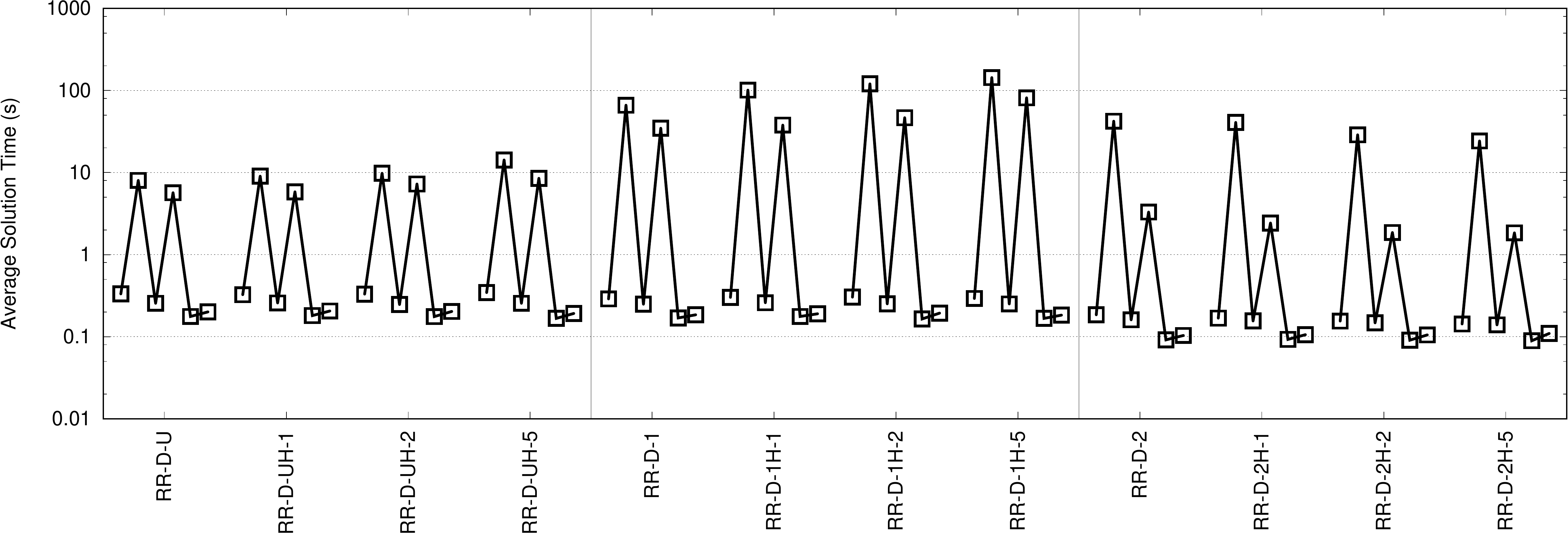}
	\caption{Recoverable robust problems with discrete uncertainty for varying values of $(p,\Delta)$, Exp2. HIRO modifies $\pmb{C}$.}
	\label{fig:rec-dis-exp2-C}
\end{figure}

\begin{figure}[htbp]
	\includegraphics[width=\textwidth]{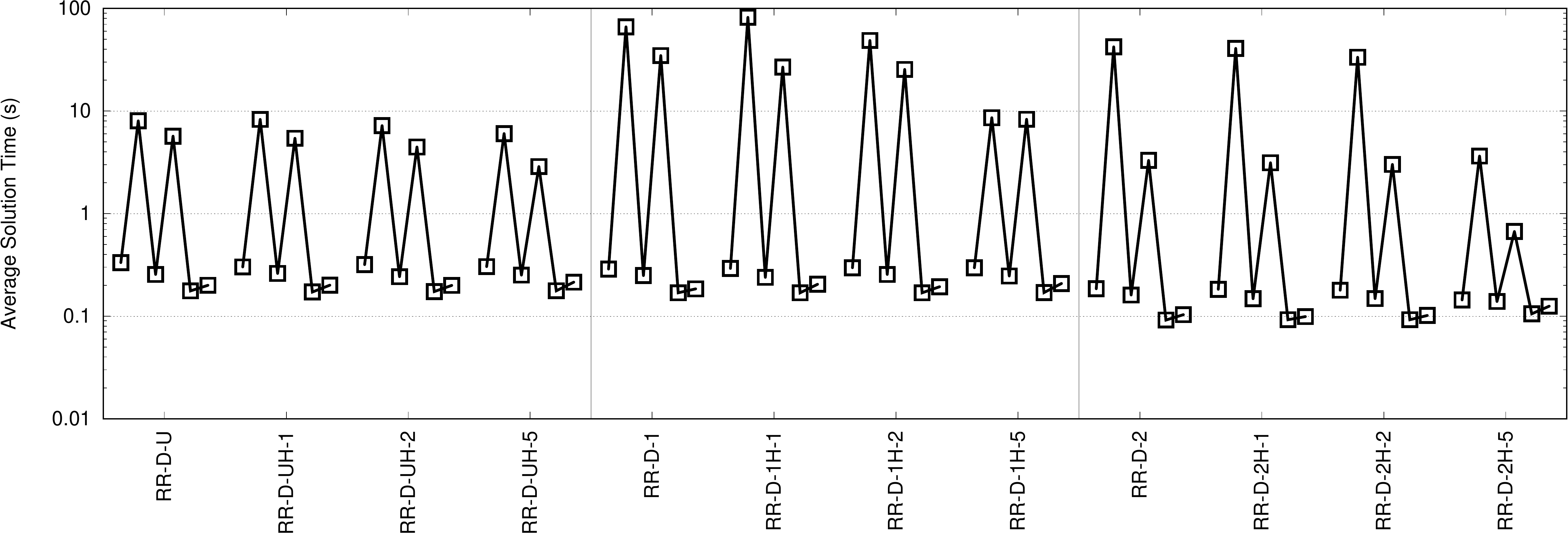}
	\caption{Recoverable robust problems with discrete uncertainty for varying values of $(p,\Delta)$, Exp2. HIRO modifies $\pmb{C}$ and $\pmb{c}^j$.}
	\label{fig:rec-dis-exp2-Cc}
\end{figure}

In Figures~\ref{fig:rec-dis-exp2-C} and \ref{fig:rec-dis-exp2-Cc}, we present results for Exp2. The HIRO approach can be effective to increase solution times, if only first-stage costs $\pmb{C}$ are modified. As in Exp1, a positive correlation between increasing $\Delta$ and solution time can be observed. While for $p=40$, all solution times remain small, increasing $\Delta$ gives a large increase in solution time for $p=25$ or $p=30$. In particular \RRD{1} gives consistently good results.
	
\begin{figure}[htbp]
	\includegraphics[width=\textwidth]{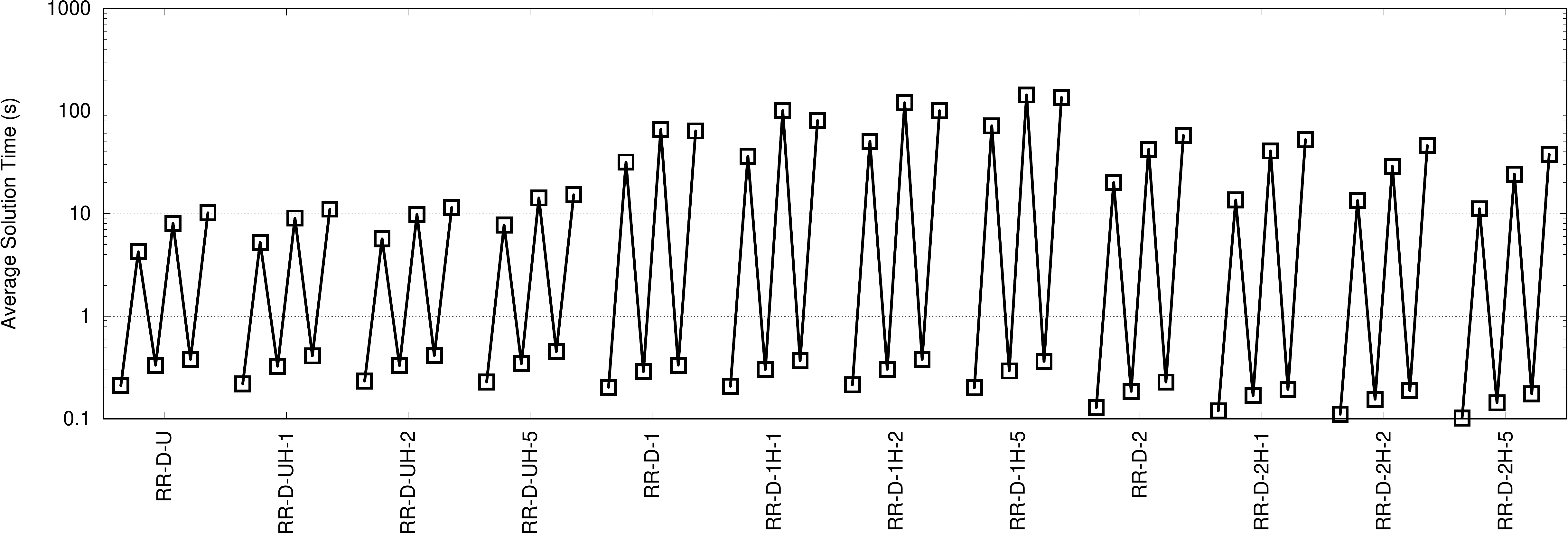}
	\caption{Recoverable robust problems with discrete uncertainty for varying values of $(N,\Delta)$, Exp3. HIRO modifies $\pmb{C}$.}
	\label{fig:rec-dis-exp3-C}
\end{figure}

\begin{figure}[htbp]
	\includegraphics[width=\textwidth]{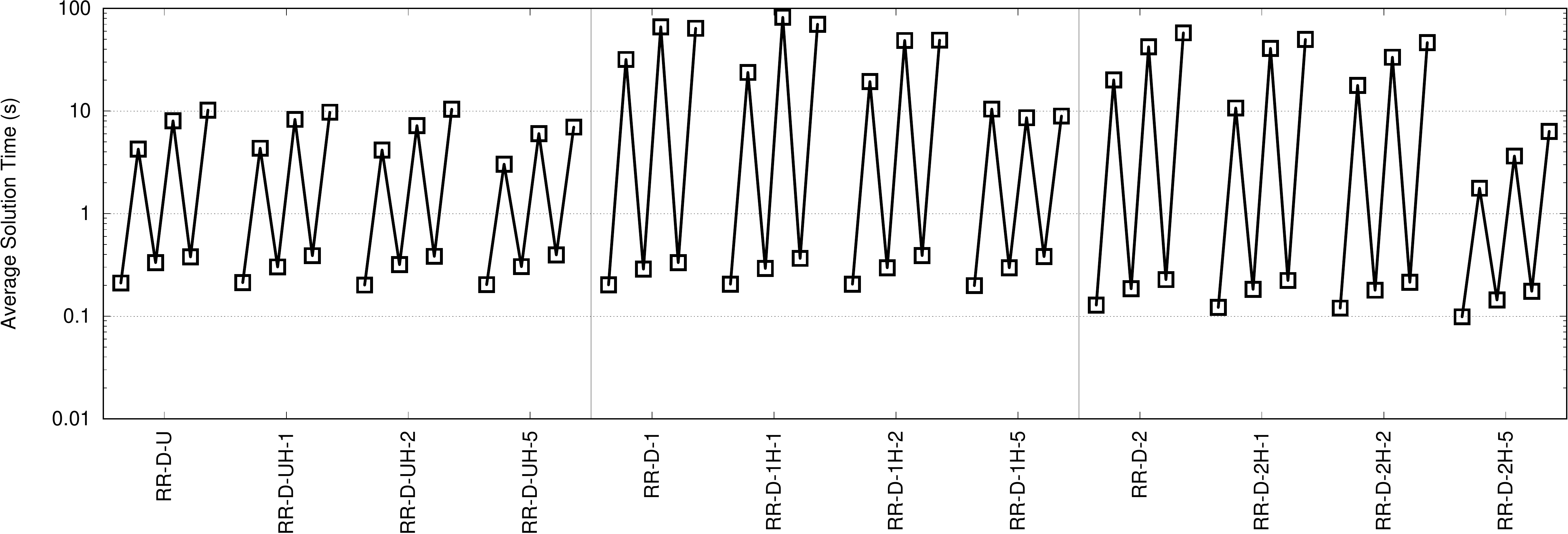}
	\caption{Recoverable robust problems with discrete uncertainty for varying values of $(N,\Delta)$, Exp3. HIRO modifies $\pmb{C}$ and $\pmb{c}^j$.}
	\label{fig:rec-dis-exp3-Cc}
\end{figure}

In Figures~\ref{fig:rec-dis-exp3-C} and \ref{fig:rec-dis-exp3-Cc}, we summarize results on Exp3. Sampling methods without HIRO become slightly harder with increasing $N$, and considerably harder when using $\Delta=20$ instead of $\Delta=13$. Both \RRD{1} and \RRD{2} perform similarly and produce instances that are harder to solve than when using \RRD{U}. For these instances, applying the HIRO method can have a large impact: for $\Delta=13$ and $N=50$, solution times increase by several orders of magnitude. This effect is not apparent for $N=40$ or $N=60$. For $\Delta=20$ and $N=50$, this effect is reversed, and solution times drop by orders of magnitude (while they can be slightly increased for other values of $N$ when only first-stage costs are modified).

\begin{figure}[htbp]
	\includegraphics[width=\textwidth]{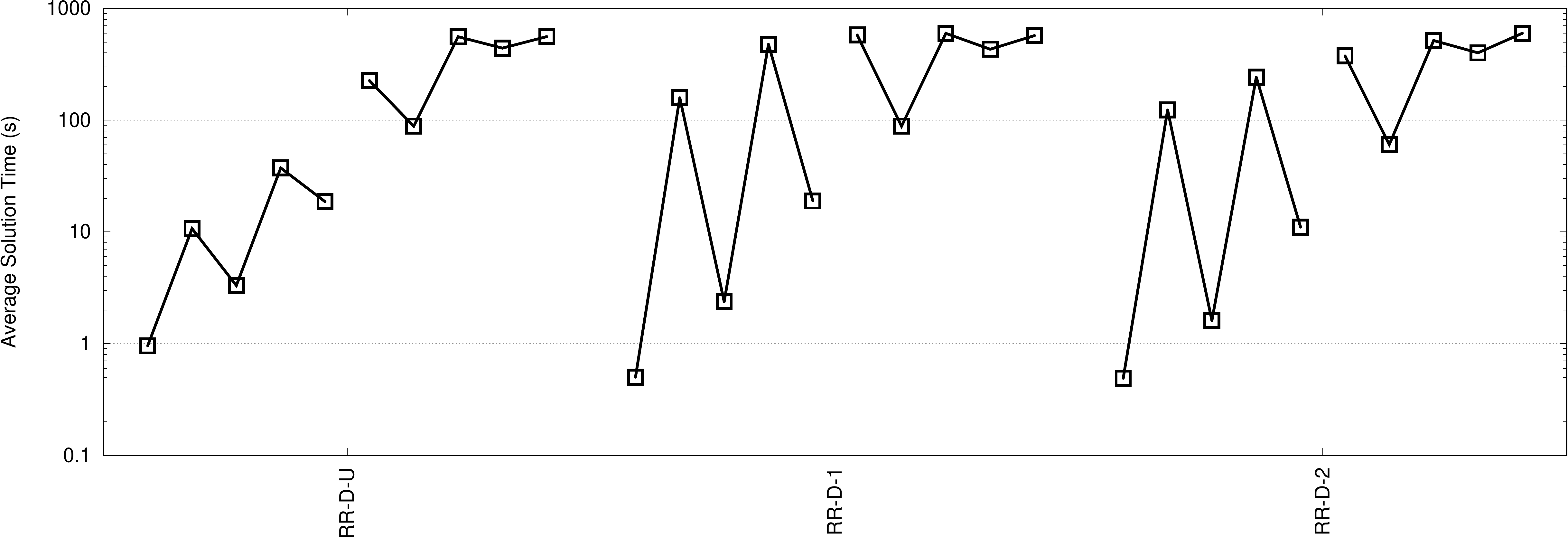}
	\caption{Recoverable robust problems with discrete uncertainty for varying values of $(N,\Delta)$, Exp4.}
	\label{fig:rec-dis-exp4}
\end{figure}

Finally, we consider large values of $N$ for Exp4 in Figure~\ref{fig:rec-dis-exp4}. Especially for $\Delta=20$, sampling methods \RRD{1} and \RRD{2} produce instances that are harder to solve than those generated by \RRD{U}. For $\Delta=13$, these differences are less pronounced. For large values of $N$, all solution times are at the time limit, which means that they cannot be further differentiated.

\subsection{Discrete Budgeted Uncertainty}\label{subsec:recoverable-discrete-budgeted}
	
\subsubsection{Problem Statement}
	
We now discuss recoverable selection problems under discrete budgeted uncertainty. To model these problems, we apply the following formulation adapted from \cite{chassein2018recoverable}. Let
\begin{align*}
\mathcal{S} &= \left\{ (\alpha, \beta) : \exists i,j \in[n], \delta_i,\delta_j\in\{0,1\} \text{ s.t. } \alpha = \underline{c}_i + \delta_id_i, \beta = [\underline{c}_j - d_j\delta_j - \alpha]_+ \right\} \\
&=: \{(\alpha^1,\beta^1),\ldots,(\alpha^K,\beta^K)\}
\end{align*}
with $K=|\cS|$. Then, a compact formulation for the robust problem is as follows.
\begin{align*}
\min\ & t\\
\text{s.t. } & t \geq \sum_{i\in[n]} C_i x_i + \Gamma \pi^k + \sum_{i\in[n]} \rho_i^k + p\alpha^k + (p-\Delta)\beta^k - \sum_{i\in[n]}[\alpha^k+x_i\beta^k-\underline{c}_i]_+ & \forall k\in[K]\\
& \pi^k + \rho^k_i \geq [\alpha^k+x_i\beta^k-\underline{c}_i]_+ - [\alpha^k+x_i\beta^k-\underline{c}_i-d_i]_+ & \forall i\in[n] , k\in[K]\\
& \sum_{i\in[n]} x_i = p\\
& x_i\in\{0,1\} & \forall i\in[n]\\
& \rho^k_i \geq 0 & \forall i\in[n] , k\in[K]\\
& \pi^k \geq 0 & \forall k\in[K]
\end{align*}
	
\subsubsection{Sampling}\label{subsubsec:recoverable-discrete-budgeted-sampling}
	
We use the same sampling methods as for two-stage robust selection with discrete budgeted uncertainty, see Section~\ref{Two Stage-Discrete Budgeted Uncertainty-Sampling}. We refer to them as \RRDB{U} for uniform sampling, \RRDB{1} for sampling with small costs $\underline{c}_i$ and high deviations $d_i$, and \RRDB{2} for the setting where $\underline{c}_i = 100 - C_i$.

\subsubsection{Experimental Setup}\label{subsubsec:recoverable-discrete-budgeted-experiment}
	
We consider problems where $n=100$ and $p \in \{25,50,75\}$. We also use different $\Gamma$ for each value of $p$. If $p=25$, we choose $\Gamma\in\{5,10,15,20\}$;  if $p=50$, we choose $\Gamma\in\{10,20,30,40\}$; finally, we choose $\Gamma\in\{15,30,45,60\}$ for $p=75$. Unlike the experiment presented in Section~\ref{Two Stage-Discrete Budgeted Uncertainty-Experiment}, we also need to consider parameter $\Delta$ which depends on the value of $p$, too. Here, we consider the similar value of $\Gamma$ for $\Delta$, so if $p=25$, we choose $\Delta\in\{5,10,15,20\}$;  if $p=50$, we choose $\Delta\in\{10,20,30,40\}$; finally, we choose $\Delta\in\{15,30,45,60\}$ for $p=75$.
We generate 50 instances using the three sampling methods for each parameter choice. As a consequence, we generate $3\times4\times4\times50\times3=7200$ instances for this experiment. All instances are solved using CPLEX-20-10 with a 600-second time limit.

\subsubsection{Experimental Results}\label{subsubsec:recoverable-discrete-budgeted-results}
	
\begin{figure}[htbp]
	\includegraphics[width=\textwidth]{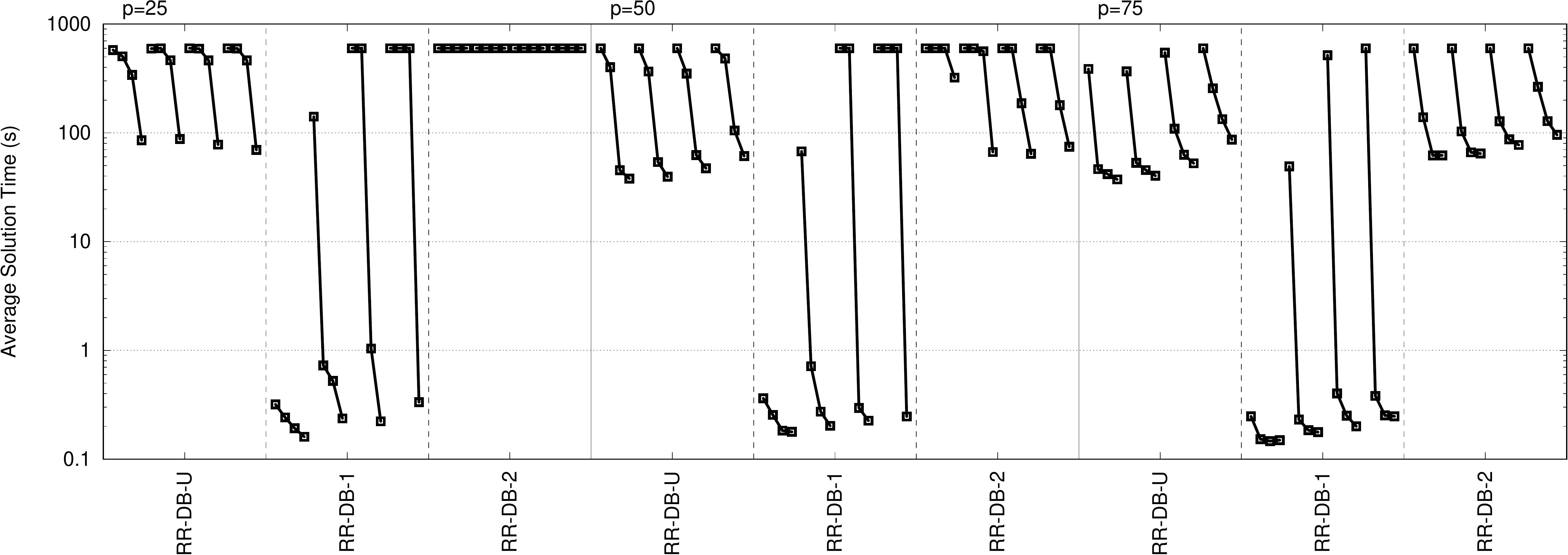}
	\caption{Recoverable robust problems with discrete budgeted uncertainty for varying $p$ (noted in top left corners) and $(\Gamma,\Delta)$.}
	\label{fig:rec-dis-bud}
\end{figure}

We present results in Figure~\ref{fig:rec-dis-bud}, where values for each $\Gamma$ are grouped together and within each group, the four values for different choices of $\Delta$ are shown. Method \RRDB{2} only slightly outperforms \RRDB{U}, while instances generated by \RRDB{1} can be amongst the easiest to solve. While for $p=25$, all parameter choices of $\Gamma$ and $\Delta$ result in problems of similar difficulty, only $\Gamma=10$ or $\Gamma=20$ results in hard problems for $p=50$, while only $\Gamma=15$ gives hard problems for $p=75$.
		
\subsection{Continuous Budgeted Uncertainty}\label{subsec:recoverable-continuous-budgeted}

\subsubsection{Problem Statement}
	
A compact formulation for recoverable robust selection problems with continuous budgeted uncertainty $\cU_{Bvar}^\Gamma$ was presented in \cite{chassein2018recoverable} as follows.
\begin{align*}
\min\ & \sum_{i\in[n]} C_i x_i + \sum_{i\in[n]} \underline{c}_i y_i + \Gamma \pi + \sum_{i\in[n]} d_i \rho_i\\
\text{s.t. } & \sum_{i\in[n]} y_i = p\\
& \sum_{i\in[n]} x_i = p\\
& \sum_{i\in[n]} z_i \geq p-\Delta\\
& z_i \leq x_i & \forall i\in[n]\\
& z_i \leq y_i & \forall i\in[n]\\
& \pi + \rho_i \geq y_i & \forall i\in[n]\\
& x_i\in\{0,1\} & \forall i\in[n]\\
& y_i\in [0,1] & \forall i\in[n]\\
& \pi \geq 0 \\
& \rho_i \geq 0 & \forall i\in[n]\\
\end{align*}
	
\subsubsection{Sampling}\label{subsubsec:recoverable-continuous-budgeted-sampling}

We use the same sampling methods as provided for two-stage problems in Section~\ref{Two Stage-Discrete Budgeted Uncertainty-Sampling}, which we denote here as \RRDB{U}, \RRDB{1}, and \RRDB{2}, respectively.

\subsubsection{Experimental Setup}\label{subsubsec:recoverable-continuous-budgeted-experiment}

We use similar parameters as in Section~\ref{subsubsec:recoverable-discrete-budgeted-experiment}, except for $\Gamma$. We set $n=100$ and consider $p\in\{25,50,75\}$ and $\Gamma \in 
\{400,800,1000,1200\}$. We consider $\Delta\in\{5,10,15,20\}$ when $p=25$, $\Delta\in\{10,20,30,40\}$ when $p=50$ and $\Delta\in\{15,30,45,60\}$ when $p=75$. For each parameter choice and sampling method, 50 instances are generated (a total of $3\times4\times4\times50\times3=7200$ instances) and solved using CPLEX-20-10 with a 600-second time limit.
	
\subsubsection{Experimental Results}\label{subsubsec:recoverable-continuous-budgeted-results}
	
\begin{figure}[htbp]
	\includegraphics[width=\textwidth]{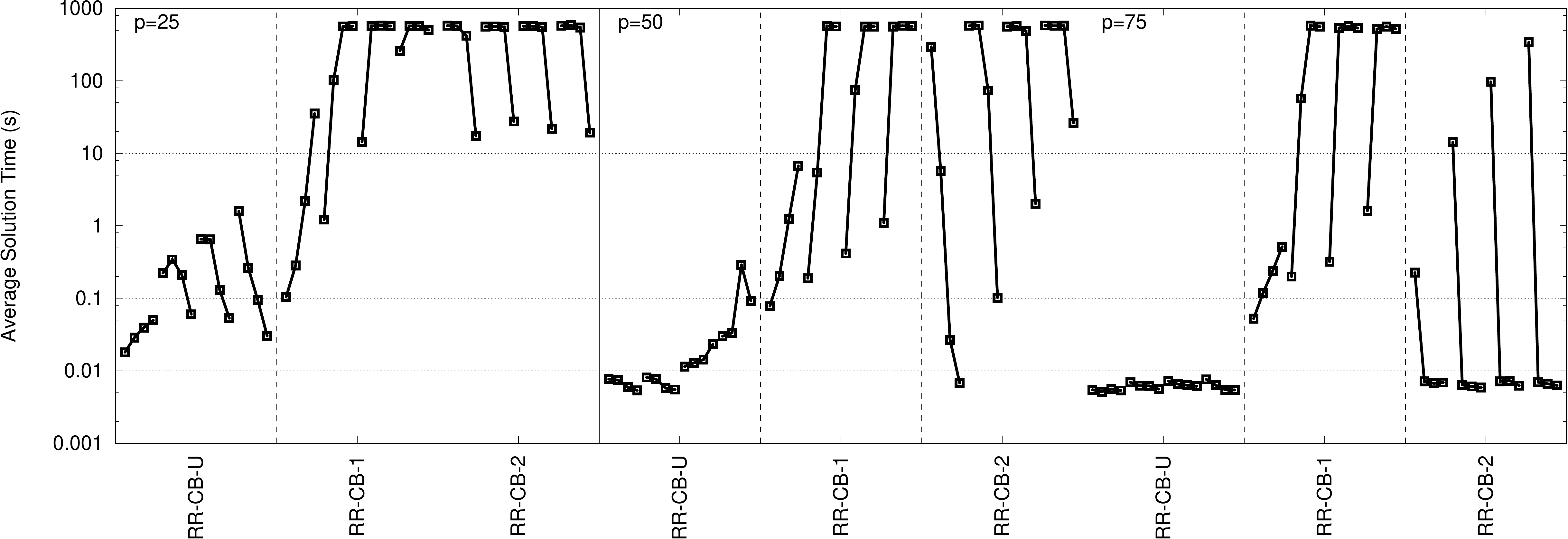}
	\caption{Recoverable robust problems with discrete budgeted uncertainty for varying $p$ (noted in top left corners) and $(\Gamma,\Delta)$.}
	\label{fig:rec-con-bud}
\end{figure}

We present results in Figure~\ref{fig:rec-con-bud}, where as before values for each $\Gamma$ are grouped together and within each group, the four values for different choices of $\Delta$ are shown.
Methods \RRDB{1} and \RRDB{2} generate instances that are harder to solve in almost all cases compared to the \RRDB{U} method. In all cases, solution times tend to increase with $\Gamma$. For $\Delta$, behavior is very different, with an increase with $\Delta$ in some cases, and a decrease in others. Furthermore, instances tend to become easier for increasing $p$. 
	
\section{Instance Library and Formats}
\label{sec:library}
	
An outcome of this paper is to create a benchmark library that collects all instances presented in this paper, along with their generating code, in order to drive further development of improved solution algorithms. The library is available at a dedicated website https://robust-optimization.com. The instances included so far will form a basic stock of problems to be extended further.
	
For ease of usability, we decided to use a format that is as simple as possible. All instances are given as comma separated values (csv). The first line of each instance represents the main parameters depending on its uncertainty set and criterion. The instances generated for the min-max and min-max regret criteria under discrete uncertainty set have the same structure. Their first line consists of $n$, $p$ and $N$, and the following $N$ lines represents each scenario. Instances generated for the two-stage and recoverable criteria under discrete uncertainty set have also similar structure. The difference is that the second line represent the first stage scenario and the the following $N$ line show each scenario. In addition, for the recoverable criterion the first line contains the forth parameter $\Delta$.
	
The interval uncertainty set is only used for the min-max regret criterion. The instances generated for the related section has only three lines, the first one provides $n$ and $p$, the second line shows the lower bound of the interval assigned to each item and the third line represents the deviation of the cited intervals. The same structure also is used for the min-max criterion under budgeted uncertainty set. The only difference is that the first line contains the third parameter $\Gamma$.
	
The two-stage and recoverable criteria under both budgeted uncertainty set have similar structure to the budgeted uncertainty for the min-max criterion. The difference is that they have four lines, the second line represents the first stage scenario and the third and fourth line shows the lower bound and deviation of each interval, respectively. Moreover, for the recoverable criterion, the first line has a fourth parameter $\Delta$. 
	
\section{Conclusions}
\label{sec:conclusions}
	
Having benchmark sets of instances available is of central importance for the development of improved solution algorithms. For robust discrete optimization, it seems that no such benchmark currently exists, despite an active research community developing specialized solution methods. In this paper, we considered a wide range of robust decision criteria (min-max, min-max regret, two-stage, recoverable robust) in combination with different uncertainty sets (interval, discrete, continuous and discrete budgeted) to generate instances for selection problems that are harder to solve than the baseline method of sampling parameters uniformly. Besides alternative sampling methods, we proposed optimization models that modify uncertainty sets to increase robust objective values with the aim to find instances that are even harder to solve. For all problem combinations we were able to generate instances that are harder to solve than those generated by the baseline method, sometimes by several orders of magnitude.
	
This work can be considered as a starting point for an ongoing community effort to standardize benchmarking of robust optimization algorithms. While the sampling methods proposed here can be applied to any other combinatorial optimization problem on principle, it is likely that more problem-tailored sampling methods will be more effective (e.g., profits that are correlated to weights for knapsack problems; or graphs that need to be generated for network problems). On a dedicated webpage (robust-optimization.com), we started collecting instances, including those generated in this paper, and cordially invite interested researchers to participate.
	
In terms of further research, it would be interesting to develop solution algorithms that make use of properties of hard instances to alleviate this hardness. For example, some difficult min-max regret instances with interval uncertainty we generated have either small $\underline{c}_i$ and $d_i$, or both large $\underline{c}_i$ and $d_i$. Along the paradigm of algorithm engineering, methods that can use such structure have the potential to perform better than the results produced by the general mixed-integer programming solver we used.

\appendix
	
\section{Overview of Sampling Generators}
\label{app:overview}
	
\subsection{Min-Max Problems with Discrete Uncertainty}

\begin{itemize}
	\item  \MMD{U}: For all $i\in[n]$ and $j\in[N]$, we choose $c_i^j\in\{1,\ldots,100\}$ iid uniformly.
		
	\item  \MMD{1}: For all $i\in[n]$ and $j\in[N]$, with probability $\frac{1}{2}$ we choose $c_i^j\in\{1,\ldots,10\}$, and with probability $\frac{1}{2}$ we choose $c_{ij}\in\{91,\ldots,100\}$ iid uniformly.
		
	\item  \MMD{2}: For each $i\in[n]$ and $j\in[N]$, if $i\leq\lfloor\frac{n}{2}\rfloor$, we choose $c_i^j\in\{1,\ldots,100\}$ iid uniformly; otherwise, we set $c_i^j = 100 - c_{i-\lfloor n/2\rfloor}^j$.
\end{itemize}

\subsection{Min-Max Problems with Budgeted Uncertainty}
	
\begin{itemize}
	\item \MMB{U}: For all $i\in[n]$, we choose both $\underline{c}_i , d_i$ uniformly as an integer number from $\{1,\ldots,100\}$.
		
	\item \MMB{1}: For all $i\in[n]$, we choose $\underline{c}_i$ uniformly as an integer number from $\{1,\ldots,100\}$ and set $d_i = 100 - \underline{c}_i$.
		
	\item \MMB{2}: For all $i\in[n]$, we choose $\underline{c}_i$ uniformly as an integer number from $\{1,\ldots,10\}$ and then choose $d_i \in \{100-\underline{c}_i-1,\ldots,100\}$.
\end{itemize}
	
\subsection{Min-Max Regret Problems with Interval Uncertainty}
	
\begin{itemize}
	\item \MMRI{U}: For all $i\in[n]$, we choose both $\underline{c}_i , d_i$ uniformly as an integer number from $\{1,\ldots,100\}$.
	\item \MMRI{1}: For all $i\in[n]$, with probability $\frac{1}{2}$, we choose $\underline{c}_i$ from $\{1,\ldots,10\}$ and $d_i$ from $\{91,\ldots,100\}$ uniformly. Otherwise, we choose $\underline{c}_i$ from $\{91,\ldots,100\}$ and $d_i$ from $\{1,\ldots,10\}$.
	\item \MMRI{2}: For all $i\in[n]$, with probability $\frac{1}{2}$, we choose $\underline{c}_i$ and $d_i$ from $\{1,\ldots,10\}$ uniformly. Otherwise, we choose $\underline{c}_i$ and $d_i$ from $\{91,\ldots,100\}$ uniformly.
\end{itemize}

\subsection{Min-Max Regret Problems with Discrete Uncertainty}

\begin{itemize}
	\item  \MMRD{U}: For all $i\in[n]$ and $j\in[N]$, we choose $c_i^j\in\{1,\ldots,100\}$ iid uniformly (=\MMD{U}).
		
	\item  \MMRD{1}: For all $i\in[n]$ and $j\in[N]$, with probability $\frac{1}{2}$ we choose $c_i^j\in\{1,\ldots,10\}$, and with probability $\frac{1}{2}$ we choose $c_{ij}\in\{91,\ldots,100\}$ iid uniformly (=\MMD{1}).
		
	\item  \MMRD{2}: For each $i\in[n]$ and $j\in[N]$, if $i\leq\lfloor\frac{n}{2}\rfloor$, we choose $c_i^j\in\{1,\ldots,100\}$ iid uniformly; otherwise, we set $c_i^j = 100 - c_{i-\lfloor n/2\rfloor}^j$ (=\MMD{2}).
\end{itemize}

\subsection{Two-Stage Problems with Discrete Uncertainty}
	
\begin{itemize}
	\item \TSTD{U}: For all $i\in[n]$, we choose $C_i , c_i^j$ iid uniformly from $\{1,\ldots,100\}$.
		
	\item \TSTD{1}: For all $i\in[n]$, with probability equal to $0.5$, we choose $C_i$ from $\{45,\ldots,55\}$ and with probability equal to $0.5$, from $\{25,\ldots,75\}$. 
	For all $i\in[n]$ and $j\in[N]$, with probability equal to $0.5$ we choose $c_i^j$ from $\{C_i - 5,\ldots,C_i + 5\}$, with probability equal to $0.25$ we choose $c_i^j$ from $\{1,\ldots,10\}$ and with probability equal to $0.25$, $c_i^j$ from $\{91,\ldots,100\}$.
		
	\item \TSTD{2}: For all $i\in[n]$, with probability equal to $0.5$ we choose $C_i$ from $\{1,\ldots,100\}$ and with probability equal to $0.5$ we set $C_i = 50$. Then, for all $i\in[n]$ and $j\in[N]$, if $C_i = 50$ with probability equal to $0.5$ we choose $c_i^j$ from $\{1,\ldots,10\}$ and with probability equal to $0.5$ we choose $c_i^j$ from $\{91,\ldots,100\}$; else, we choose $c_i^j$ from $\{C_i - 5,\ldots,C_i + 5\}$ and set negative values to be equal to zero.
\end{itemize}

\subsection{Two-Stage Problems with Discrete Budgeted Uncertainty}
\begin{itemize}
	\item \TSTDB{U}: For all $i\in[n]$, we choose $C_i, \underline{c}_i, d_i$ iid uniformly from $\{1,\ldots,100\}$.
		
	\item \TSTDB{1}: For all $i\in[n]$, we choose $C_i$ iid uniformly from $\{1,\ldots,100\}$, $\underline{c}_i$ from $\{1,\ldots,10\}$ and $d_i$ from $\{100-c_i,\ldots,100\}$.
		
	\item \TSTDB{2}: For all $i\in[n]$, we choose $C_i$ iid uniformly from $\{1,\ldots,100\}$, set $c_i = 100 - C_i$ and choose $d_i$  iid uniformly from $\{\underline{c}_i,\ldots,100\}$.
\end{itemize}

\subsection{Two-Stage Problems with Continuous Budgeted Uncertainty}
\begin{itemize}
	\item \TSTCB{U}: For all $i\in[n]$, we choose $C_i, \underline{c}_i, d_i$ iid uniformly from $\{1,\ldots,100\}$ (=\TSTDB{U}).
		
	\item \TSTCB{1}:  For all $i\in[n]$, we choose $C_i$ iid uniformly from $\{1,\ldots,100\}$, $\underline{c}_i$ from $\{1,\ldots,10\}$ and $d_i$ from $\{100-c_i,\ldots,100\}$ (=\TSTDB{1}).
		
	\item \TSTCB{2}:  For all $i\in[n]$, we choose $C_i$ iid uniformly from $\{1,\ldots,100\}$, set $c_i = 100 - C_i$ and choose $d_i$  iid uniformly from $\{\underline{c}_i,\ldots,100\}$ (=\TSTDB{2}).
\end{itemize}

\subsection{Recoverable Robust Problems with Discrete Uncertainty}
	
\begin{itemize}
	\item \RRD{U}: For all $i\in[n]$, we choose $C_i , c_i^j$ iid uniformly from $\{1,\ldots,100\}$ (=\TSTD{U}).
		
	\item \RRD{1}: For all $i\in[n]$, with probability equal to $0.5$, we choose $C_i$ from $\{45,\ldots,55\}$ and with probability equal to $0.5$, from $\{25,\ldots,75\}$. 
	For all $i\in[n]$ and $j\in[N]$, with probability equal to $0.5$ we choose $c_i^j$ from $\{C_i - 5,\ldots,C_i + 5\}$, with probability equal to $0.25$ we choose $c_i^j$ from $\{1,\ldots,10\}$ and with probability equal to $0.25$, $c_i^j$ from $\{91,\ldots,100\}$ (=\TSTD{1}).
		
	\item \RRD{2}: For all $i\in[n]$, with probability equal to $0.5$ we choose $C_i$ from $\{1,\ldots,100\}$ and with probability equal to $0.5$ we set $C_i = 50$. Then, for all $i\in[n]$ and $j\in[N]$, if $C_i = 50$ with probability equal to $0.5$ we choose $c_i^j$ from $\{1,\ldots,10\}$ and with probability equal to $0.5$ we choose $c_i^j$ from $\{91,\ldots,100\}$; else, we choose $c_i^j$ from $\{C_i - 5,\ldots,C_i + 5\}$ and set negative values to be equal to zero (=\TSTD{2}).
\end{itemize}

\subsection{Recoverable Robust Problems with Discrete Budgeted Uncertainty}
\begin{itemize}
	\item \RRDB{U}: For all $i\in[n]$, we choose $C_i, \underline{c}_i, d_i$ iid uniformly from $\{1,\ldots,100\}$ (=\TSTDB{U}).
		
	\item \RRDB{1} For all $i\in[n]$, we choose $C_i$ iid uniformly from $\{1,\ldots,100\}$, $\underline{c}_i$ from $\{1,\ldots,10\}$ and $d_i$ from $\{100-c_i,\ldots,100\}$ (=\TSTDB{1}).
		
	\item \RRDB{2}: For all $i\in[n]$, we choose $C_i$ iid uniformly from $\{1,\ldots,100\}$, set $c_i = 100 - C_i$ and choose $d_i$  iid uniformly from $\{\underline{c}_i,\ldots,100\}$ (=\TSTDB{2}).
\end{itemize}

\subsection{Recoverable Robust Problems with Continuous Budgeted Uncertainty}
\begin{itemize}
	\item \RRCB{U}: For all $i\in[n]$, we choose $C_i, \underline{c}_i, d_i$ iid uniformly from $\{1,\ldots,100\}$ (=\TSTDB{U}).
		
	\item \RRCB{1}: For all $i\in[n]$, we choose $C_i$ iid uniformly from $\{1,\ldots,100\}$, $\underline{c}_i$ from $\{1,\ldots,10\}$ and $d_i$ from $\{100-c_i,\ldots,100\}$ (=\TSTDB{1}).
		
	\item \RRCB{2}: For all $i\in[n]$, we choose $C_i$ iid uniformly from $\{1,\ldots,100\}$, set $c_i = 100 - C_i$ and choose $d_i$  iid uniformly from $\{\underline{c}_i,\ldots,100\}$ (=\TSTDB{2}).
\end{itemize}

\section{Generation Times}
\label{app:times}
	
In addition to solution times of instances, the generation time of these instances is of importance, too. Preferably, there should be a balance between both solution and generation time of instances, as spending a lot of time on generating instances that might only slightly increase the solution time in comparison to the corresponding sampling method is not acceptable.
	
In this sense, it should be noticed that some HIRO methods hit the 600-second time limit. In the following tables we illustrate the generation time of our HIRO methods for all cited criteria in this paper. It also can be easily observed that the generation time of all sampling methods can be ignored as they are considerably fast. 
	
\begin{table}[]
	\begin{center}
		\begin{tabular}{r|r|r|r|r|r|r}
			$n$ & $p$ & $N$ & $b$ & \texttt{UH} & \texttt{1H} & \texttt{2H} \\
			\hline
			20 & 10 & 20 & 1 &   1.21 &  37.71 & 586.87 \\
			   &    &    & 2 &  27.17 & 346.83 & 599.17 \\
			   &    &    & 5 & 596.35 & 508.87 & 599.22 \\
			20 & 11 & 20 & 1 &   1.15 &  16.93 & 195.01 \\
			   &    &    & 2 &  27.61 & 242.24 & 472.15 \\
			   &    &    & 5 & 591.55 & 409.46 & 599.85 \\
			25 & 13 & 25 & 1 &  10.19 & 125.55 & 553.69 \\
			   &    &    & 2 & 293.57 & 390.39 & 599.84 \\
			   &    &    & 5 & 599.91 & 479.44 & 599.78 \\
			30 & 15 & 30 & 1 &  62.94 & 329.47 & 573.91 \\
			   &    &    & 2 & 570.01 & 488.36 & 599.87 \\
			   &    &    & 5 & 599.92 & 534.38 & 599.89 \\
			35 & 17 & 35 & 1 & 288.29 & 342.12 & 599.29 \\
			   &    &    & 2 & 599.91 & 442.22 & 599.52 \\
			   &    &    & 5 & 599.99 & 556.39 & 599.71 \\
			40 & 20 & 40 & 1 & 553.75 & 469.73 & 599.92 \\
			   &    &    & 2 & 599.98 & 515.58 & 599.86 \\
			   &    &    & 5 & 600.01 & 589.19 & 599.44 \\
			40 & 21 & 40 & 1 & 530.51 & 441.24 & 596.07 \\
			   &    &    & 2 & 599.99 & 534.46 & 597.15 \\
			   &    &    & 5 & 600.01 & 598.27 & 596.71 
		\end{tabular}
		\caption{Min-Max problem under discrete uncertainty set - Exp1}
	\end{center}
\end{table}
	
\begin{table}[]
	\begin{center}
		\begin{tabular}{r|r|r|r|r}
			$p$ & $b$ & \texttt{UH} & \texttt{1H} & \texttt{2H} \\
			\hline
			5  & 1 &   0.81 &   7.47 & 189.19 \\
			   & 2 &   5.75 & 106.33 & 457.56 \\
			   & 5 & 436.99 & 234.42 & 599.99 \\
			10 & 1 &   8.19 & 202.81 & 599.94 \\
			   & 2 & 212.83 & 445.85 & 599.97 \\
			   & 5 & 600.03 & 469.48 & 599.92 \\
			11 & 1 &  21.89 & 216.37 & 560.04 \\
			   & 2 & 385.08 & 383.39 & 599.96 \\
			   & 5 & 600.01 & 436.44 & 599.92 \\
			15 & 1 &  62.94 & 329.47 & 573.91 \\
			   & 2 & 570.01 & 488.36 & 599.87 \\
			   & 5 & 599.92 & 534.38 & 599.89 \\
			20 & 1 &  12.33 & 119.56 & 599.92 \\
			   & 2 & 275.58 & 372.14 & 599.92 \\
			   & 5 & 599.98 & 555.09 & 599.83 \\
			21 & 1 &   8.36 &  57.71 & 533.68 \\
			   & 2 & 161.51 & 317.91 & 588.92 \\
			   & 5 & 573.73 & 510.98 & 599.83 \\
			25 & 1 &   0.99 &   6.29 & 297.04 \\
			   & 2 &   8.44 &  74.96 & 529.49 \\
			   & 5 & 240.71 & 434.21 & 599.89 
		\end{tabular}
		\caption{Min-Max problem under discrete uncertainty set - Exp2 ($n=N=30$)}
	\end{center}
\end{table}
	
\begin{table}[]
	\begin{center}
		\begin{tabular}{r|r|r|r|r}
			$N$ & $b$ & \texttt{UH} & \texttt{1H} & \texttt{2H} \\
			\hline
			5  & 1 &   0.41 &   7.86 &  16.26 \\
			   & 2 &  34.03 & 210.78 & 146.79 \\
			   & 5 & 551.48 & 430.74 & 572.03 \\
			10 & 1 &   2.43 &  89.53 & 148.61 \\
			   & 2 & 149.27 & 339.52 & 450.97 \\
			   & 5 & 599.88 & 505.95 & 599.77 \\
			15 & 1 &   7.08 & 181.06 & 259.81 \\
			   & 2 & 335.47 & 398.85 & 548.47 \\
			   & 5 & 599.91 & 461.67 & 599.67 \\
			20 & 1 &  17.81 & 170.39 & 438.02 \\
			   & 2 & 439.65 & 405.91 & 587.37 \\
			   & 5 & 599.88 & 449.26 & 599.69 \\
			25 & 1 &  47.27 & 202.91 & 528.54 \\
			   & 2 & 557.81 & 419.55 & 597.05 \\
			   & 5 & 599.91 & 512.57 & 599.79 \\
			30 & 1 &  62.94 & 329.47 & 573.91 \\
			   & 2 & 570.01 & 488.36 & 599.87 \\
			   & 5 & 599.92 & 534.38 & 599.89 \\
			35 & 1 & 142.93 & 294.53 & 593.41 \\
			   & 2 & 584.25 & 489.83 & 599.89 \\
			   & 5 & 599.96 & 543.21 & 599.91 \\
			40 & 1 & 149.74 & 329.88 & 599.87 \\
			   & 2 & 596.37 & 493.11 & 599.89 \\
			   & 5 & 599.96 & 507.77 & 599.91 
		\end{tabular}
		\caption{Min-Max problem under discrete uncertainty set - Exp3 ($n=30,p=15$)}
	\end{center}
\end{table}
	
\begin{table}[]
	\begin{center}
		\begin{tabular}{r|r|r|r|r|r|r|r|r|r|r}
			\multicolumn{2}{}{} &\multicolumn{3}{|c|}{HIRO-$\pmb{c}$} &\multicolumn{3}{|c|}{HIRO-$\pmb{d}$} &\multicolumn{3}{|c}{HIRO-$\pmb{c}$-$\pmb{d}$} \\
			\cline{3-11}
			$\Gamma$ & $b$ & \texttt{UH} & \texttt{1H} & \texttt{2H} & \texttt{UH} & \texttt{1H} & \texttt{2H} & \texttt{UH} & \texttt{1H} & \texttt{2H} \\
			\hline
			5  &  1 & 0.017 & 0.012 & 0.017 & 0.013 & 0.008 & 0.006 & 0.014 & 0.009 & 0.008 \\
			   &  2 & 0.015 & 0.012 & 0.021 & 0.014 & 0.009 & 0.007 & 0.015 & 0.009 & 0.008 \\
			   &  5 & 0.016 & 0.012 & 0.054 & 0.015 & 0.009 & 0.008 & 0.019 & 0.009 & 0.009 \\
			   & 10 & 0.019 & 0.014 & 0.046 & 0.019 & 0.009 & 0.008 & 0.031 & 0.010 & 0.012 \\
			   & 20 & 0.027 & 0.021 & 0.045 & 0.024 & 0.009 & 0.023 & 0.051 & 0.039 & 0.046 \\
			10 &  1 & 0.016 & 0.011 & 0.019 & 0.013 & 0.009 & 0.006 & 0.015 & 0.010 & 0.008 \\
			   &  2 & 0.016 & 0.012 & 0.027 & 0.014 & 0.009 & 0.006 & 0.016 & 0.010 & 0.008 \\
			   &  5 & 0.016 & 0.013 & 0.051 & 0.016 & 0.011 & 0.007 & 0.021 & 0.014 & 0.008 \\
			   & 10 & 0.021 & 0.017 & 0.038 & 0.023 & 0.012 & 0.007 & 0.041 & 0.042 & 0.009 \\
			   & 20 & 0.028 & 0.047 & 0.039 & 0.035 & 0.017 & 0.011 & 0.076 & 0.057 & 0.009 \\
			15 &  1 & 0.013 & 0.014 & 0.021 & 0.012 & 0.011 & 0.006 & 0.013 & 0.012 & 0.007 \\
			   &  2 & 0.013 & 0.013 & 0.027 & 0.013 & 0.011 & 0.006 & 0.012 & 0.014 & 0.008 \\
			   &  5 & 0.015 & 0.018 & 0.045 & 0.013 & 0.014 & 0.007 & 0.016 & 0.056 & 0.008 \\
			   & 10 & 0.015 & 0.041 & 0.041 & 0.018 & 0.035 & 0.008 & 0.026 & 0.081 & 0.008 \\
			   & 20 & 0.021 & 0.066 & 0.038 & 0.028 & 0.052 & 0.011 & 0.066 & 0.101 & 0.009 \\
			20 &  1 & 0.012 & 0.016 & 0.013 & 0.011 & 0.009 & 0.006 & 0.015 & 0.027 & 0.009 \\
			   &  2 & 0.012 & 0.015 & 0.017 & 0.009 & 0.011 & 0.006 & 0.015 & 0.031 & 0.011 \\
			   &  5 & 0.012 & 0.017 & 0.019 & 0.011 & 0.011 & 0.007 & 0.017 & 0.028 & 0.012 \\
			   & 10 & 0.014 & 0.019 & 0.018 & 0.011 & 0.015 & 0.008 & 0.024 & 0.052 & 0.011 \\
			   & 20 & 0.015 & 0.024 & 0.018 & 0.016 & 0.018 & 0.011 & 0.059 & 0.112 & 0.011 
		\end{tabular}
		\caption{Min-Max problem under budgeted uncertainty set}
	\end{center}
\end{table}

\begin{table}[]
	\begin{center}
		\begin{tabular}{r|r|r|r|r}
			$p$ & $b$ & \texttt{UH} & \texttt{1H} & \texttt{2H} \\
			\hline
			10 & 1 &  10.27 &   3.56 &   7.32 \\
			   & 2 &  13.34 &   9.23 &  18.91 \\
			   & 5 &  32.11 & 342.77 & 592.56 \\
			20 & 1 &   9.82 &   2.33 &   7.08 \\
			   & 2 &  12.27 &  18.59 &  18.01 \\
			   & 5 &  32.55 & 599.88 & 592.41 \\
			30 & 1 &  16.33 &   1.76 &   7.85 \\
			   & 2 &  19.36 &   3.97 &  18.43 \\
			   & 5 &  42.98 & 174.09 & 578.04 \\
			40 & 1 &  17.68 &   1.94 &   6.28 \\
			   & 2 &  21.51 &   4.29 &  21.54 \\
			   & 5 &  44.54 &   8.22 & 340.26 \\
			50 & 1 &  17.94 &   3.29 &   5.54 \\
			   & 2 &  24.54 &   8.23 &  16.93 \\
			   & 5 &  61.32 &  47.44 & 115.22 \\
			60 & 1 &  31.28 &   7.79 &   6.86 \\
			   & 2 &  48.32 &  57.67 &  25.38 \\
			   & 5 &  96.75 &  82.52 & 452.97 \\
			70 & 1 &  40.91 &   7.55 &   7.54 \\
			   & 2 &  59.51 &  34.48 &  18.66 \\
			   & 5 & 161.15 &  74.91 & 599.88 \\
			80 & 1 &  44.31 &   9.07 &   7.14 \\
			   & 2 &  67.74 &  24.51 &  21.46 \\
			   & 5 & 180.21 &  31.91 & 593.85 \\
			90 & 1 &  39.34 &   9.14 &   9.92 \\
			   & 2 &  59.58 &  19.11 &  41.36 \\
			   & 5 & 171.41 & 144.04 & 452.59
		\end{tabular}
		\caption{Min-Max Regret problem under interval uncertainty set ($n=100$)}
	\end{center}
\end{table}

\begin{table}[]
	\begin{center}
		\begin{tabular}{r|r|r|r|r|r|r}
			$n$ & $p$ & $N$ & $b$ & \texttt{UH} & \texttt{1H} & \texttt{2H} \\
			\hline
			30 & 15 & 30 & 1 & 573.87 & 519.23 & 572.47 \\		
			   &    &    & 2 & 599.89 & 512.64 & 599.95 \\
			   &    &    & 5 & 599.88 & 490.18 & 599.92 \\
			40 & 20 & 40 & 1 & 599.86 & 558.33 & 599.89 \\
			   &    &    & 2 & 599.88 & 585.66 & 599.91 \\
			   &    &    & 5 & 599.87 & 556.48 & 599.89 \\ 
			40 & 21 & 40 & 1 & 599.91 & 599.94 & 599.90 \\
			   &    &    & 2 & 599.91 & 599.95 & 599.93 \\
		   	   &    &    & 5 & 599.93 & 599.93 & 599.93
   	  	\end{tabular}
		\caption{Min-Max Regret problem under discrete uncertainty set - Exp1}
	\end{center}
\end{table}

\begin{table}[]
	\begin{center}
		\begin{tabular}{r|r|r|r|r}
			$p$ & $b$ & \texttt{UH} & \texttt{1H} & \texttt{2H} \\
			\hline
			10 & 1 & 557.76 & 275.91 & 567.07 \\
			   & 2 & 533.49 & 462.03 & 575.28 \\
			   & 5 & 573.87 & 519.23 & 572.47 \\
			11 & 1 & 576.31 & 483.14 & 563.74 \\
			   & 2 & 571.21 & 564.83 & 593.08 \\
			   & 5 & 599.84 & 456.07 & 599.85 \\
			15 & 1 & 599.89 & 574.49 & 599.85 \\
			   & 2 & 599.86 & 558.33 & 599.89 \\
			   & 5 & 599.92 & 564.22 & 599.86 \\
			20 & 1 & 599.91 & 579.45 & 599.88 \\
			   & 2 & 599.94 & 585.79 & 599.91 \\
			   & 5 & 599.92 & 599.98 & 599.93 \\
			21 & 1 & 599.91 & 599.94 & 599.91 \\
			   & 2 & 599.93 & 595.42 & 599.93 \\
			   & 5 & 599.93 & 599.95 & 599.91 
		\end{tabular}
		\caption{Min-Max Regret problem under discrete uncertainty set - Exp2 ($n=N=30$)}
	\end{center}
\end{table}

\begin{table}[]
	\begin{center}
		\begin{tabular}{r|r|r|r|r}
			$N$ & $b$ & \texttt{UH} & \texttt{1H} & \texttt{2H} \\
			\hline
			20 & 1 & 530.04 & 402.71 & 524.44 \\
			   & 2 & 573.87 & 519.23 & 572.47 \\
			   & 5 & 599.92 & 540.41 & 599.46 \\
			30 & 1 & 599.91 & 472.94 & 596.06 \\
			   & 2 & 599.86 & 558.33 & 599.89 \\
			   & 5 & 599.91 & 575.77 & 599.91 \\
			40 & 1 & 599.92 & 564.02 & 599.92 \\
			   & 2 & 599.91 & 599.94 & 599.91 \\
			   & 5 & 599.92 & 599.94 & 599.92 
		\end{tabular}
		\caption{Min-Max Regret problem under discrete uncertainty set - Exp3 ($n=30,p=15$)}
	\end{center}
\end{table}

\begin{table}[]
	\begin{center}
		\begin{tabular}{r|r|r|r|r|r|r|r|r|r}
			\multicolumn{4}{}{}&\multicolumn{3}{|c|}{HIRO-$\pmb{C}$}&\multicolumn{3}{|c}{HIRO-$\pmb{C}$-$\pmb{c}^j$} \\
			\cline{5-10}
			$n$ & $p$ & $N$ & $b$ & \texttt{UH} & \texttt{1H} & \texttt{2H} & \texttt{UH} & \texttt{1H} & \texttt{2H} \\
			\hline
			50  & 25 & 50  & 1 &  1.63 &   8.37 & 182.54 & 160.74 & 546.89 & 114.88 \\
			    &    &     & 2 &  2.91 &  17.15 & 241.46 & 339.47 & 599.99 & 157.22 \\
			    &    &     & 5 &  4.64 & 441.08 & 388.92 & 559.85 & 599.94 & 599.99 \\
			100 & 50 & 100 & 1 & 18.05 & 131.22 & 404.11 & 599.97 & 599.97 &  71.41 \\
				&    &     & 2 & 31.21 & 202.34 & 475.57 & 600.01 & 599.97 & 118.48 \\
				&    &     & 5 & 55.03 & 546.54 & 552.05 & 600.03 & 599.93 & 599.99 
		\end{tabular}
		\caption{Two-stage problem under discrete uncertainty set - Exp1}
	\end{center}
\end{table}

\begin{table}[]
	\begin{center}
		\begin{tabular}{r|r|r|r|r|r|r|r}
			\multicolumn{2}{}{} &\multicolumn{3}{|c|}{HIRO-$\pmb{C}$}&\multicolumn{3}{|c}{HIRO-$\pmb{C}$-$\pmb{c}^j$} \\
			\cline{3-8}
			$p$ & $b$ & \texttt{UH} & \texttt{1H} & \texttt{2H} & \texttt{UH} & \texttt{1H} & \texttt{2H} \\
			\hline
			10 & 1 &  0.65 &   0.85 &  96.25 &  30.08 &  69.27 &   5.36 \\
			   & 2 &  0.96 &   1.53 & 108.35 &  97.83 & 212.83 &  32.34 \\
			   & 5 &  2.21 &   5.12 & 108.56 & 351.17 & 526.31 & 352.21 \\
			20 & 1 &  1.09 &   5.12 & 152.27 &  74.85 & 497.68 &  67.16 \\
			   & 2 &  1.54 &   9.23 & 175.87 & 222.41 & 599.99 & 135.35 \\
			   & 5 & 15.32 & 198.44 & 241.21 & 547.27 & 599.99 & 574.21 \\
			25 & 1 &  1.63 &   8.37 & 182.54 & 160.74 & 546.89 & 114.88 \\
			   & 2 &  2.91 &  17.15 & 241.46 & 339.47 & 599.99 & 157.22 \\
			   & 5 &  4.64 & 441.08 & 388.92 & 559.85 & 599.94 & 599.99 \\
			30 & 1 &  2.15 &  19.52 & 241.31 & 226.81 & 587.84 & 124.45 \\
			   & 2 & 14.98 &  62.57 & 427.35 & 454.03 & 599.99 & 159.71 \\
			   & 5 & 10.56 & 517.46 & 590.31 & 579.28 & 599.99 & 599.99 \\
			40 & 1 &  5.01 & 422.42 & 464.64 & 356.72 & 599.99 & 119.85 \\
			   & 2 &  7.96 & 588.31 & 525.73 & 490.83 & 599.99 & 145.79 \\
			   & 5 & 37.23 & 600.01 & 593.88 & 570.59 & 599.99 & 585.47 
		\end{tabular}
		\caption{Two-stage problem under discrete uncertainty set - Exp2 ($n=N=50$)}
	\end{center}
\end{table}

\begin{table}[]
	\begin{center}
		\begin{tabular}{r|r|r|r|r|r|r|r}
			\multicolumn{2}{}{} &\multicolumn{3}{|c|}{HIRO-$\pmb{C}$}&\multicolumn{3}{|c}{HIRO-$\pmb{C}$-$\pmb{c}^j$} \\
			\cline{3-8}
			$p$ & $b$ & \texttt{UH} & \texttt{1H} & \texttt{2H} & \texttt{UH} & \texttt{1H} & \texttt{2H} \\
			\hline
			10 & 1 & 0.18 &   1.13 & 111.96 &   7.86 & 152.35 &  67.71 \\
			   & 2 & 0.27 &   3.21 & 167.35 &  22.12 & 485.25 & 168.04 \\
			   & 5 & 0.68 & 283.64 & 293.59 & 243.51 & 587.99 & 568.67 \\
			20 & 1 & 0.41 &   2.43 & 132.94 &  44.54 & 438.42 & 114.21 \\
			   & 2 & 0.63 &   6.08 & 192.45 & 151.48 & 577.16 & 168.15 \\
			   & 5 & 1.66 & 312.48 & 323.89 & 459.03 & 599.99 & 599.99 \\
			30 & 1 & 0.68 &   3.68 & 167.11 &  98.25 & 484.89 & 111.01 \\ 
			   & 2 & 1.08 &  19.56 & 188.54 & 240.41 & 599.41 & 168.47 \\ 
			   & 5 & 2.21 & 364.68 & 378.91 & 499.55 & 599.98 & 599.99 \\
			40 & 1 & 0.75 &   5.51 & 174.31 & 124.59 & 532.34 & 120.91 \\
			   & 2 & 1.23 &   9.99 & 226.56 & 274.29 & 599.96 & 156.78 \\
			   & 5 & 2.49 & 382.72 & 406.91 & 531.51 & 599.96 & 599.99 \\
			50 & 1 & 1.63 &   8.37 & 182.54 & 160.74 & 546.89 & 114.88 \\
			   & 2 & 2.91 &  17.15 & 241.46 & 339.47 & 599.99 & 157.22 \\
			   & 5 & 4.64 & 441.08 & 388.92 & 559.85 & 599.94 & 599.99 \\
			60 & 1 & 1.46 &  10.99 & 221.81 & 175.53 & 588.94 & 122.13 \\
			   & 2 & 2.25 &  22.88 & 287.59 & 373.31 & 599.97 & 157.84 \\
			   & 5 & 4.48 & 412.46 & 409.51 & 552.21 & 599.97 & 599.99 
		\end{tabular}
		\caption{Two-stage problem under discrete uncertainty set - Exp3 ($n=50,p=25$)}
	\end{center}
\end{table}

\begin{table}[]
	\begin{center}
		\begin{tabular}{r|r|r|r|r|r|r|r|r|r|r}
			\multicolumn{5}{}{} &\multicolumn{3}{|c|}{HIRO-$\pmb{C}$}&\multicolumn{3}{|c}{HIRO-$\pmb{C}$-$\pmb{c}^j$} \\
			\cline{6-11}
			$n$ & $p$ & $N$ & $\Delta$ & $b$ & \texttt{UH} & \texttt{1H} & \texttt{2H} & \texttt{UH} & \texttt{1H} & \texttt{2H} \\
			\hline
			50  & 25 & 50  & 13 & 1 &   2.04 &  47.91 & 464.73 & 142.56 & 529.29 & 120.15 \\
			    &    &     &    & 2 &   8.81 & 216.72 & 517.62 & 341.47 & 599.98 & 156.14 \\
			    &    &     &    & 5 &  57.94 & 599.15 & 581.39 & 563.49 & 599.97 & 599.98 \\
			50  & 25 & 50  & 20 & 1 &   1.97 &  42.73 & 439.31 &  92.38 & 512.24 & 120.09 \\
			    &    &     &    & 2 &   6.51 & 216.47 & 501.02 & 256.51 & 599.97 & 156.08 \\
			    &    &     &    & 5 &  65.99 & 593.15 & 574.47 & 559.54 & 599.96 & 599.98 \\
			100 & 50 & 100 & 25 & 1 &  29.11 & 498.21 & 599.99 & 539.75 & 599.95 &  13.12 \\
			    &    &     &    & 2 & 109.04 & 599.91 & 599.98 & 599.99 & 599.97 &  61.04 \\
			    &    &     &    & 5 & 541.29 & 599.83 & 599.98 & 599.99 & 599.96 & 599.95 \\
			100 & 50 & 100 & 40 & 1 &  30.51 & 516.23 & 599.98 & 497.07 & 599.94 &  12.69 \\
			    &    &     &    & 2 & 122.73 & 599.93 & 599.99 & 591.01 & 599.98 &  60.64 \\
			    &    &     &    & 5 & 525.51 & 599.81 & 599.98 & 599.99 & 599.98 & 599.95 
		\end{tabular}
		\caption{Recoverable problem under discrete uncertainty set - Exp1}
	\end{center}
\end{table}

\begin{table}[]
	\begin{center}
		\begin{tabular}{r|r|r|r|r|r|r|r|r}
			\multicolumn{3}{}{} &\multicolumn{3}{|c|}{HIRO-$\pmb{C}$}&\multicolumn{3}{|c}{HIRO-$\pmb{C}$-$\pmb{c}^j$} \\
			\cline{4-9}
			$p$ & $\Delta$ & $b$ & \texttt{UH} & \texttt{1H} & \texttt{2H} & \texttt{UH} & \texttt{1H} & \texttt{2H} \\
			\hline
			25 & 13 & 1 &  2.04 &  47.91 & 464.73 & 142.56 & 529.29 & 120.15 \\
			   &    & 2 &  8.81 & 216.72 & 517.62 & 341.47 & 599.98 & 156.14 \\
			   &    & 5 & 57.94 & 599.15 & 581.39 & 563.49 & 599.97 & 599.98 \\
			   & 20 & 1 &  1.97 &  42.73 & 439.31 &  92.38 & 512.24 & 120.09 \\
			   &    & 2 &  6.51 & 216.47 & 501.02 & 256.51 & 599.97 & 156.08 \\
			   &    & 5 & 65.99 & 593.15 & 574.47 & 559.54 & 599.96 & 599.98 \\
			30 & 15 & 1 &  1.46 &  15.44 & 351.81 &  37.23 & 254.71 & 120.13 \\
			   &    & 2 &  2.75 & 271.54 & 441.22 & 141.73 & 570.93 & 156.12 \\
			   &    & 5 & 45.39 & 596.62 & 548.04 & 328.28 & 599.96 & 599.98 \\
			   & 25 & 1 &  1.22 &  13.07 & 339.32 &  27.18 & 209.19 & 120.08 \\
			   &    & 2 &  2.24 & 285.11 & 425.85 & 100.48 & 587.53 & 156.08 \\
			   &    & 5 & 42.16 & 588.54 & 531.46 & 319.18 & 599.96 & 599.98 \\
			40 & 20 & 1 &  0.85 &   1.24 & 158.35 &  20.26 & 46.921 &  36.35 \\
			   &    & 2 &  1.12 &   3.21 & 182.96 &  30.67 & 246.08 &  36.95 \\
			   &    & 5 &  2.38 & 362.68 & 209.19 & 163.83 & 546.48 & 261.71 \\
			   & 30 & 1 &  0.76 &   1.12 & 156.33 &  15.71 &  41.19 &  36.33 \\
			   &    & 2 &  0.96 &   2.75 & 184.68 &  29.24 & 219.76 &  36.91 \\
			   &    & 5 &  2.11 & 366.69 & 207.55 & 150.08 & 548.66 & 261.34 
		\end{tabular}
		\caption{Recoverable problem under discrete uncertainty set - Exp2 ($n=N=50$)}
	\end{center}
\end{table}

\begin{table}[]
	\begin{center}
		\begin{tabular}{r|r|r|r|r|r|r|r|r}
			\multicolumn{3}{}{} &\multicolumn{3}{|c|}{HIRO-$\pmb{C}$}&\multicolumn{3}{|c}{HIRO-$\pmb{C}$-$\pmb{c}^j$} \\
			\cline{4-9}
			$N$ & $\Delta$ & $b$ & \texttt{UH} & \texttt{1H} & \texttt{2H} & \texttt{UH} & \texttt{1H} & \texttt{2H} \\
			\hline
			40 & 13 & 1 &  2.06 &  38.92 & 413.02 &  87.71 & 505.09 & 120.11 \\
			   & 13 & 2 &  8.01 & 168.72 & 491.43 & 309.91 & 599.96 & 156.11 \\
			   & 13 & 5 & 45.21 & 577.48 & 577.88 & 556.12 & 599.97 & 599.98 \\
			   & 20 & 1 &  1.96 &  38.88 & 413.02 &  66.57 & 495.11 & 120.11 \\
			   & 20 & 2 &  6.58 & 196.81 & 491.43 & 230.15 & 599.96 & 156.11 \\
			   & 20 & 5 & 49.69 & 569.71 & 577.88 & 545.59 & 599.97 & 599.98 \\
			50 & 13 & 1 &  2.04 &  47.91 & 464.73 & 142.56 & 529.29 & 120.15 \\
			   & 13 & 2 &  8.81 & 216.72 & 517.62 & 341.47 & 599.98 & 156.14 \\
			   & 13 & 5 & 57.94 & 599.15 & 581.39 & 563.49 & 599.97 & 599.98 \\
			   & 20 & 1 &  1.97 &  42.73 & 439.31 &  92.38 & 512.24 & 120.09 \\
			   & 20 & 2 &  6.51 & 216.47 & 501.02 & 256.51 & 599.97 & 156.08 \\
			   & 20 & 5 & 65.99 & 593.15 & 574.47 & 559.54 & 599.96 & 599.98 \\
			60 & 13 & 1 &  2.81 &  48.62 & 496.17 & 185.21 & 555.02 & 120.18 \\
			   & 13 & 2 &  9.91 & 231.19 & 531.54 & 375.75 & 599.96 & 156.16 \\
			   & 13 & 5 & 82.23 & 599.98 & 585.02 & 569.62 & 599.96 & 599.98 \\
			   & 20 & 1 &  2.71 &  47.78 & 466.45 & 133.82 & 530.02 & 120.12 \\
			   & 20 & 2 &  9.89 & 236.18 & 513.78 & 299.86 & 599.96 & 156.11 \\
			   & 20 & 5 & 70.25 & 599.32 & 579.41 & 568.27 & 599.96 & 599.98 
		\end{tabular}
		\caption{Recoverable problem under discrete uncertainty set - Exp3 ($n=50,p=25$)}
	\end{center}
\end{table}


\begin{thebibliography}{BTEGN09}

\bibitem[ABV09]{aissi2009min}
Hassene Aissi, Cristina Bazgan, and Daniel Vanderpooten.
\newblock Min--max and min--max regret versions of combinatorial optimization
  problems: A survey.
\newblock {\em European journal of operational research}, 197(2):427--438,
  2009.

\bibitem[Ave01]{averbakh2001complexity}
Igor Averbakh.
\newblock On the complexity of a class of combinatorial optimization problems
  with uncertainty.
\newblock {\em Mathematical Programming}, 90(2):263--272, 2001.

\bibitem[BG21]{bold2021recoverable}
Matthew Bold and Marc Goerigk.
\newblock Recoverable robust single machine scheduling with interval
  uncertainty.
\newblock {\em arXiv preprint arXiv:2107.09310}, 2021.

\bibitem[BGK18]{bertsimas2018data}
Dimitris Bertsimas, Vishal Gupta, and Nathan Kallus.
\newblock Data-driven robust optimization.
\newblock {\em Mathematical Programming}, 167(2):235--292, 2018.

\bibitem[BGKK19]{busing2019formulations}
Christina B{\"u}sing, Sebastian Goderbauer, Arie~MCA Koster, and Manuel
  Kutschka.
\newblock Formulations and algorithms for the recoverable $\gamma$-robust
  knapsack problem.
\newblock {\em EURO Journal on Computational Optimization}, 7(1):15--45, 2019.

\bibitem[BKK11]{busing2011recoverable}
Christina B{\"u}sing, Arie~MCA Koster, and Manuel Kutschka.
\newblock Recoverable robust knapsacks: the discrete scenario case.
\newblock {\em Optimization Letters}, 5(3):379--392, 2011.

\bibitem[BL11]{birge2011introduction}
John~R Birge and Francois Louveaux.
\newblock {\em Introduction to stochastic programming}.
\newblock Springer Science \& Business Media, 2011.

\bibitem[BS03]{bertsimas2003robust}
Dimitris Bertsimas and Melvyn Sim.
\newblock Robust discrete optimization and network flows.
\newblock {\em Mathematical programming}, 98(1):49--71, 2003.

\bibitem[BSS19]{bertsimas2019two}
Dimitris Bertsimas, Shimrit Shtern, and Bradley Sturt.
\newblock Two-stage sample robust optimization.
\newblock {\em arXiv preprint arXiv:1907.07142}, 2019.

\bibitem[BTEGN09]{ben2009robust}
Aharon Ben-Tal, Laurent El~Ghaoui, and Arkadi Nemirovski.
\newblock {\em Robust optimization}.
\newblock Princeton university press, 2009.

\bibitem[BTN98]{ben1998robust}
Aharon Ben-Tal and Arkadi Nemirovski.
\newblock Robust convex optimization.
\newblock {\em Mathematics of operations research}, 23(4):769--805, 1998.

\bibitem[B{\"u}s11]{busing2011phd}
Christina B{\"u}sing.
\newblock {\em Recoverable robustness in combinatorial optimization}.
\newblock Cuvillier Verlag, 2011.

\bibitem[CG16]{chassein2016recoverable}
Andr{\'e} Chassein and Marc Goerigk.
\newblock On the recoverable robust traveling salesman problem.
\newblock {\em Optimization Letters}, 10(7):1479--1492, 2016.

\bibitem[CG17]{chassein2017minmax}
Andr{\'e} Chassein and Marc Goerigk.
\newblock Minmax regret combinatorial optimization problems with ellipsoidal
  uncertainty sets.
\newblock {\em European Journal of Operational Research}, 258(1):58--69, 2017.

\bibitem[CG18]{chassein2018scenario}
Andr{\'e} Chassein and Marc Goerigk.
\newblock On scenario aggregation to approximate robust combinatorial
  optimization problems.
\newblock {\em Optimization Letters}, 12(7):1523--1533, 2018.

\bibitem[CGKZ18]{chassein2018recoverable}
Andr{\'e} Chassein, Marc Goerigk, Adam Kasperski, and Pawe{\l} Zieli{\'n}ski.
\newblock On recoverable and two-stage robust selection problems with budgeted
  uncertainty.
\newblock {\em European Journal of Operational Research}, 265(2):423--436,
  2018.

\bibitem[CGKZ20]{chassein2020approximating}
Andr{\'e} Chassein, Marc Goerigk, Adam Kasperski, and Pawe{\l} Zieli{\'n}ski.
\newblock Approximating combinatorial optimization problems with the ordered
  weighted averaging criterion.
\newblock {\em European Journal of Operational Research}, 286(3):828--838,
  2020.

\bibitem[DGR20]{dokka2020mixed}
Trivikram Dokka, Marc Goerigk, and Rahul Roy.
\newblock Mixed uncertainty sets for robust combinatorial optimization.
\newblock {\em Optimization Letters}, 14(6):1323--1337, 2020.

\bibitem[DK12]{dolgui2012min}
Alexandre Dolgui and Sergey Kovalev.
\newblock Min--max and min--max (relative) regret approaches to representatives
  selection problem.
\newblock {\em 4OR}, 10(2):181--192, 2012.

\bibitem[DW13]{deineko2013complexity}
Vladimir~G Deineko and Gerhard~J Woeginger.
\newblock Complexity and in-approximability of a selection problem in robust
  optimization.
\newblock {\em 4OR}, 11(3):249--252, 2013.

\bibitem[FHLW20]{fischer2020investigation}
Dennis Fischer, Tim~A Hartmann, Stefan Lendl, and Gerhard~J Woeginger.
\newblock An investigation of the recoverable robust assignment problem.
\newblock {\em arXiv preprint arXiv:2010.11456}, 2020.

\bibitem[FIMY15]{furini2015heuristic}
Fabio Furini, Manuel Iori, Silvano Martello, and Mutsunori Yagiura.
\newblock Heuristic and exact algorithms for the interval min--max regret
  knapsack problem.
\newblock {\em INFORMS Journal on Computing}, 27(2):392--405, 2015.

\bibitem[GH19]{goerigk2019representative}
Marc Goerigk and Martin Hughes.
\newblock Representative scenario construction and preprocessing for robust
  combinatorial optimization problems.
\newblock {\em Optimization Letters}, 13(6):1417--1431, 2019.

\bibitem[GKZ20]{goerigk2020combinatorial}
Marc Goerigk, Adam Kasperski, and Pawe{\l} Zieli{\'n}ski.
\newblock Combinatorial two-stage minmax regret problems under interval
  uncertainty.
\newblock {\em Annals of Operations Research}, pages 1--28, 2020.

\bibitem[GM20]{goerigk2020generating}
Marc Goerigk and Stephen~J Maher.
\newblock Generating hard instances for robust combinatorial optimization.
\newblock {\em European Journal of Operational Research}, 280(1):34--45, 2020.

\bibitem[GS16]{goerigk2016algorithm}
Marc Goerigk and Anita Sch{\"o}bel.
\newblock Algorithm engineering in robust optimization.
\newblock In {\em Algorithm engineering}, pages 245--279. Springer, 2016.

\bibitem[HDJPR21]{hashemi2021exploiting}
Hossein Hashemi~Doulabi, Patrick Jaillet, Gilles Pesant, and Louis-Martin
  Rousseau.
\newblock Exploiting the structure of two-stage robust optimization models with
  exponential scenarios.
\newblock {\em INFORMS Journal on Computing}, 33(1):143--162, 2021.

\bibitem[HRS18]{hansknecht2018fast}
Christoph Hansknecht, Alexander Richter, and Sebastian Stiller.
\newblock Fast robust shortest path computations.
\newblock In {\em 18th Workshop on Algorithmic Approaches for Transportation
  Modelling, Optimization, and Systems (ATMOS 2018)}. Schloss
  Dagstuhl-Leibniz-Zentrum fuer Informatik, 2018.

\bibitem[JWZG13]{jiang2013two}
Ruiwei Jiang, Jianhui Wang, Muhong Zhang, and Yongpei Guan.
\newblock Two-stage minimax regret robust unit commitment.
\newblock {\em IEEE Transactions on Power Systems}, 28(3):2271--2282, 2013.

\bibitem[KKZ15]{kasperski2015approximability}
Adam Kasperski, Adam Kurpisz, and Pawe{\l} Zieli{\'n}ski.
\newblock Approximability of the robust representatives selection problem.
\newblock {\em Operations Research Letters}, 43(1):16--19, 2015.

\bibitem[KMZ12]{kasperski2012tabu}
Adam Kasperski, Mariusz Makuchowski, and Pawe{\l} Zieli{\'n}ski.
\newblock A tabu search algorithm for the minmax regret minimum spanning tree
  problem with interval data.
\newblock {\em Journal of Heuristics}, 18(4):593--625, 2012.

\bibitem[KY95]{klir1995fuzzy}
George Klir and Bo~Yuan.
\newblock {\em Fuzzy sets and fuzzy logic}, volume~4.
\newblock Prentice hall New Jersey, 1995.

\bibitem[KY13]{kouvelis2013robust}
Panos Kouvelis and Gang Yu.
\newblock {\em Robust discrete optimization and its applications}, volume~14.
\newblock Springer Science \& Business Media, 2013.

\bibitem[KZ09]{kasperski2009randomized}
Adam Kasperski and Pawe{\l} Zieli{\'n}ski.
\newblock A randomized algorithm for the min-max selecting items problem with
  uncertain weights.
\newblock {\em Annals of Operations Research}, 172(1):221--230, 2009.

\bibitem[KZ16]{kasperski2016robust}
Adam Kasperski and Pawe{\l} Zieli{\'n}ski.
\newblock Robust discrete optimization under discrete and interval uncertainty:
  A survey.
\newblock In {\em Robustness analysis in decision aiding, optimization, and
  analytics}, pages 113--143. Springer, 2016.

\bibitem[KZ17]{kasperski2017robust}
Adam Kasperski and Pawe{\l} Zieli{\'n}ski.
\newblock Robust recoverable and two-stage selection problems.
\newblock {\em Discrete Applied Mathematics}, 233:52--64, 2017.

\bibitem[MBMG07]{montemanni2007robust}
Roberto Montemanni, J{\'a}nos Barta, Monaldo Mastrolilli, and Luca~Maria
  Gambardella.
\newblock The robust traveling salesman problem with interval data.
\newblock {\em Transportation Science}, 41(3):366--381, 2007.

\bibitem[MG04]{montemanni2004exact}
Roberto Montemanni and Luca~Maria Gambardella.
\newblock An exact algorithm for the robust shortest path problem with interval
  data.
\newblock {\em Computers \& Operations Research}, 31(10):1667--1680, 2004.

\bibitem[MHS01]{muller2001algorithm}
Matthias M{\"u}ller-Hannemann and Stefan Schirra.
\newblock {\em Algorithm Engineering}.
\newblock Springer, 2001.

\bibitem[MPS13]{monaci2013exact}
Michele Monaci, Ulrich Pferschy, and Paolo Serafini.
\newblock Exact solution of the robust knapsack problem.
\newblock {\em Computers \& operations research}, 40(11):2625--2631, 2013.

\bibitem[San09]{sanders2009algorithm}
Peter Sanders.
\newblock Algorithm engineering--an attempt at a definition.
\newblock In {\em Efficient algorithms}, pages 321--340. Springer, 2009.

\bibitem[SHY17]{shang2017data}
Chao Shang, Xiaolin Huang, and Fengqi You.
\newblock Data-driven robust optimization based on kernel learning.
\newblock {\em Computers \& Chemical Engineering}, 106:464--479, 2017.

\bibitem[SLTW12]{song2012incomplete}
Xiang Song, Rhyd Lewis, Jonathan Thompson, and Yue Wu.
\newblock An incomplete m-exchange algorithm for solving the large-scale
  multi-scenario knapsack problem.
\newblock {\em Computers \& operations research}, 39(9):1988--2000, 2012.

\bibitem[Sub21]{subramanyam2021lagrangian}
Anirudh Subramanyam.
\newblock A lagrangian dual method for two-stage robust optimization with
  binary uncertainties.
\newblock {\em arXiv preprint arXiv:2112.13138}, 2021.

\bibitem[TYK08]{taniguchi2008heuristic}
Fumiaki Taniguchi, Takeo Yamada, and Seiji Kataoka.
\newblock Heuristic and exact algorithms for the max--min optimization of the
  multi-scenario knapsack problem.
\newblock {\em Computers \& Operations Research}, 35(6):2034--2048, 2008.

\bibitem[WIMY18]{wu2018exact}
Wei Wu, Manuel Iori, Silvano Martello, and Mutsunori Yagiura.
\newblock Exact and heuristic algorithms for the interval min-max regret
  generalized assignment problem.
\newblock {\em Computers \& Industrial Engineering}, 125:98--110, 2018.

\bibitem[YK12]{yager2012ordered}
Ronald~R Yager and Janusz Kacprzyk.
\newblock {\em The ordered weighted averaging operators: theory and
  applications}.
\newblock Springer Science \& Business Media, 2012.

\bibitem[ZZ13]{zeng2013solving}
Bo~Zeng and Long Zhao.
\newblock Solving two-stage robust optimization problems using a
  column-and-constraint generation method.
\newblock {\em Operations Research Letters}, 41(5):457--461, 2013.

\end{thebibliography}
\end{document}